\RequirePackage{fix-cm}
\documentclass[envcountsect,envcountsame]{svjour3}
\usepackage{amsmath,amssymb,mathtools,algorithm,bm,bbm,xcolor}
\usepackage{latexsym,eqnarray}
\usepackage{booktabs}
\usepackage[colorlinks=true,citecolor=blue,linkcolor=blue,urlcolor=blue,bookmarks,bookmarksopen,bookmarksdepth=2,backref=page]{hyperref}

\usepackage{multirow}

\usepackage{bookmark}
\bookmarksetup{
	numbered, 
	open
}

\usepackage{threeparttable}

\usepackage{enumerate}
\usepackage{algorithm,algorithmic}
\floatname{algorithm}{Algorithm}

\DeclareMathOperator{\gjb}{GJB}
\DeclareMathOperator{\ind}{ind}
\DeclareMathOperator{\rind}{r-ind}

\smartqed

\newcommand*{\doi}[1]{\href{\detokenize{#1}}{doi: \detokenize{#1}}}

\renewcommand*{\backref}[1]{}
\renewcommand*{\backrefalt}[4]{
	\ifcase #1 (Not cited.)
	\or        (Cited on page~#2.)
	\else      (Cited on pages~#2.)
	\fi}

\newcommand{\fF}[1]{\mathbb{F}_{2^{#1}}}

\def\cM{{\mathcal{M}}}
\def\cC{{\mathcal{C}}}
\def\cD{{\mathcal{D}}}
\newcommand{\B}{\mathcal{B}}
\def \F {{\mathbb F}}

\newcommand{\FF}{{\mathbb F}}

\newcommand{\MrsClass}{\mathcal{M}_{r,s}^{\#}}

\newtheorem{cons}[theorem]{Construction}
\newtheorem{openproblem}{Open Problem}
\usepackage[T1]{fontenc}

\numberwithin{equation}{section}
\numberwithin{table}{section}

\begin{document}
\title{Bent Functions and the Completed Maiorana-McFarland Class}

\titlerunning{Bent Functions and the Completed Maiorana-McFarland Class}

\author{Enes Pasalic \and Alexandr Polujan \and Fengrong Zhang \and Sadmir Kudin}

\authorrunning{E. Pasalic, A. Polujan, F. Zhang, S. Kudin}

\institute{ 
            Enes Pasalic \at 
            University of Primorska, FAMNIT \& IAM, Glagolja\v{s}ka 8, 6000 Koper, Slovenia \\
            \email{enes.pasalic6@gmail.com}
            \and
            Alexandr Polujan \at
            Otto von Guericke University, Universit\"{a}tsplatz 2, 39106, Magdeburg, Germany \\
            \email{alexandr.polujan@ovgu.de}
            \and 
           Fengrong Zhang \at
           School of Cyber Engineering, Xidian University, Xi'an, 710071, P.R. China\\
              \email{zhfl203@163.com}
            \and
			Sadmir Kudin \at 
			University of Primorska, FAMNIT \& IAM, Glagolja\v{s}ka 8, 6000 Koper, Slovenia \\
			\email{sadmir.kudin@iam.upr.si}       
}

\date{Received: date / Accepted: date}

\maketitle

\begin{abstract}
In this article, we survey some fundamental results and recent advances in the design and analysis of Boolean bent functions from the perspective of their relation to the completed Maiorana-McFarland class.
\end{abstract}

\tableofcontents 

\section{Introduction}
Bent functions form a fundamental class of Boolean functions in an even number of variables that achieve the maximum possible distance from affine functions and therefore play an important role in cryptography, coding theory, and related areas of discrete mathematics. They were introduced by Rothaus in 1966 and published in 1976~\cite{Rothaus76_JCTA}. In the meantime, several fundamental works dedicated to their design and analysis appeared, and the results presented there largely shaped the development of the research area for decades. Among these works are the survey article of Dillon~\cite{Dillon72_SURVEY}, his PhD thesis~\cite{Dillon74_PhD}, and McFarland's construction of difference sets~\cite{McFarland73_JCTA}, which, in the binary elementary abelian setting, corresponds to bent functions of the form $f(x,y)=x\cdot\pi(y)+h(y)$, for all $x,y\in\F_2^m$, where $\pi$ is a permutation of $\F_2^m$ and $h$ is an arbitrary Boolean function on $\F_2^m$. The set of bent functions on $\F_2^n$ with $n=2m$ of this form is denoted by $\mathcal{M}$ and is called the \emph{Maiorana-McFarland class}, in honour of Maiorana and McFarland, who independently discovered this family. The closure of this set under the action of the affine general linear group $\operatorname{AGL}(n,2)$ on the input, together with the addition of affine functions on $\F_2^n$ to the output, is called the \emph{completed Maiorana-McFarland class} and is denoted by $\mathcal{M}^\#$. Since their introduction, Maiorana-McFarland bent functions have played a central role in the theory of bent functions and, more generally, in the study of cryptographic Boolean functions. They constitute the first and most thoroughly understood explicit class and have strongly influenced subsequent constructions, classifications, and the analysis of algebraic, combinatorial, and cryptographic properties~\cite{Mesnager16_BOOK,Carlet21_BOOK}.

In the past five decades, research on bent functions has developed in various directions, including $p$-ary bent functions and generalized bent functions~\cite{KumarSV1985,Nyberg91,Agievich_2008,Schmidt2009,HodzicMP18,AnbarM22,Meidl22_CCDS,AKM_Survey}, bent functions over arbitrary groups~\cite{CarletD04,Pott04_DAM,KoelschP24_COMB,Xu2025}, and connections with design and coding theory~\cite{Bending93_PhD,CarletCZ98_DCC,Ding:CodesfromDS,Ding:DesignsfromLC,Limesnager2020,DingMT19,DingT20,PolujanP21,MeidlPP23}. Several books~\cite{Carlet10_CHAPTER_BOOLEAN,Carlet10_CHAPTER_VECTORIAL,LogachevSY12_BOOK,Budaghyan14_BOOK,Tokareva15_BOOK,Mesnager16_BOOK,CusickS17_BOOK,Carlet21_BOOK} on bent functions and closely related classes of cryptographic functions have appeared, as well as numerous survey articles~\cite{DobbertinL04_survey,HellesethK13,KholoshaP13_HANDBOOK,Carlet14_OP,CarletM16_DCC,Meidl22_CCDS,AKM_Survey} devoted to various aspects of bent functions, including open problems, many of which remain largely unresolved.

While many research directions on bent functions exist, the main aim of this survey article is to provide a comprehensive overview of \emph{Boolean bent functions} from the perspective of their relation to the \emph{completed Maiorana-McFarland class}, with a particular focus on exclusion from $\mathcal{M}^\#$, recent developments in this area, and possible directions for future research.

The paper is organized as follows. After giving some preliminary definitions on Boolean and bent functions in Section~\ref{sec 1.1}, we summarize in Section~\ref{sec MM} several fundamental results on the completed Maiorana-McFarland class of Boolean bent functions. These include the characterization of functions in this class in terms of vanishing second-order derivatives (Lemma~\ref{lem M-M second}), as well as known computational results in eight variables~\cite{LangevinL11_DCC} indicating that most bent functions in the first dimension, in which they are not completely classified under (EA-)equivalence, do not belong to this class. This observation provides one of the main motivations for investigating the \emph{$\mathcal{M}^\#$ class exclusion problem}, that is, determining the intersection of a given class of bent functions with the completed Maiorana-McFarland class.

In Section~\ref{sec C and D classes}, we consider in detail recent advances in explicit construction methods and in the analysis of the $\mathcal{M}^\#$-exclusion problem for the $\mathcal{C}$ and $\mathcal{D}$ classes of bent functions (as well as their superclasses). These two classes were introduced by Carlet in~\cite{Carlet93_EUROCRYPT} by adding to Maiorana-McFarland functions of the form $x\cdot \pi(y)$ indicators of suitably chosen vector subspaces: $1_{L^{\perp}}(x)$ in the case of the $\mathcal{C}$ class, where $L$ is a vector subspace of $\F_2^m$ and $\pi$ is a permutation of $\F_2^m$ such that $\pi^{-1}(a+L)$ is a flat (i.e., an affine subspace) for every $a\in\F_2^m$ (this condition is called property~\eqref{eq: C property}); $1_{E_1}(x)1_{E_2}(y)$ in the case of the $\mathcal{D}$ class, where $\pi$ is a permutation of $\F_2^m$ and $E_1,E_2$ are vector subspaces of $\F_2^m$ such that $\pi(E_2)=E_1^{\perp}$; and $\delta_0(x)$ in the case of the subclass $\mathcal{D}_0\subset\mathcal{D}$, where $\delta_0$ denotes the indicator of the zero vector in $\F_2^m$. In Section~\ref{sec C class}, we first summarize sufficient conditions guaranteeing that pairs $(\pi,L)$ satisfy the property~\eqref{eq: C property} and hence yield bent functions in the $\mathcal{C}$ class, and then present known sufficient conditions ensuring that such functions lie outside $\mathcal{M}^\#$. Similar questions for the $\mathcal{D}_0$ class are discussed in Section~\ref{sec D class}, followed by a complete characterization of bent functions in $\mathcal{D}_0\cap\mathcal{M}^\#$ in Section~\ref{sec:D0class}. To conclude the section, we consider superclasses of the $\mathcal{C}$, $\mathcal{D}$, and $\mathcal{D}_0$ classes obtained by adding suitable linear combinations of the corresponding indicators, discuss their relation to $\mathcal{M}^\#$, and indicate applications to the design and analysis of further cryptographically significant classes of (vectorial) Boolean functions.

In Section~\ref{sec Further construcitons outside MM}, we summarize known results on classes of bent functions constructed using the multiplicative structure of finite fields and describe their relation to the completed Maiorana-McFarland class. These include monomial, binomial, and related bent functions (Section~\ref{sec monomial binomial etc}), bent functions arising as indicators of AB functions (Section~\ref{sec bent from AB and APN}), the partial spread class $\mathcal{PS}$, Dillon's $\mathcal{H}$ class and its variants (Section~\ref{subsec:classH}), as well as a recent approach based on the so-called $\mathcal{P}_\tau$ property (Section~\ref{sec Pi tau property}).

The fact that many of the families of bent functions summarized in Sections~\ref{sec C and D classes} and~\ref{sec Further construcitons outside MM} contain only a comparable number of functions to that of the $\mathcal{M}$ class, or even only a few examples in each fixed dimension, makes it difficult to explain theoretically the discrepancy observed in dimension $8$ between the total number of bent functions ($\approx 2^{106}$ obtained in~\cite{LangevinL11_DCC}) and those arising from known construction methods. This observation motivates further investigation of bent functions in the generalized Maiorana-McFarland class $\mathcal{GMM}$. These are bent functions in $n=r+s$ variables of the form $x\cdot\phi(y)+h(y)$, where $x\in\F_2^r$, $y\in\F_2^s$, and $\phi\colon\F_2^s\to\F_2^r$ and $h\colon\F_2^s\to\F_2$ are suitably chosen. 

In Section~\ref{sec GMM class}, we summarize recent advances in the design and analysis of bent functions in the generalized Maiorana-McFarland class $\mathcal{GMM}$. While the $\mathcal{GMM}$ class provides a broad framework for constructing bent functions, its algebraic structure is considerably less explicit than that of the classical Maiorana-McFarland class. After introducing some fundamental results and definitions related to the generalized completed Maiorana-McFarland class in Section~\ref{sec GMM basics} and surveying known characterizations of bentness in Section~\ref{sec Characterizing bentness in GMM}, we study two special subclasses of $\mathcal{GMM}$ in more detail. In Section~\ref{sec Almost MM}, we consider almost Maiorana-McFarland bent functions, which are, in some sense, as close as possible to classical Maiorana-McFarland functions, and provide a complete characterization showing that certain choices of building blocks may lead to functions both inside and outside $\mathcal{M}^\#$. In Section~\ref{sec Minimal linearity index}, we study $\ell$-optimal bent functions, which can be viewed as opposite to Maiorana-McFarland functions. Similarly, we show that certain known constructions yield such functions when the building blocks are properly chosen.

In Section~\ref{sec M-subspaces}, we investigate $\mathcal{M}$-subspaces of Boolean (bent) functions, which are special vector subspaces of $\F_2^n$ serving as witnesses of membership in the corresponding generalized Maiorana-McFarland class and are therefore fundamental algebraic objects for understanding the relation of bent functions to the generalized completed Maiorana-McFarland class. In Section~\ref{sec M-subspaces basics}, we summarize algebraic and combinatorial properties of these subspaces and discuss constructions of Maiorana-McFarland bent functions with many or few $\mathcal{M}$-subspaces. Bent functions with a unique $\mathcal{M}$-subspace of the maximal dimension are studied in detail in Section~\ref{se M-subspace unique}; such functions play an important role in understanding the cardinality of $\mathcal{M}^\#$ and in constructing bent functions outside this class. In Section~\ref{sec M-subspaces primary secondary}, we summarize recent results on $\mathcal{M}$-subspaces of secondary constructions of bent functions and discuss their implications for the theory of bent functions.

Finally, the paper is concluded in Section~\ref{sec Conclusion} with a list of further open problems, complementing those stated earlier, and a discussion of possible directions for future research.

\subsection{Preliminaries}\label{sec 1.1}
Let $\mathbb{F}_2^n$ be the vector space of all binary $n$-tuples $x=(x_1,\ldots,x_n)$, where $x_i \in \mathbb{F}_2$. For two elements $x=(x_1,\ldots,x_n)$ and $y=(y_1,\ldots,y_n)$ of $\mathbb{F}^n_2$, we define the scalar product over $\mathbb{F}_2$ as $x\cdot y=x_1 y_1 + \cdots +  x_n y_n$. Throughout this paper, $\operatorname{wt}(\cdot)$ denotes the \emph{Hamming weight}. For a vector in $\F_2^n$, it refers to the number of its nonzero coordinates, while for an integer it refers to the Hamming weight of its binary representation. We also denote by $0_n=(0,0,\ldots,0)\in \mathbb{F}^n_2$ the all-zero vector with $n$ coordinates. If necessary, we endow $\F_2^n$ with the structure of the finite field $\left(\F_{2^{n}},+,\cdot \right)$. An element $\alpha \in \mathbb{F}_{2^n}$ is called a \emph{primitive element}, if it is a generator of the multiplicative group $\mathbb{F}_{2^n}^*$. In this case, we define the scalar product on $\F_{2^n}$ by $x\cdot y={\rm Tr}(xy)$, where ${\rm Tr}(z)={\rm Tr}_1^n(z)$ is the \emph{absolute trace} and  ${\rm Tr}^n_m(z)=\sum_{i=0}^{\frac{n}{m}-1} z^{2^{i\cdot m}}$ is the \emph{relative trace} of $z \in \mathbb{F}_{2^n}$ from $\mathbb{F}_{2^n}$ into the subfield $\mathbb{F}_{2^m}$.

In this paper, we denote the set of all Boolean functions in $n$ variables by $\mathcal{B}_n$. One can uniquely represent any Boolean function $f\in\mathcal{B}_n$ using the \emph{algebraic normal form (ANF, for short)}, which is given by $f(x_1,\ldots,x_n)=\sum_{u\in \mathbb{F}^n_2}{\lambda_u}{(\prod_{i=1}^n{x_i}^{u_i})}$, where $x_i, \lambda_u \in \mathbb{F}_2$ and $u=(u_1, \ldots,u_n)\in \mathbb{F}^n_2$. The \emph{support} of $f \in \mathcal{B}_n$, denoted by $\textnormal{supp}(f)$, is defined as $\textnormal{supp}(f)=\{x \in \F_2^n\colon f(x)=1 \}.$
	The \emph{algebraic degree} of $f\in\mathcal{B}_n$, denoted by $\deg(f)$, is the maximum Hamming weight of $u \in \F_2^n$ for which $\lambda_u \neq 0$ in its ANF. By convention, the degree of the identically zero function is defined as $\deg(0)=-\infty$. A Boolean function $f\in\mathcal{B}_n$ is called \emph{homogeneous} if all the monomials in its ANF have the same algebraic degree.

    A mapping $F\colon\F_2^n\rightarrow\F_2^m$ is called a \emph{vectorial Boolean function}. The special case $m=1$ corresponds to the set $\mathcal{B}_n$ of Boolean functions. Every vectorial function $F\colon\F_2^n\rightarrow\F_2^m$ can be written as $F(x)=(f_1(x),\dots,f_m(x))$, for all $x\in\F_2^n$, where each Boolean function $f_i$ on $\F_2^n$ is called a \emph{coordinate function} of $F$. The \emph{components} of $F$ are the Boolean functions obtained as nonzero linear combinations of its coordinate functions, i.e., $b\cdot F(x)=b_1f_1(x)+\cdots+b_mf_m(x)$, for $b=(b_1,\ldots,b_m)\in\F_2^m\setminus\{0_m\}$. The \emph{algebraic normal form (ANF)} of a vectorial Boolean function is defined coordinate-wise. Consequently, the \emph{degree} of $F$ is the maximum degree among its coordinate functions, i.e., $\deg(F)=\max_{1\leq i\leq m}\deg(f_i)$. The \emph{first-order derivative} of $F$ in the direction $a\in\F_2^n$ is the mapping $D_aF(x)=F(x+a)+F(x)$. Derivatives of higher orders are defined recursively, i.e., the \emph{$k$-th order derivative} of $F$ is
$$
D_{a_k}D_{a_{k-1}}\cdots D_{a_1}F(x)
=
D_{a_k}(D_{a_{k-1}}\cdots D_{a_1}F)(x),
$$
for $a_1,\ldots,a_k\in\F_2^n$. A vector $a\in\F_2^n$ is called a \emph{linear structure} of $F$ if there exists a vector $v\in\F_2^m$ such that
$D_aF(x)=F(x+a)+F(x)=v$, for all $x\in\F_2^n$. We say that $F$ \emph{has no nonzero linear structures} if $0_n$ is its only linear structure. A mapping $F\colon\F_2^n\rightarrow\F_2^m$ is called \emph{$k$-to-$1$} if every element of its image has exactly $k$ preimages. In the particular case $n=m$ and $k=1$, such a mapping is called a \emph{permutation}.

Throughout this paper, we frequently identify $\F_2^n$ with the finite field $\F_{2^n}$ via a fixed $\F_2$-basis. Under this identification, any function $F\colon\F_2^n\rightarrow\F_2^n$ can be uniquely represented by a univariate polynomial $F(x)=\sum_{i=0}^{2^n-1}a_ix^i$, with coefficients $a_i\in\F_{2^n}$. In this representation, permutations correspond to permutation polynomials, and linear mappings correspond to linearized polynomials. When $m|n$, any function $F\colon\F_2^n\to\F_2^m$ can be written as a polynomial $F\colon\F_{2^n}\to\F_{2^m}$ given by $F(x)=\operatorname{Tr}^{n}_{m}\left(\sum_{i=0}^{2^{n}-1} a_{i} x^{i}\right)$. This representation is called the \emph{univariate (trace) representation}, however, it is not unique in general.

    The \emph{Walsh-Hadamard transform (WHT)} of a Boolean function $f\in\mathcal{B}_n$ is the mapping
$$W_f\colon\F_2^n\rightarrow\mathbb{Z},\qquad W_f(u)=\sum_{x\in\F_2^n}(-1)^{f(x)+u\cdot x}.$$
The multiset $\{*\, W_f(u)\colon u\in\F_2^n\, *\}$ is called the \emph{Walsh spectrum} of $f$. The Boolean function $f\in\mathcal{B}_n$ can be recovered from its Walsh-Hadamard transform via the \emph{inverse Walsh-Hadamard transform},
$$ (-1)^{f(x)} = 2^{-n}\sum_{u\in\F_2^n}W_f(u)(-1)^{u\cdot x}.$$
		
	\begin{definition}
	    For even $n$, a function $f\in\mathcal{B}_n$ is called \emph{bent} if $W_f(u)=\pm2^{\frac{n}{2}}$ for all $u\in\F_2^n$. 
	\end{definition}
    
    For a bent function $f\in\mathcal{B}_n$, a Boolean function $f^*\in \mathcal{B}_n$ defined by $W_f(u)=2^{\frac{n}{2}}(-1)^{f^*(u)}$ for all $u\in\F_2^n$ is a bent function, called the \emph{dual} of $f$. Another important related notion is that of a \emph{vectorial bent function}. A function $F\colon\F_2^n\to\F_2^m$, with even $n$, is called \emph{vectorial bent} if, for every nonzero $b\in\F_2^m$, the component function $F_b(x)=b\cdot F(x)$ is bent. It is well-known that vectorial bent functions exist only when $m\le n/2$; this result is known as the \emph{Nyberg bound}~\cite{Nyberg91}. For $m>n/2$, one can instead consider \emph{MNBC functions}, namely vectorial Boolean functions $F\colon\F_2^n\to\F_2^m$ having the \emph{maximum number of bent components}, equal to $2^m-2^{m-n/2}$~\cite{PottPMB18,ZhengPengDCC2020}. Thus, the notion of MNBC functions can naturally be seen as an extension of that of vectorial bent functions, which attain the maximal possible number $2^m-1$ of bent components when $m\le n/2$.

    Another two important classes of Boolean functions closely related to bent functions are \emph{semi-bent} and \emph{$5$-valued spectrum} functions, which naturally arise in concatenation constructions of bent functions. A function $f\in \mathcal{B}_n$ is called \emph{semi-bent} if its Walsh transform takes only the values $0$ and $\pm 2^{\frac{n+s}{2}}$, where $s=1$ if $n$ is odd, and $s=2$ if $n$ is even. A function $f\in \mathcal{B}_n$, where $n$ is even, is called a \emph{$5$-valued spectrum function} if its Walsh transform takes only the values $0$, $\pm 2^{\frac{n}{2}}$, and $\pm 2^{\frac{n+2}{2}}$. 

    One of the possible ways to construct new Boolean functions is by concatenating the known ones. Two types of concatenations will be particularly useful in our context. The first is the \emph{2-concatenation}, denoted by $h=f_1||f_2\in\mathcal{B}_{n+1}$ and defined by
$$h(x,y)=(1+y)f_1(x)+y f_2(x),$$
where $x\in\F_2^n$ and $y\in\F_2$. A Boolean function $h=f_1||f_2$ is bent if and only if $f_1$ and $f_2$ are semi-bent functions with disjoint Walsh supports~\cite{ZhengZhang_IEEE_2001}, i.e., $\operatorname{supp}(W_{f_1}) \cap \operatorname{supp}(W_{f_2}) = \emptyset$, where $\operatorname{supp}(W_{f_i})=\{u\in\F_2^n\colon W_{f_i}(u)\neq 0\}$. The second is the \emph{4-concatenation}, denoted by $h=f_1||f_2||f_3||f_4\in\mathcal{B}_{n+2}$ and defined by
\begin{equation}\label{eq:4-concat canonical}
    h(x,y,z)=(y+1)(z+1)f_1(x)+(y+1)zf_2(x)+y(z+1)f_3(x)+yzf_4(x),
\end{equation}
where $x\in\F_2^n$ and $y,z\in\F_2$, and the restrictions $h(x,0,0)$, $h(x,0,1)$, $h(x,1,0)$, and $h(x,1,1)$ correspond, respectively, to the functions $f_1$, $f_2$, $f_3$, and $f_4$. More generally, let $f\in\mathcal{B}_{n+2}$ be a Boolean bent function and let $S$ be an $n$-dimensional subspace of $\F_2^{n+2}$. If $f_1$ is the restriction of $f$ to $S$, and $f_2,f_3,f_4$ are the restrictions of $f$ to the remaining three cosets of $S$, then the quadruple $(f_1,f_2,f_3,f_4)_S$ is referred to as the \emph{4-decomposition} of $f$ with respect to $S$. By~\cite[Theorem~7]{CanteautC03_IEEE}, the functions in such a 4-decomposition are necessarily either all bent, all semi-bent, or all $5$-valued spectrum functions. The 4-concatenation defined by Eq.~\eqref{eq:4-concat canonical} is called \emph{canonical} and corresponds to the decomposition with respect to the subspace $S=\langle (0,\ldots,1,0),(0,\ldots,0,1)\rangle^\perp$, where $\langle u,v\rangle$ denotes the subspace spanned by the vectors $u$ and $v$, and $\langle u,v\rangle^\perp$ denotes the \emph{orthogonal complement} of $\langle u,v\rangle\subset\F_2^{n+2}$.

    Finally, we introduce the following natural equivalence relation on Boolean functions, which will be used in the next section to define the completed Maiorana-McFarland class. Boolean functions $f,f'\in\mathcal{B}_n$ are called \emph{(extended-affine) equivalent}, if there exists an affine permutation $A$ of $\F_2^n$ and an affine function $l\in\mathcal{B}_n$, such that $f\circ A + l= f'$. It is well-known that extended-affine equivalence preserves the bent property. A class of bent functions $\mathit{B}_n \subset \mathcal{B}_n$ is \emph{complete} if it is globally invariant under extended-affine equivalence. By $\mathit{B}_n^\#$ we denote the \emph{completed $\mathit{B}_n$ class}, which is the smallest possible complete class that  contains $\mathit{B}_n$.

\subsection{The completed Maiorana-McFarland class \texorpdfstring{$\mathcal{M}^\#$}{M\#}}\label{sec MM}

The \emph{Maiorana-McFarland class} $\cM$ is the set of $n$-variable ($n=2m$) Boolean bent functions of the form
\begin{equation*}\label{eq: MM original}
	f(x,y)=x \cdot \pi(y)+ h(y), \mbox{ for all } x, y\in\F_2^m,
\end{equation*}
where $\pi$ is a permutation of $\F_2^m$, and $h$ is an arbitrary Boolean function on $\F_2^m$. Throughout this paper, we will also use the finite-field representation of the $\mathcal{M}$ class, given by $f(x,y)=\operatorname{Tr}_1^m(x\pi(y))+h(y)$, where $x,y\in\F_{2^m}$, $\pi$ is a permutation of $\F_{2^m}$ and $h\in\mathcal{B}_m$.

The Maiorana-McFarland class is a natural generalization of Rothaus' family $R$ of the form $x\cdot y + h(y)$, which was discovered independently by Maiorana and McFarland~\cite{McFarland73_JCTA}, as mentioned by Dillon~\cite[p.~90]{Dillon74_PhD} and in~\cite{Dillon72_SURVEY}. The set of all bent functions that are extended-affine equivalent to bent functions from $\mathcal{M}$ is called the \emph{completed Maiorana-McFarland class} and is denoted by $\mathcal{M}^\#$, i.e.,
$$\mathcal{M}^\#=\{f(xA+b) + c \cdot x + d  \; :  f \in \mathcal{M}, A \in \operatorname{GL}(n,\F_2), b,c \in \F_2^n,  d \in \F_2\}.$$
As we will see in Remark~\ref{rem: GMM iff remarks}, every function that belongs to $\mathcal{M}^\#$ can be expressed in the form $f(xA)$, for some $f \in \mathcal{M}$ and $A \in \operatorname{GL}(n,\F_2)$.

The following characterization of bent functions in $\mathcal{M}^\#$ (called \emph{the second-order derivatives criterion}) was provided by Dillon in his thesis.
\begin{lemma} \cite[p. 102]{Dillon74_PhD}\label{lem M-M second}
	Let $n=2m$. A Boolean bent function $f\in\mathcal{B}_n$ belongs to $\cM^{\#}$ if and only if
	there exists (at least one) $m$-dimensional vector subspace $U$ of $\F_2^n$ such that, for all $a,b \in U$,
	\begin{equation*}\label{eq: 2nd order criterion}
	    D_{a}D_{b}f(x)=f(x) +  f(x +  a) +  f(x +  b) +  f(x +  a +  b)=0,
	\end{equation*}
	for all $ x \in \F_2^n$.
\end{lemma} 
Following the terminology introduced in~\cite{PolujanP20_DCC}, a vector subspace $U$ of $\F_2^n$ is called an \emph{$\mathcal{M}$-subspace} of $f\in\mathcal{B}_n$, if for all $a ,b  \in U$ it holds that $D_a D_b f=0$. In turn, the maximal dimension of an $\mathcal{M}$-subspace $U$ of $f\in\mathcal{B}_n$ is called the \emph{linearity index} of $f$, and is denoted by $\ind(f)$. Boolean bent functions with the minimal \emph{linearity index} (which is equal to 1) are called \emph{$\ell$-optimal}~\cite{KudinPPZ25_JOC}. These notions provide the main tools for studying the relationship between bent functions and the generalized Maiorana-McFarland class, an essential extension of the original Maiorana-McFarland class, and will be developed in detail in Sections~\ref{sec GMM class} and~\ref{sec M-subspaces}.

Since the introduction of the Maiorana-McFarland class in the early 1970s, many questions about the membership of various constructions and families in $\mathcal{M}^\#$ have been raised. For a summary of the initial results in this direction, we refer to Dillon's survey~\cite{Dillon72_SURVEY}. For example, thanks to the classification of quadratic forms~\cite[Chapter 15]{MacWilliamsSloane}, all quadratic bent functions are contained in the $\mathcal{M}^\#$ class. Due to the computer classification of cubic bent functions in $6$ variables~\cite{Rothaus76_JCTA,Dillon74_PhD} and in $8$ variables~\cite{Braeken2006}, it is also well-known that all cubic bent functions up to $8$ variables belong to the class $\mathcal{M}^\#$.

The first provable examples of bent functions outside $\mathcal{M}^\#$ were obtained by Dillon in his thesis~\cite{Dillon74_PhD}, where he identified concrete examples in $8$ variables and extended them to infinite families. Only years later, in~\cite{LangevinL11_DCC}, it was shown, thanks to computer search, that Maiorana-McFarland bent functions form only a tiny fraction of all bent functions.  It was estimated that the number of bent functions equivalent to Maiorana-McFarland in dimension eight is at most $2^{81}$ (with recent improvements in~\cite{LP_BFA_2024} and~\cite{KolomeecB25}), among all bent functions in dimension $8$, of which there are approximately $2^{106}$ in total~\cite{LangevinL11_DCC}. 

The Maiorana-McFarland construction is the simplest known explicit construction of bent functions, and due to the flexibility of its building blocks, it can be used efficiently to construct bent functions with additional cryptographic properties that are usually difficult to obtain. For example, for many years, the only known constructions of bent-negabent functions~\cite{SchmidtPP08}, rotation-symmetric bent functions~\cite{GaoZLC12}, and (homogeneous) cubic bent functions~\cite{CharpinK08_SIAM_JDM} were those obtained from the $\mathcal{M}^\#$ class. A drawback of bent functions in $\mathcal{M}^\#$ from a cryptographic viewpoint is that they can be written as a concatenation of affine functions defined on ``large'' vector spaces, a feature that may be used in certain attacks~\cite{CanteautV02}. Due to recent developments in the study of bent functions, several constructions of bent functions with the above mentioned  property of being outside $\mathcal{M}^\#$ have been obtained; see~\cite{ZhFR2015IT,PasalicKPP23}, \cite{PolujanKP2024WCC,PolujanKP2025_IEEE}, and~\cite{PolujanP20_DCC,ZhangPBW24_INFCOMP}.

\section{Carlet's \texorpdfstring{$\mathcal{C}$}{C} and \texorpdfstring{$\mathcal{D}$}{D} classes of bent functions and their superclasses}\label{sec C and D classes}
In this section, we will discuss two classes of bent functions introduced by Carlet \cite{Carlet93_EUROCRYPT} and provide an overview of some specific design methods. The main motivation for dedicating a whole section to these classes is that many explicit families of bent functions outside $\cM^\#$ were deduced within these classes, e.g., in \cite{ZhangPCW17_C2SI,ZhangPWC17_IEEE,LiuZPXZ20,ZhangCPW20_DAM,PasalicZKW21_DAM,ZhaoWZPC22_DCC,BapicP22_DAM,BapicPZH22_CCDS,PasalicBZW23_DCC,ZhangPBW24_INFCOMP,PasalicBZW24_DCC}. 

By modifying the $\cM$ class, as originally proposed by Dillon \cite{Dillon74_PhD}, Carlet~\cite{Carlet93_EUROCRYPT} derived two new classes of bent functions. Denoting by $1_E$ the indicator of $E\subseteq\F_2^n$, the class $\cal D$ consists of bent functions in $n=2m$ variables of the form
\begin{equation*}
\label{eqD}
f(x, y)=x\cdot\pi(y) + 1_{E_1}(x) 1_{E_2}(y),
\end{equation*}
where $\pi$ is a permutation of $\F_2^{m}$ and $E_1,E_2$ are two vector subspaces
of $\F_2^m$ such that
\begin{equation}
    \label{eq: D property} \tag{$D$}
    \pi(E_2)=E_1^{\perp}.
\end{equation}
In this case, we say that $\pi$ \emph{has property~\eqref{eq: D property}} with respect to $E_1$ and $E_2$. Notice that $\dim(E_1) + \dim(E_2)=m$, which then implies that $f$ is at Hamming distance $2^m$ from the function $x \cdot \pi(y)$ which is known to be minimal in terms of preserving bentness \cite{KolomeetsPavlov2009,KolomeetsPavlov2010}. More precisely, in~\cite{KolomeetsPavlov2009,KolomeetsPavlov2010}, it is shown that two bent functions in $2m$ variables are at the minimum possible distance between distinct bent functions $2^m$ if and only if they differ on an affine subspace of dimension $m$, on which both functions are affine (see also the works~\cite{Kolomeets2012,KolomeecB25}).

An explicit subclass of $\cal D$ that uses the indicator  $E_1 \times E_2=\{0_m\} \times \F_2^m$, denoted by
${\cal D}_0$,  contains all elements of the form $x\cdot\pi(y)+\delta_0(x)$, where $\delta_0(x)$ is the \emph{Dirac symbol}, which is equal to 1 if $x=0_m$, and 0, otherwise. It has been shown that ${\cal D}_0$ contains bent functions outside the $\cal M^{\#}$ and $\cal PS^{\#}$ classes~\cite{Carlet93_EUROCRYPT}.

The second Carlet class $\cal C$ of bent functions contains all functions of the form
\begin{equation*}
\label{eqC}
f(x,y)=x\cdot\pi(y)+ 1_{L^{\perp}}(x),
\end{equation*}
where $L$ is any vector subspace of $\F_2^m$ and $\pi$ is any permutation
of $\F_2^m$ such that:
\begin{equation}
    \label{eq: C property} \tag{$C$}
    \phi(a + L) \; \mbox{is a flat (i.e., an affine subspace) for any } a\in \F_2^m \mbox{ (where } \phi=\pi^{-1} \mbox{)}.
\end{equation}
Then, the pair {\em $(\phi,L)$ is said to have property~\eqref{eq: C property}}, which characterizes the bentness of $f$.

\subsection{\texorpdfstring{$\mathcal{C}$}{C} class of bent functions}\label{sec C class} 
We now discuss the property~\eqref{eq: C property} in more depth, together with sufficient conditions under which bent functions in the $\mathcal{C}$ class are provably outside $\cM^\#$.

\subsubsection{Satisfying the \texorpdfstring{$(C)$}{C} property}

The initial work on specifying bent functions within the $\cC$ class was presented in \cite{MandalSGP16_FUIN}, where it was observed that satisfying the~\eqref{eq: C property} property is not difficult. 
\begin{remark}\label{rem:Cproperty}
The fact that we can construct many bent functions in the $\cC$ class of bent functions is not difficult. One could take two subspaces $L, M$ in $\F_2^m$ of the same dimension and partition $\F_2^m$ into $\bigcup_{a \in A}(a+L)$ and $\bigcup_{b \in B}(b+M)$, with $A, B$ being subsets of $\F_2^m$ of the same cardinality $|A| =|B|$, and then take any permutation $\phi$ that maps the elements of $\{a + L \colon  a \in  A\}$ onto the elements of $\{b + M \colon  b \in  B\}$. The pair $(\phi, L)$ would satisfy property~\eqref{eq: C property}.
\end{remark}

However, the problem with the above construction method is two-fold. Firstly, there is no nice algebraic description of the mappings suggested in Remark \ref{rem:Cproperty}, and even more importantly, there is no guarantee that the constructed bent function is outside $\cM^\#$. Moreover, it also became evident that certain permutations do not admit any subspace $L$ of fixed dimension such that 
$f(x,y)=x\cdot\pi(y)+ 1_{L^{\perp}}(x)$ is a bent function in $\cC$, as the following result demonstrates (see also Theorem~\ref{th:nonexist}).
\begin{theorem}\label{th:Houquadratic}
\cite{MandalSGP16_FUIN}
Let $m\geq 1$ and $\phi(y)=ay+b y^{2^m}+ y^{2^{m+1}-1}$ be a  permutation polynomial over $\F_{2^{2m}}$ (see~Hou~\cite[Theorem B]{Hou13} for explicit criteria). Then, there exists no $2$-dimensional vector subspace of $\F_2^{2m}$ such that $(\phi,L)$ has property~\eqref{eq: C property}.
\end{theorem}

In addition to various (general) sufficient conditions that the permutation $\pi$ must meet, as discussed in \cite{MandalSGP16_FUIN}, a family of so-called $k$-linear split permutations was considered for the purpose of specifying bent functions in $\cC$ explicitly. Two polynomials $\phi, \psi \in \F_{2^m}[y]$ are said to be {\em linearly equivalent} if there exist bijective $\ell_1, \ell_2 \in \mathcal{L}(m)$ such that
$\ell_1 \circ \phi \circ \ell_2 = \psi,$ where $\mathcal{L}(m)$ is a set of linearized polynomials (of the form $l(y) = \sum_{i=0}^{m-1}a_i y^{2^i} \mbox{ with } a_i \in \F_{2^m}$).
\begin{lemma}\cite{MandalSGP16_FUIN}
\label{lem:connection-pi-phi}
Suppose $\pi$ and $\phi$ are two linearly equivalent permutations on $\F_{2^m}$
such that $\phi = \ell_1 \circ \pi \circ \ell_2$ where $\ell_1, \ell_2 \in \mathcal{L}(m)$, and
$L$ is a vector subspace of $\F_{2^m}$.
If $\pi( a + L)$ is a flat for all $a \in \F_{2^m}$, then $\phi(a + \ell_2^{-1}(L))$ is
also a flat for all $a \in \F_{2^m}$.
\end{lemma}

It should be noticed that certain generic families of bent functions in $\cC$ could be derived from a so-called {\em $k$-linear split}  polynomial  $\phi \in \F_{2^m}[y]$  of the form
\begin{equation*}
\phi(y) = \ell_1(y)\ell_2(y)\cdots\ell_k(y) \mbox{ with } \ell_i \in \mathcal{L}(m), 1\leq i\leq k.
\end{equation*}
In particular, when $k=2$, such a polynomial is called bilinear, and the most typical representative is a quadratic permutation of the form $\phi(y)=y^{2^r+1}$, for a suitably chosen $r$. Nevertheless, assuming that $\gcd(m,r)=e$ and $m/e$ is odd (to ensure that $\phi(y)=y^{2^r+1}$ is a permutation of $\F_{2^m}$), it was shown in \cite{MandalSGP16_FUIN} that if $\dim(L)=2$ then the case $e=1$ cannot lead to a bent function 
$f(x,y)=x \cdot \pi(y) + 1_{L^{\perp}}(x)$. However, there is a possibility of specifying bent functions in $\cC$ when $e=2$ and $\dim(L)=2$, as demonstrated by Example 5.9 in \cite{MandalSGP16_FUIN}. For clarity, we recall this example as we will refer to it in the next section.
\begin{example}
    \label{ex:explicitCclass}
Let $m = 2p$ where $p$ is any odd prime, $r = 2$ and
$e = \gcd(m, r) = 2$. Since $m/e$ is odd, it is known that
$\gcd(2^r+1, 2^m-1)=1$. Therefore $\phi(y) = y^{2^r+1}$ is a
permutation on~$\F_{2^m}$. Let $\zeta$ be a primitive element of
$\F_{2^m}$. Therefore, $\lambda = \zeta^{\frac{2^m-1}{2^e-1}}= \zeta^{\frac{2^m-1}{3}}$
is a generator of~$\F_{2^e}$. Suppose that
$\pi(y) = \phi^{-1}(y) = y^{\gamma}$, where
$\gamma (2^r+1) \equiv 1 \pmod{2^m-1}$. Given $r$ and $m$, $\gamma$
can be computed easily by the Euclidean algorithm. Consider the Maiorana-McFarland bent function $f(x,y)=\operatorname{Tr}^m_1(x\pi(y))$, with $x,y \in \F_{2^m}$. According to Theorem 5.8 in \cite{MandalSGP16_FUIN}, if we choose
$L = \langle 1, \lambda \rangle$, then the function
$f^*(x,y) = \operatorname{Tr}^m_1(x\pi(y)) + 1_{L^{\perp}}(x)$ is in $\cC$. The bent function
$f^*$ can be explicitly written as
\begin{equation*}
\label{c-class-ex}
\begin{split}
f^*(x,y) &= {\rm Tr}_1^m(xy^{\gamma}) + ({\rm Tr}_1^m(x) + 1)({\rm Tr}_1^m(\lambda x) + 1)\\
&=  {\rm Tr}_1^m(xy^{\gamma})
+ {\rm Tr}_1^m(x){\rm Tr}_1^m(\lambda x) + {\rm Tr}_1^m((1+\lambda) x) + 1.\\
\end{split}
\end{equation*}
\end{example}

The above discussion emphasizes the difficulty of satisfying the property~\eqref{eq: C property} when using permutations with nice algebraic descriptions. We note that the property of being outside $\cM^\#$ was not discussed in \cite{MandalSGP16_FUIN}.

\subsubsection{Some efficient design methods of bent functions in \texorpdfstring{$\mathcal{C}$}{C} outside \texorpdfstring{$\mathcal{M}^\#$}{M\#}}
As pointed out in the previous section, specifying tuples $(\phi,L)$ that satisfy the property \eqref{eq: C property} (assuming that $\phi$ has a simple algebraic representation) is, in general, a difficult task, and in certain cases it is even impossible, see Theorem \ref{th:Houquadratic}. However, certain progress in this context has been achieved, which is the main topic of this section.
The main tool for analysing the exclusion from the class $\cM^\#$ is the result of Dillon, namely Lemma \ref{lem M-M second}, which is almost exclusively utilized in such proofs. 

Since we are using bivariate representation to represent $f \in \cC$, it is convenient to represent an element $a \in \F_2^{2m}$ as $a=(a_1,a_2) \in \F_2^m \times \F_2^m$. Then, to  prove  that a bent function $f \in \mathcal{B}_{2m}$  does not belong to $\mathcal{M}^{\#}$, applying Lemma \ref{lem M-M second}, one  needs to show that there does not exist  an  $m$-dimensional  subspace $V$ such that
   $$D_{(a_1, a_2)}D_{(b_1, b_2)}f(x,y)=0,$$
   for all  $a, b \in V$. Such computations are commonly tedious (see, for instance, the proof of Theorem \ref{th main Ctype} below in \cite{ZhangPCW17_C2SI}) because one needs to cover several cases with respect to the structure of $V$. For a bent function $f(x,y) = x \cdot \pi(y) + 1_{L^{\perp}}(x)$ it is easily confirmed that
   \begin{equation}\label{eq:2ndderivCclass}
       \begin{split}
           D_{(a_1, a_2)}D_{(b_1, b_2)}f(x,y)
        &=x\cdot\left( D_{a_2}D_{b_2}\pi(y)\right)
 +  a_1\cdot D_{b_2}\pi(y + a_2)\\
     &+ \;  b_1 \cdot D_{a_2}\pi(y + b_2)+ D_{a_1}D_{b_1}1_{L^{\bot}}(x).
       \end{split}
   \end{equation}

The first explicit condition imposed on $(\phi,L)$ in order to satisfy the property~\eqref{eq: C property} and produce a bent function outside $\mathcal{M}^{\#}$ was given in~\cite{ZhangPCW17_C2SI}. Notice that due to the term $a_1\cdot D_{b_2}\pi(y + a_2)$ (or $b_1 \cdot D_{a_2}\pi(y + b_2)$) in Eq.~\eqref{eq:2ndderivCclass}, the condition below is quite natural.
\begin{theorem}\label{th main Ctype}
  Let $n=2m >4$ be an even integer
 and let $f(x,y)=x\cdot \pi(y) + 1_{L^{\perp}}(x)$, where $L$ is any vector subspace of  $\F_2^{m}$ and $\pi$ is a permutation on  $\F_2^{m}$ such that $(\pi, L)$ has property~\eqref{eq: C property}.
 If $\pi$  satisfies:
\begin{enumerate}[i)]
  \item  $\dim(L)\geq 2$,
  \item  $b \cdot \pi(y)$ has no nonzero linear structure for any nonzero $b \in \F_2^m$,
\end{enumerate}
then $f$ does not  belong to $\mathcal{M}^{\#}$.
\end{theorem}
\begin{remark}\label{rem:wronguse}
    Unfortunately, the authors in \cite{ZhangPCW17_C2SI} used an incorrect formulation in $ii)$ ``$\pi$ has no nonzero linear structure'' which means that $\pi(y)+\pi(y+a)$ is not a constant vector in $ \F_2^m$, for any nonzero $a \in \F_2^m$, instead of the correct one above. 
\end{remark}
It was exactly Example \ref{ex:explicitCclass} (originally Example 5.9 in \cite{MandalSGP16_FUIN}) to which Theorem \ref{th main Ctype} was applied to demonstrate the existence of infinite instances of bent functions in $\cC$ that are outside $\cM^\#$. 
\begin{lemma}\cite{ZhangPCW17_C2SI}
	For $r\neq 0$, the function $f^*$ from Example \ref{ex:explicitCclass} does not belong to $\cM^\#$.
\end{lemma}
This method was later generalized in \cite{ZhangCPW20_DAM}, using the fact that the compositional inverses of quadratic monomial permutations $y^d$ (thus $\operatorname{wt}(d)=2$) are not necessarily quadratic. Indeed, setting  $m=6$ and $r=2$ in Example \ref{ex:explicitCclass} one can easily verify that  $\gamma (2^r+1) \equiv 1 \pmod{2^m-1}$, which gives $\gamma=38$ and $\pi(y)=y^{38}.$ Then, the permutation $\pi$ over $\F_2^6$ is thus cubic and by Theorem 5 in \cite{CharpinGohar2010} none of its component functions admit linear structures. Therefore, by Theorem \ref{th main Ctype}, the function $f^*$ is outside $\mathcal{M}^{\#}$. 

To ensure that the inverse permutation of $y^{2^r+1}$ is non-quadratic, the following result is useful.  
\begin{lemma}\label{lem:cubicexpo}\cite{ZhangCPW20_DAM}
Let $m >  4$ be a positive integer and define the permutation $\phi(y)=y^{2^r+1}$, where $\gcd(m,r)=e$ and $m/e$ is odd. Then, the unique $1 < \gamma <2^m-1$ which satisfies $\gamma (2^r+1) \equiv 1 \pmod{2^m-1}$ is such that $\operatorname{wt}(\gamma)>2$. In other words, for any  $1 \leq i \neq j\leq m-1$, we always have
 \begin{equation*}\label{equ leminv0}
 (2^i+2^j)(2^r+1)~/\!\!\!\!\!\equiv 1\mod (2^m-1).
 \end{equation*}
\end{lemma}

\paragraph{Designing permutations whose components do not admit linear structures}
As already pointed in Remark \ref{rem:Cproperty} the~\eqref{eq: C property} property of a tuple $(\pi,L)$ can be easily achieved by splitting the ambient space $\F_2^m$ into cosets and ensuring that $\pi^{-1}(a + L)$ is a flat for any $a \in \F_2^m$. However, the partition of $\F_2^m$ must be performed carefully to ensure that the components of $\pi$ do not admit linear structures, as the following result suggests. 
\begin{proposition}\label{prop:trivialdecompCclass} \cite{KudinPCZ22_CCDS} Let $L$ and $P$ be two $k$-dimensional subspaces of $\mathbb{F}_2^m$, and let $\pi \colon \mathbb{F}_2^m \rightarrow \mathbb{F}_2^m$ be a permutation such that $\pi$ maps cosets of $L$ to cosets of $P$, i.e., $\pi(a+L)= \pi(a) + P$, for every $a \in \mathbb{F}_2^m$. Then, $\pi$ has at least $2^{m-\dim(L)}$ component functions such that their space of linear structures contains $L$.
\end{proposition}
This shows that the simple idea of taking two (not necessarily different) subspaces of $\mathbb{F}_2^m$, and mapping their cosets in some bijective way to construct a permutation with the~\eqref{eq: C property} property, will necessarily produce a permutation whose components admit non-trivial linear structures. Therefore,  more sophisticated (non-trivial) partitions of $\F_2^m$ are required; this problem was originally addressed in \cite{BaNe1975}. 

We briefly recall the discussion in Section 3.2 in \cite{KudinPCZ22_CCDS}, where the main goal was to establish an isomorphism between a finite field and a vector subspace. One starts by choosing an arbitrary subspace $E$ of $\F_2^m$ whose dimension $d=\dim(E)$ satisfies $1\le d\le (m-1)/2$. Then, we can take  another subspace $V$ of $\mathbb{F}_2^m$ such that $\mathbb{F}_2^m$ is the direct sum of $E$ and $V$, i.e., $E \cap V = \lbrace 0_m \rbrace $, and $E \oplus V = \mathbb{F}_2^m$. Moreover, we will need $d+1$ permutations $I_V, \sigma_1, \ldots, \sigma_d$ on $V$, (here $I_V$ is the identity permutation on $V$) such that any non-trivial linear combination of them is again a permutation of $V$. To achieve this, the subspace $V$ is interpreted as the finite field $\mathbb{F}_{2^{m-d}}$ in the following way. We define a mapping $\mathcal{L}$ between the vector subspace $V$ and a finite field $\FF_{2^{m-d}}$. Let ${u_1,\ldots , u_{m-d}}$ be the basis of the subspace $V$ and let $\alpha$ be a primitive element of the finite field $\FF_{2^{m-d}}$. Then, we  define
$\mathcal{L}\colon u_i \mapsto \alpha^{i-1}$, whereas the remaining elements are mapped in such a way as to preserve isomorphism between additive structures. That is,
\begin{equation*}\label{eq:isomorp}
\mathcal{L}\left(\sum _{i = 1}^{m-d} \beta_i u_i \right) =  \sum _{i = 1}^{m-d} \beta_i \alpha ^{i-1},
\end{equation*}
where $\beta_i \in \FF_2$.
 Using such a mapping, one can induce a multiplication operation on elements of the subspace $V$ by defining
\begin{equation*} \label{eq:mul}
v_i \bullet v_j = \mathcal{L}^{-1}(\mathcal{L}(v_i)\mathcal{L}(v_j)),  \text{ for any }v_i,v_j\in V.
\end{equation*}
With $`` \bullet "$ defined in this way, $(V, + , \bullet )$ has the structure of a finite field.  Moreover, $\mathcal{L}$ becomes a field isomorphism between $V$ and $\FF_{2^{m-d}}$. The above-mentioned permutations 
$I_V, \sigma_1, \ldots, \sigma_d$ on $V$ can also be defined efficiently; for details see Section 3.2 in \cite{KudinPCZ22_CCDS}. It was shown in~\cite{KudinPCZ22_CCDS} that the constructed isomorphism can be used to obtain permutations that satisfy the~\eqref{eq: C property} property and have no components with linear structures. Moreover, the following result holds (note that $1_V$ denotes the multiplicative identity in $(V,+,\bullet)$).
\begin{theorem} \label{th:constrnontrivial} Let $E$ and $V$ be two subspaces of  $\mathbb{F}_2^m$ such that $ E \cap V = \lbrace 0_m \rbrace$, $E \oplus V = \mathbb{F}_2^m$, and $d= \dim(E)$ satisfies $1 \leq d \leq (m-1)/2$. Let  $1_V,w_1, \ldots , w_d \in V$ be linearly independent vectors and $\sigma_i$ be a permutation over $V$ defined by $\sigma_i(v)=w_i\bullet v$. Let $r$ be such that $\gcd (2^{m-d}-1,r)=1$, $\operatorname{wt}(r) \geq 3$, and let  $\varphi(v)=  \underbrace{v \bullet v \bullet \cdots \bullet v}_{r \text{-times}}=v^{r}$, for all $ v \in V$. Let $\psi$ be a permutation of $E$ whose component functions do not admit linear structures. Then, the function $\phi \colon \mathbb{F}_2^m \rightarrow \mathbb{F}_2^m$ defined by  Eq. (6) and Eq. (7) in \cite{KudinPCZ22_CCDS} is a permutation of $\mathbb{F}_2^n$ whose component functions do not admit (nonzero) linear structures.
\end{theorem}

This result allows us to construct infinite families of $2m$-variable bent functions in $\cC$ and outside $\cM^\#$, though for relatively large $n=2m$. As remarked in \cite{MandalSGP16_FUIN}, if $L$ has dimension 1  then $\pi^{-1}(a+L) = \phi(a + L)$ is always a one-dimensional
flat and in this case $1_{L^{\perp}}(x)$ is an affine function. Therefore, the necessary condition that $f(x,y)=x \cdot \pi(y) + 1_{L^\perp}(x)$ is outside $\cM^\#$ is that $\dim(L)\geq 2$. However, for a permutation $\pi$ of $\mathbb{F}_2^m$  there is a trade-off between the property that no component functions of $\pi$ admit linear structures, and the requirement that $\pi^{-1}$ satisfies the \eqref{eq: C property} property for some subspace $L$. 

\begin{theorem} \label{th:nonexist} \cite{KudinPCZ22_CCDS} Let $m$ and $d$ be two integers such that $d \geq 1$ and $m > 2^d \left( d-1 \right) +1$. Let $\pi \colon \mathbb{F}_2^m \rightarrow \mathbb{F}_2^m$ be a permutation. If there exists a vector subspace $L$ with dimension $\dim(L) = m-d$, such that $ \left( \pi ^{-1} , L \right)$ has the \eqref{eq: C property} property, i.e., $ \pi ^{-1} ( a + L )$ is an affine subspace for every $a \in \mathbb{F}_2^m$, then $\pi$ has at least one component function that admits (nonzero)  linear structures.
\end{theorem}

\paragraph{Employing Rothaus' construction}
A more complicated approach was taken in \cite{ZhangPWC17_IEEE,LiuZPXZ20}, where a secondary construction method of Rothaus of obtaining new bent functions in $n+2$ variables starting with three suitable bent functions in $n$ variables was employed. 
\begin{definition}[Rothaus' construction]\cite{Rothaus76_JCTA} Let
$x=(x_1,x_2,\ldots, x_n)\in \F_2^n$ and $x_{n+1},x_{n+2}\in
\F_2$. Let $A,B,C$ be bent
functions on $\F_2^{n}$ such that  $A+
B + C$ is bent as well. Then, the function
defined for every element $(x,x_{n+1},x_{n+2})\in \F_2^{n+2}$
by
\begin{equation}\label{eq:Rothaus}
	\begin{split}
		 f(x,x_{n+1},x_{n+2})&= A(x)B(x)+ A(x)C(x)+ B(x)C(x)+ x_{n+1}x_{n+2} \\
  							 &+ [A(x)+ B(x)]x_{n+1}+ [A(x)+C(x)]x_{n+2}
	\end{split}
\end{equation} 
is a bent function in $n+2$ variables.
\end{definition}
Notice that the question of whether the above condition is also necessary was raised in \cite{ZhangPWC17_IEEE}, and was positively answered in \cite{GuoWangGong2025}. 

In \cite{ZhangPWC17_IEEE}, the above construction was used in a particular way. The following lemma gives more precision to a known result on the minimum distance of bent functions, where the connection to Eq. \eqref{eq:Rothaus} is the assignment $$A \to f_0,\quad B \to f_1 \quad \mbox{and} \quad C \to f_0 + 1_{\Delta}.$$ 
Notice that ensuring $A+B+C$ is bent is then equivalent to $f_1+1_{\Delta}$ is bent.
\begin{lemma}\label{lemma 1}\cite{ZhangPWC17_IEEE}
Let $f_0,f_1$ be two bent functions in $n$ variables. Let $\Delta\subseteq\FF_2^n$ be a set and $1_{\Delta}(x)$ be the function that equals 1  if $x$ is in $\Delta$, otherwise it equals  0.
If  either $\Delta\subseteq (\operatorname{supp}(f_0)\cap \operatorname{supp}(f_1))$ or $\Delta\subseteq (\operatorname{supp}(1 + f_0)\cap \operatorname{supp}(1 + f_1))$
and $f_0 + 1_{\Delta}, f_1 + 1_{\Delta}$ are also bent, then we have $|\Delta |=2^{\frac{n}{2}}$.
\end{lemma}

The following result, as a special case of Rothaus' construction, was heavily employed in the design of bent functions outside $\cM^\#$.
\begin{theorem}\label{theo:RothausMain} \cite{ZhangPWC17_IEEE}
   Let
$x=(x_1,\ldots, x_n)\in \F_2^n$ and $x_{n+1},x_{n+2}\in
\F_2$. Let $f_0,f_1$ be two Boolean functions in $n$ variables such that  $f_0(x)=f_1(x)$ for $x\in \Delta$, where $\Delta\subseteq\FF_2^n$.
Then, the function $f$ defined as
\begin{equation}\label{eq:Rothaus-usingf0andf1}
	\begin{split}
		   f(x,x_{n+1},x_{n+2}) &=(x_{n+1}+ x_{n+2}+ 1)f_0(x)+ (x_{n+1}+ x_{n+2})f_1(x) \\
 		& + (x_{n+2}+ 1)1_{\Delta}(x) + x_{n+1}x_{n+2}+ x_{n+1},
	\end{split}
\end{equation}  
is bent if and only if all the functions $f_0,f_1$, $f_0 + 1_{\Delta}$, $f_1 + 1_{\Delta}$ are bent. 
\end{theorem}

As noted in \cite{ZhangPWC17_IEEE}, the existence of functions satisfying the conditions in Theorem~\ref{theo:RothausMain} is easily confirmed by considering the $\mathcal{PS}^-$ class (see Section~\ref{subsec:classH} for the formal definition) and defining $\Delta$ to be an $\frac{n}{2}$-dimensional vector subspace. Moreover, the connection to the class $\mathcal{D}$ is easily seen through the specification of the indicator $1_{\Delta}(x,y)=1_{E_1}(x)1_{E_2}(y)$ so that both $f_0 + 1_{\Delta}$ and $f_1 + 1_{\Delta}$ belong to $\mathcal{D}$, where $|\Delta|=2^{n/2}$ and additionally assuming that $f \in \cM$. A similar reasoning applies when the class $\cC$ is considered, simply by defining $\Delta =L^\perp \times \F_2^{n/2}$, assuming that $f_0(x,y)=x \cdot \pi(y)$, with $x,y \in \F_2^{n/2}$. The main result, given below, related to the class $\cD$ (quite similarly for $\cC$) is a set of sufficient conditions that $f$ defined by Eq. \eqref{eq:Rothaus-usingf0andf1} is outside $\cM^\#$ \cite[Theorem 3]{ZhangPWC17_IEEE}. 

\begin{theorem}\label{th:mainDoutsideMM}
 Let $n=2m >4$ be an even integer
 and let $f_0(x,y)=x \cdot \pi(y)$,   $f_1(x,y)=x \cdot \phi(y)$, where  $\pi$ and  $\phi$  are two permutations of  $\F_2^{m}$.
 Then, for any $x, y \in \F_2^m$, the function $f \in \mathcal{B}_{n+2}$ defined as
 \begin{equation}\label{eq:corol-twofunctions1 plus00}
 \begin{split}
 f(x,y, x_{n+1},x_{n+2})=&(x_{n+1}+ x_{n+2}+ 1)f_0(x,y)+ (x_{n+1} + x_{n+2})f_1(x,y) \\
+&(x_{n+2}+ 1)1_{E_1}(x)1_{E_2}(y)  +  x_{n+1}x_{n+2}+ x_{n+1}, 
\end{split}
\end{equation}
is  bent, where $E_1, E_2$ are
 vector subspaces of  $\F_2^{m}$ such that $\pi(E_2)=E_1^{\perp}$ (resp. $\phi(E_2)=E_1^{\perp}$).
 Further, if $\pi$ and  $\phi$ satisfy:
\begin{enumerate}[i)]
  \item  $u \cdot \pi(y)$ (or $u \cdot \phi(y)$)  has no nonzero linear structure, for any $u \in \F_2^{m}\backslash \{0_{m}\}$,
  \item  $ \nu\cdot (\pi + \phi)\neq constant$ for $\nu\in \F_2^{m}\backslash \{0_{m}\}$,
  \item   $\max\limits_{\nu\in \F_2^{m}}\deg(\nu\cdot (\pi + \phi))\geq 2$,
   \item    $ E_1 \subset\F_2^{m}$  and $\dim(E_1)\leq m-2$~~(that is, $\deg(1_{E_1}(x))\geq 2$),
\end{enumerate}
then $f$ does not  belong to $\cM^{\#}$.
\end{theorem}

\begin{remark}\label{rem:Complicatedanalysis}
 Compared to the simple ANF expression for the second-order derivative given by Eq. \eqref{eq:2ndderivCclass}, the related expression concerning  $f$ defined by Eq. \eqref{eq:corol-twofunctions1 plus00} is more complicated and the proof that $f \not \in \cM^\#$ is more tedious. Similarly to Remark \ref{rem:wronguse}, 
 as pointed out in \cite[Remark 1]{LiuZPXZ20}, the incorrect statement that $\pi$ (or $\phi$) does not admit linear structures is replaced by the correct one in item $i)$.  
\end{remark}
\begin{remark}
    In Section IV of \cite{ZhangPWC17_IEEE}, the authors proposed an iterative use of the Rothaus construction using the necessary and sufficient condition that $A(x)B(x) + A(x)C(x) + B(x)C(x)=A(x)+B(x)+C(x)$, for three bent functions $A(x),B(x)$ and $C(x)$. However, it was shown in \cite{GuoWangGong2025} that the only solution to this equation is $A(x)=B(x)=C(x)$, which then reduces this method to the well-known direct sum method.
\end{remark}
In the follow-up article \cite{LiuZPXZ20}, the sufficient conditions in Theorem \ref{theo:RothausMain} were significantly simplified. More precisely, item $i)$ remained the same whereas $ii)-iv)$ were replaced by the conditions that $\deg(\pi)>2$ and $\deg(\phi)=2$, see Theorem 3 in \cite{LiuZPXZ20} for more details.

\begin{openproblem}
    Specify the exact value of the linearity index of bent functions in $\mathcal{C}^\#\setminus\mathcal{D}_0^\#$ that do not belong to $\mathcal{M}^\#$. In particular, prove or disprove that $\mathcal{C}^\#\setminus\mathcal{D}_0^\#$ contains $\ell$-optimal bent functions.
\end{openproblem}

\subsection{Designing bent functions in \texorpdfstring{$\mathcal{D}$}{D} outside \texorpdfstring{$\mathcal{M}^\#$}{M\#}}\label{sec D class}
The main difference between the classes $\cC$ and $\cD$ is that the latter class (with members $f(x,y)=x \cdot \pi(y) + 1_{E_1}(x) 1_{E_2}(y)$) only allows for the modifications on $n/2$ dimensional subspaces $E_1 \times E_2$ since $\dim(E_1) + \dim(E_2)=n/2$. Therefore, the sets of sufficient conditions derived in, e.g., \cite{ZhangPCW17_C2SI,ZhangPWC17_IEEE,LiuZPXZ20}, related to the exclusion of $f$ from $\cM^\#$, are quite similar to those that regard the $\cC$ class, for instance see Theorem \ref{th:mainDoutsideMM}. Notice, however, that the~\eqref{eq: D property} property (i.e., $\pi(E_2)=E_1^\perp$) is less restrictive compared to the~\eqref{eq: C property} property since no cosets are considered, thus not requiring that $\pi(E_2 +a)=b+ E_1^\perp$. For instance, the following result captures the main conditions related to the $\cD$ class. 
\begin{theorem}\label{th main Dtype}\cite[Theorem 2]{ZhangCPW20_DAM}
Let $n=2m\geq 8$ be an even integer and let $f(x,y)=x\cdot \pi(y) + 1_{E_1}(x)1_{E_2}(y)$, where $\pi$ is a permutation of $\F_2^{m}$, and $E_1,E_2$ are two vector subspaces of $\F_2^{m}$ such that $\pi(E_2)=E_1^{\perp}$. If $(\pi,E_1,E_2)$ satisfies:
\begin{enumerate}[i)]
\item $\dim(E_1)\geq 2$ and $\dim(E_2)\geq 2$,
\item $u\cdot \pi$ has no nonzero linear structure for all $u\in\F_2^{m}\setminus \{0_m \}$,
\item $\deg(\pi)\leq m-\dim(E_2)$,
\end{enumerate} 
then $f$ is a bent function in $\cD$ outside $\cM^{\#}$.
\end{theorem}
The above result still assumes that $\pi$ satisfies the~\eqref{eq: D property} property so that $f$ is bent. In the special case, 
when $m$ is even and $\dim(E_2)=2$, we can be more precise about the~\eqref{eq: D property} property. 
\begin{proposition}\label{pro:existDeven}\cite{ZhangPCW17_C2SI}
	Let $m$ be even. Then, any monomial permutation $\pi(y)=y^d$ of $\F_{2^m}$, where  $3\leq deg(\pi )\leq m-2$, satisfies the required conditions in Theorem~\ref{th main Dtype} for the $2$-dimensional vector subspace $E_2=\langle \zeta^{\frac{2^m-1}{3}}, \zeta^{\frac{2(2^m-1)}{3}} \rangle$, where $\zeta$ is a primitive element of  $\F_{2^m}$. Thus, $f(x,y)=x\cdot \pi(y) + 1_{E_1}(x)1_{E_2}(y)$, where $\pi(E_2)=E_2=E_1^{\perp}$, is a bent function in $\mathcal{D}$ and outside $\mathcal{M}^{\#}$.
\end{proposition}

The following example illustrates this approach. 
\begin{example}
	Let $m=6$ and $d=11$. Since $(2^6-1,11)=1$ and  $\operatorname{wt}(11)=3$, then $\pi(y)=y^d$ is a cubic permutation over $\F_2^6$. Using the programming package Magma~\cite{BosmaCP97},  the vector space representation on $\F_2^6$ of the subspace $E_2=\langle \zeta^{21}, \zeta^{42}\rangle $, where $\zeta$ is the primitive element of the field $\F_{2^6}$, is:
	$$ E_2 = \left\lbrace   \begin{array}{c}
	(0,0,0,0,0,0), \\
	(1, 0, 0, 0, 0, 0),\\
	(1, 1, 1, 1, 0, 0), \\
	(0, 1, 1, 1, 0, 0)
	\end{array} \right\rbrace.   $$
Since $1^{11}=1, (\zeta^{21})^{11}=\zeta^{42},$ and $(\zeta^{42})^{11}=\zeta^{21}$, the subspace $E_2$ is indeed mapped to itself. This implies that $\pi(E_2)=E_2=E_1^{\perp}$ and
$$ E_1 = \left\langle    \begin{array}{c}
	(0,1,0,1,0,0), \\
	(0, 0, 1, 1, 0, 0),\\
	(0,0,0,0,1, 0), \\
	(0,0,0,0,0,1)
	\end{array} \right\rangle .$$
Thus, all the requirements of  Theorem~\ref{th main Dtype} are satisfied and the permutation $\pi$ gives rise to a bent function $f(x,y)=x\cdot \pi(y)+ 1_{E_1}(x)1_{E_2}(y)$ on $\F_2^{12}$ contained in $\cD$ but outside $\cM ^\#$.
\end{example}

\paragraph{A generic design approach with provable exclusion from \texorpdfstring{$\mathcal{M}^\#$}{M\#}}
 A more general approach to constructing bent functions outside $\cM^\#$ was considered in \cite{ZhangCPW20_DAM}, where a special class of permutations was employed. The authors denoted by $\cD_k$ the subclass of bent functions in $\cD$ having the property that
$$\dim(E_1)=k\quad \mbox{and} \quad \dim(E_2)=m-k,$$
where the original \emph{$\mathcal{D}_0$ class} of bent functions introduced by Carlet~\cite{Carlet93_EUROCRYPT} is obtained for the special case $k=0$. This class will be considered in detail in the following section. By fixing $\dim(E_1)=k$ for $k \in \{1,\ldots, m-2\}$,  an explicit class $\cD_k^\star \subset \cD_k$ (where ``$\star$'' denotes the use of a special class of permutations) of bent functions was proposed. This ``special'' class of permutations on $\F_2^m$ is obtained from the identity permutation by swapping the values at two distinct points, and is given by
\begin{equation}\label{eq:pistar}
\sigma _{\mathbbm{e}_l , \mathbbm{e}_t}(y)=\left \{ \begin{array}{ll}
y, & y \not \in \{\mathbbm{e}_l,\mathbbm{e}_t\}, \\
\mathbbm{e}_l, & y= \mathbbm{e}_t, \\
\mathbbm{e}_t, & y= \mathbbm{e}_l,
\end{array} \right .
\end{equation}
where $l,t\in \{1,2,\ldots,m\}$ with $l\neq t$, and furthermore $\mathbbm{e}_l,\mathbbm{e}_t \in \F_2^m$ denote elements in the canonical basis of $\F_2^m$. More precisely, $ ( \mathbbm{e}_{l})_i=1$ if and only if $i=l$, otherwise $( \mathbbm{e}_{l})_i=0$.  W.l.o.g., we assign $l=1,t=2$, and then it can be readily verified that selecting $E_1=\langle \mathbbm{e}_1, \ldots,\mathbbm{e}_k \rangle$ and $E_2=\langle \mathbbm{e}_{k+1}, \ldots,\mathbbm{e}_m \rangle=E^{\perp}_1$ implies that $ \sigma _{\mathbbm{e}_1 , \mathbbm{e}_2}(E_2)=E_2=E^{\perp}_1$, for $k \in \{2,\ldots,m-2\}$. Thus,  the function  $f(x,y)=  x\cdot \sigma _{\mathbbm{e}_1, \mathbbm{e}_2}(y) + 1_{E_1}(x)1_{E_2}(y)$ is a bent function which belongs to $\cD$, where $x,y\in \F_2^m$. The following lemma was crucial for the design to work. 
\begin{lemma}\label{lem weight1}\cite{ZhangCPW20_DAM}
	Let $n>2$ be an integer and $f\in \mathcal{B}_n$. If 
    $|\operatorname{supp}(f)|=1$, so that $f(x)=1_{\{ \alpha\}}(x)$ for some  $\alpha\in \F_2^n$,  then we have $\deg(D_{a}D_{b}f)=n-2$ for any two different vectors $a, b\in \F_2^n \setminus \{ 0_n \}$. Furthermore, if	
	$ f(x)=\prod\limits_{j=1}^t(x_{i_j}+\alpha_{i_j}+1)$, where $\alpha_{i_j}\in \F_2, \{i_1,\ldots,i_t\}\subseteq \{1,2,\ldots,n\} $ and $2< t < n-2$, then either $\deg(D_{a}D_{b}f)= t-2$ or $ \deg(D_{a}D_{b}f)= 0$ for any two different vectors $a, b\in \F_2^n \setminus \{ 0_n \}$.
\end{lemma}
Using the above lemma, the following result from~\cite{ZhangCPW20_DAM} provides an explicit design method of bent functions in $\cD$ outside $\cM^\#$, without any initial conditions. 
\begin{theorem} \label{th:explicitclass plus}
	Let $m \geq 5$ and $k$ be  positive integers such that  $1\leq k\leq m-2$. Let $  \sigma _{\mathbbm{e}_1, \mathbbm{e}_2}(y)$ be a permutation over $\F_2^m$ given by Eq.~\eqref{eq:pistar}. 
	Define $f(x,y)=  x\cdot \sigma _{\mathbbm{e}_1 , \mathbbm{e}_2}(y) + 1_{E_1}(x)1_{E_2}(y)$, with $x,y \in \F_2^m$, 	where $E_1=\langle \mathbbm{e}_1,\ldots, \mathbbm{e}_k\rangle$
	 and $E_2=E_1^\perp$ so that $\dim(E_2)=n-k$. Then, $f$ is a bent function  outside $\mathcal{M}^\#$.
\end{theorem}
 Quite likely, other permutations (possibly of similar form) can be used for the same purpose, which was left as an open problem in \cite{ZhangCPW20_DAM}. 
 
  Indeed, in  a recent article \cite{ZhaoDLQ_2026_DAM} the authors provided a complete characterization of $\cD_k \cap \cM^\#$, where $k < m-3$, thus extending the result of Kudin and Pasalic in \cite{KudinP22_DCC}, who considered the case $\cD_0 \cap \cM^\#$. Moreover, two infinite classes of bent functions in $\cD_1$ that are provably outside $\cM^\#$ were provided in \cite{ZhaoDLQ_2026_DAM}.
 \begin{openproblem}
    Specify the exact value of the linearity index of bent functions in $\mathcal{D}$ that do not belong to $\mathcal{M}^\#$. In particular, prove or disprove that $\mathcal{D}^\#$, when $\dim(E_1)>0$, contains $\ell$-optimal bent functions.
 \end{openproblem}

\subsection{\texorpdfstring{$\mathcal{D}_0$}{D0} class of bent functions, a subclass of \texorpdfstring{$\mathcal{D}$}{D}} \label{sec:D0class}
We recall that a bent function of the form $f(x,y)=x \cdot \pi(y) + 1_{E_1}(x) 1_{E_2}(y)$, $x \in \F_2^m, y \in \F_2^m$, belongs to the class $\cD$ when $\pi(E_2)=E_1^\perp$ and $\dim(E_1) + \dim(E_2)=m$. Then, setting $E_1=0_m$ and $E_2=\F_2^m$, Carlet already observed \cite{Carlet93_EUROCRYPT} that the above condition is automatically satisfied. Moreover, Carlet derived a sufficient condition on the permutation $\pi$ so that $f$ defined above is outside $\cM^\#$. More precisely, for $m \geq 4$, it was deduced \cite[Proposition 2]{Carlet93_EUROCRYPT} that if the restriction of $\pi$ to any linear hyperplane of $\F_2^m$ is not affine, then $f \not \in \cM^\#$. 

However, the question of whether the above sufficient condition is also necessary has remained unanswered for more than two decades. In  \cite{KudinP22_DCC}, a full characterization of the $\cD_0$ class, in terms of its intersection with $\cM^\#$, was provided. This result is based on an important connection between the second-order derivatives and the algebraic degree, which we recall first. 
\begin{lemma} \label{lem:algdeg}\cite{KudinP22_DCC}
Let $g$ be a function in $\mathcal{B}_m$. If there exists an $m-k$ dimensional subspace $H$ of $\F_2^m$ such that $D_aD_b g=0$ for all $a,b \in H$, then the algebraic degree of $g$ is at most $k+1$.
\end{lemma}
Based on this lemma, the following was deduced in \cite{KudinP22_DCC} for the special case when $\deg(\pi)\geq 3$.
\begin{theorem} \label{th:DoutM}
Let $m$ be an integer, $m \geq 4$. Let $\pi$ be a permutation of $\F_2^m$ with algebraic degree $\deg( \pi) \geq 3$. Then, the function $f\colon \F_2^m \times \F_2^m \to \F_2$ defined by $f(x,y)=x \cdot \pi(y)+ \delta_0(x) \in \cD_0$ is a bent function outside $\cM^{\#}$.
\end{theorem}
Thus, it remained to provide a characterization of quadratic permutations $\pi$ such that $x \cdot \pi(y)+\delta_0(x)$ is in/outside the $\cM^{\#}$ class, which is given below.

\begin{theorem} \label{quadout} \cite{KudinP22_DCC} Let $\pi$ be a quadratic permutation of $\F_2^m$, $m \geq 4$. The function $f\colon \F_2^m \times \F_2^m \to \F_2$, defined by $f(x,y) = x \cdot \pi(y) + \delta_0(x)$, is in $\cM^{\#}$ class if and only if there is an affine hyperplane of $\F_2^m$ on which $\pi$ is affine. 
\end{theorem}

The above theorem essentially shows that Carlet's condition is not only sufficient, but also necessary. Notice that when $\deg(\pi) \geq 3$, even though the restriction of $\pi$ might be affine, by Theorem \ref{th:DoutM}, $f$ is still outside $\cM^\#$. 

The fact that the total cardinality of bent functions on $\F_2^8$ of all the classes mentioned in this section (derived from the $\cM$ class) is not even close to the total number of bent functions on $\F_2^8$ motivates the investigation of a broader class of Boolean functions, known as the generalized Maiorana-McFarland class ($\mathcal{GMM}$); it will be considered in detail in Section~\ref{sec GMM class}.
\begin{remark}
    We note that vectorial bent functions $F\colon\F_2^n\rightarrow\F_2^k$, for $2\leq k\leq m=n/2$, have been analyzed with respect to the membership of their component functions in $\cM^\#$~\cite{PasalicZKW21_DAM}. In the case when all nonzero component functions lie outside $\cM^\#$, one speaks of vectorial bent functions \emph{strongly outside} $\cM^\#$. Otherwise, if only a prescribed number of nonzero component functions lie outside $\cM^\#$, one speaks of vectorial bent functions \emph{weakly outside} $\cM^\#$. 
    
    For instance, taking $f(x,y)=\operatorname{Tr}^m_1(x \pi(y)) + \delta_0(x) \in \cD_0$, where $x,y \in \F_{2^m}$, the main idea is to define a set of permutations $\{\alpha_i \pi(y): i=1, \ldots,m \}$ where $\alpha_i$ are linearly independent elements in $\F_{2^m}$ over $\F_2$. Then, the vectorial bent function $F=(f_1, \ldots,f_{m})$ is weakly outside $\cM^\#$, where $f_i(x,y) =\operatorname{Tr}^m_1(x (\alpha_i \pi(y))) + \delta_0(x)$, see~\cite[Theorem 2]{PasalicZKW21_DAM}. In the same article, one family of vectorial bent functions $F\colon\F_2^n\to\F_2^2$, defined for doubly even dimensions $n$, that are strongly outside $\cM^\#$ is given. Some other research efforts in this direction can be found in \cite{BapicP22_DAM,KudinPCZ22_CCDS}. Similar ideas were applied in~\cite{BapicPPP24} to construct vectorial Boolean functions with the maximum number of bent components, and to construct vectorial bent-negabent functions with components outside $\mathcal{M}^\#$ in~\cite{PasalicKPP23}.
\end{remark}
 
\subsection{Superclasses of \texorpdfstring{$\mathcal{C}$}{C}, \texorpdfstring{$\mathcal{D}$}{D} and \texorpdfstring{$\mathcal{D}_0$}{D0}}\label{sec Superclasses}
Quite recently, the idea of combining the indicators that characterize the $\cC$ and $\cD$ classes was considered in \cite{BapicP22_DAM,BapicPZH22_CCDS}, for the same purpose of modifying bent functions in $\cM$ so that other (super) classes of bent functions are derived. Although the bent conditions become harder to satisfy and these classes are not rich, it is interesting that the modification of a bent function in $\cM$ can be performed on subsets (union of two subspaces) instead of on the vector subspaces. 

The following ``superclass'', denoted by $\mathcal{SC}$, was deduced in \cite{BapicP22_DAM} by combining the indicators of the $\cC$ and $\cD_0$ classes. Its formal definition was given in \cite{BapicPZH22_CCDS} as follows.
\begin{definition}
    Let $\pi$ be a permutation on $\F_{2^m}$ and let $L \subset \F_{2^m}$  be a linear subspace
of $\F_{2^m}$ such that $(\pi^{-1}, L)$ satisfies the~\eqref{eq: C property} property. Then, the class of bent functions
$f\colon\F_{2^m} \times \F_{2^m} \to \F_2$ containing all functions of the form
\begin{equation*}\label{eq:definofSC}
 f(x,y)={\rm Tr}_1^m(x\pi(y))+ 1_{L^\perp}(x) + \delta_0(x)   
\end{equation*}
is called $\mathcal{SC}$ and is a superclass of $\cD_0$ and $\cC$.
\end{definition}
\begin{theorem}\cite{BapicP22_DAM}\label{th:SCclass}
Let $\pi$ be a permutation on $\F_{2^{m}}$, $L$ be a vector subspace of $\F_{2^{m}}$ such that $(\pi^{-1},L)$ satisfies the~\eqref{eq: C property} property, $2\leq \dim(L)<m$ and $\mu \pi$ has no nonzero linear structures for all $\mu\in\F^*_{2^{m}}$. Then, the function $f\colon\F_{2^{m}}\times\F_{2^{m}}\to\F_2$ defined with $$f(x,y)={\rm Tr}_1^m(x\pi(y))+1_{S}(x),\ x,y\in\F_{2^{m}},$$
where $1_{S}(x)=1_{L^{\perp}}(x)+\delta_0(x)$, is a bent function in the class $\mathcal{SC}$ outside $\cM^{\#}$.
\end{theorem}
Especially, an explicit family of bent functions outside $\cM^\#$ was derived in \cite{BapicP22_DAM}.

\begin{corollary}
Let $n=2m$ and $s\geq 3$ be a positive divisor of $m$ such that $m/s$ is odd. Let $\pi(y)=y^d$ be a permutation on $\fF{m}$ such that $d(2^s+1)\equiv 1\ (\text{mod}\ 2^m-1)$ and $\operatorname{wt}(d)\geq 3$. Let $U=\{1,\alpha,\ldots,\alpha^{s-1}\}$, where $\alpha$ is a primitive element of $\fF{s}$. If $L=\langle K \rangle$, where $K\subseteq U$, then the Boolean function $f\colon\fF{m}\times\fF{m}\to\fF{}$ defined by $$f(x,y)={\rm Tr}_1^m(xy^d)+1_{L^{\perp}}(x)+\delta_0(x), \ x,y\in\fF{2^m},$$
is a bent function in $\mathcal{SC}$ outside $\cM^{\#}$. 
\end{corollary}
However, it is not possible to combine the typical indicators for the $\cD$ and $\cD_0$ classes. 
\begin{theorem}\label{th SD}\cite{BapicPZH22_CCDS}
Let $\pi$ be a permutation on $\fF{m}$ and let $E_1,E_2\subset\fF{m}$ be two vector subspaces of $\fF{m}$ such that $\pi(E_2)=E_1^{\perp}$ and $\dim(E_1) \geq 2$. Then, the function $f\colon\fF{m}\times\fF{m}\to\fF{}$ defined by $$f(x,y)={\rm Tr}_1^m(x\pi(y))+1_{E_1}(x)1_{E_2}(y)+\delta_0(x)$$
is not bent.
\end{theorem}

The class of functions described in Theorem \ref{th SD} was named  $\mathcal{SD}$ (superclass of $\cD$ and $\cD_0$), and similarly $\mathcal{CD}$ (superclass of $\cC$ and $\cD$) and $\mathcal{SCD}$ (superclass of $\cC,\ \cD$ and $\cD_0$) were examined in \cite{BapicPZH22_CCDS}. It turned out that, in addition to the $\mathcal{SC}$ class, the class $\mathcal{CD}$ given in Theorem \ref{th:CDclass} also contains bent functions, whilst the other classes do not.

\begin{theorem}\label{th:CDclass}\cite{BapicPZH22_CCDS}
Let $\pi$ be a permutation on $\F_2^{m}$, $L\subset\F_2^m$ be a vector subspace of $\F_2^m$ such that $(\pi^{-1},L)$ satisfies the~\eqref{eq: C property} property, and let $E_1,E_2\neq \{0_m\}$ be two vector subspaces of $\F_2^m$ such that $\pi(E_2)=E_1^{\perp}$ and $\dim(E_1\cap L^{\perp})=\dim(E_1)-1$. Then, the function  $f\colon\F_2^{m}\times\F_2^{m}\to\F_2$ defined by $$f(x,y)=x \cdot \pi(y)+1_{L^{\perp}}(x)+1_{E_1}(x)1_{E_2}(y),$$ 
is bent. 
\end{theorem}

This approach was considered in a wider framework in \cite{LiKMPTZ22_IEEE}, where the bentness of functions obtained by adding indicators of vector subspaces to a bent function was analyzed. More precisely, the authors considered 
\begin{equation*}
    h(x)=f(x) + \sum_{i=1}^{\tau}a_i 1_{V_i}(x) \;\; a_i \in \F_2,
\end{equation*}
where $f \in \mathcal{B}_n$ is bent and $V_i$ are vector subspaces of $\F_2^n$. Then, imposing the conditions on the dual $f^*$ (corresponding to the $P_\tau$ property, cf. Section \ref{sec Pi tau property}) so that $D_{v_j^{(i)}}D_{v_k^{(i)}}f=0$, where $\{v_1^{(i)}, \ldots,v_{e_i}^{(i)}\}$ is a basis of $V_i^\perp$, the bentness of $h$ could be established. In particular, when $\tau=1$, the function $h$ reduces to $h(x)=f(x)+1_V(x)$, and it was mentioned in \cite[Remark 6]{LiKMPTZ22_IEEE} that an explicit example of bent functions of this form was given in \cite[Corollary 15]{CanteautDDL06_DAM}. Moreover, the following result provides an explicit class of bent functions outside $\cM^\#$. 
\begin{corollary}\label{cor:Genericmethods2022}\cite{LiKMPTZ22_IEEE}
Let $n =2m$ and $t$ be two integers such that $m/r$ is odd, where $r = \gcd(m, t)$, and $d(2t+1) \equiv 1 \pmod{2^m-1}$ with $\operatorname{wt}(d) \geq  3$. Let $v_1,\ldots,v_\tau \in \F_{2^m}^*$ be $\tau$ ($2 \leq \tau \leq m$)
linearly independent elements such that $v_i/
v_j \in \F_{2^r}$ for any $1 \leq i < j \leq \tau$. Then, the function
\begin{equation*}
h(x, y) = {\rm Tr}_1^m(xy^d)+ \prod_{i=1}^{\tau}
(1+{\rm Tr}_1^m (v_ix)), \; \; \; x, y \in  \F_{2^m},
\end{equation*}
is a bent function outside $\cM^\#$.
\end{corollary}

We summarize this section in Table~\ref{tab:summaryCDclasses}, which covers the known families of bent functions outside $\cM^\#$ obtained by modifying bent functions in $\cM$. 

\begin{table}[htbp!]
    \centering
    \footnotesize\caption{$\cC,\cD,\cD_0$ and their superclasses of bent functions outside $\mathcal{M}^{\#}$}
    \label{tab:summaryCDclasses}
\begin{tabular}{|m{35pt}|m{90pt}|m{170pt}|}
\hline
Family & Function & Conditions and Sources \\ \hline
$\mathcal{D}_0 $ & $ x\cdot \pi(y)+\delta_0(x)$ & $\pi$ is a permutation whose restriction to any linear hyperplane is not affine (a sufficient but not necessary condition, but necessary holds for $\deg(\pi)=2$), see Theorem 6.5 and Corollary 6.6 of \cite{KudinPPZ25_JOC} \\ \cline{3-3} & & $n\geq 4,$ $\pi$ is a permutation with $\deg(\pi)\geq 3$ \\ \hline
$\mathcal{C}$ & $x\cdot \pi(y)+1_{L^{\bot}}(x)$ & Theorem 1 of \cite{ZhangPCW17_C2SI} (a general sufficient but not necessary condition, Lemma 2 of \cite{ZhangPCW17_C2SI}, Theorems 11 of \cite{ZhangCPW20_DAM}, Corollary 3.5 of \cite{KudinPCZ22_CCDS} are obtained from this condition) \\ \cline{3-3} & & see Theorem 12 and Corollary 13 of \cite{KudinP22_DCC} (contrast to Theorem 1 of \cite{ZhangPCW17_C2SI}) \\ \hline
$\mathcal{SC}$ & $x\cdot \pi(y)+1_{L^{\bot}}(x)+\delta_0(x)$  & Theorem 5 of \cite{BapicP22_DAM} (Corollary 1 of \cite{BapicP22_DAM} is a particular case) \\ \hline
$\mathcal{D}$ & $x\cdot \pi(y)+1_{E_1}(x)1_{E_2}(y)$ & Theorem 2 of \cite{ZhangPCW17_C2SI} (Examples 2, 3 of \cite{ZhangPCW17_C2SI}, Theorems 7, 13, 18, Corollaries 8, 19, and Example 2 of \cite{ZhangCPW20_DAM} are particular cases of this result) \\ \hline
$\mathcal{CD}$ & $x\cdot \pi(y)+1_{E_1}(x)1_{E_2}(y)+1_{L^{\bot}}(x)$ & Theorem 8, Propositions 2 and 3 of \cite{BapicPZH22_CCDS} \\[2mm] \hline
Modified Rothaus' function, $\mathcal{C}$ and $\mathcal{D}$ & The functions defined in Theorem 3, Corollaries 1, 2 of \cite{ZhangPWC17_IEEE}, and in Theorems 2, 3, 4, Lemma 2 of \cite{LiuZPXZ20} & Theorem 3, Corollaries 1, 2 of \cite{ZhangPWC17_IEEE}, and in Theorems 2, 3, 4, Lemma 2 of \cite{LiuZPXZ20} \\ \hline
\end{tabular}
\end{table}

\section{Bent functions with a few trace terms and the \texorpdfstring{$P_\tau$}{Ptau} property}\label{sec Further construcitons outside MM}

In this section, we first consider bent functions with a few trace terms, including monomial and binomial bent functions, as well as their relation to the $\mathcal{M}^\#$ class. Then, we consider an approach based on the so-called $P_\tau$ property, originally considered in \cite{Mesnager2014} and generalized later in \cite{Xu2017}, that aims to add further trace terms while preserving bentness. This provides a structural way to investigate the exclusion of the constructed bent functions from the $\mathcal{M}^\#$ class.
 
\subsection{Monomial, binomial, and related bent functions}\label{sec monomials}\label{sec monomial binomial etc}
In general,  monomial (trace) functions are of the form $f(x)={\rm Tr}^n_1(\alpha x^d)$, with $\alpha,x \in \F_{2^n}$, for a suitably chosen nonzero $\alpha \in \F_{2^n}$ and the exponent $d$. In particular, certain choices of exponents $d$, such as  Dillon \cite{Dillon74_PhD} and Kasami \cite{DILLON2004342,LangevinL_2008_FFAP} for instance, ensure  bentness of $f(x)={\rm Tr}^n_1(\alpha x^d)$ through a proper choice of $\alpha \in {\F}^*_{2^n}$, and furthermore these bent functions are provably outside $\cM^\#$. The \emph{Dillon's exponent $d=k(2^{n/2}-1)$} is the earliest demonstrated instance of bent functions not contained in $\cM^\#$, where the bentness is achieved if and only if $\alpha$ is a zero of the Kloosterman Sum $$K(\alpha) = \sum_{v \in \mathbb{F}_{2^{n/2}}} (-1)^{{\rm Tr}^{n/2}_1(v^{-1} + \alpha v)}.$$ Another interesting point is that the Kasami exponent gives rise to bent functions which are not weakly normal \cite{CanteautDDL06_DAM}. We also mention the pioneering work on the decomposition of bent functions \cite{CanteautC03_IEEE}, where in particular the monomial Boolean function $f(x)={\rm Tr}^n_1(\alpha x^{d})$, with even $n \geq 8$ and $15 \leq d \leq 2^{n/2}-1$, is shown to have the second-order derivatives that are not constant (thus $D_aD_b f \neq constant$ for any linearly independent $a,b \in \F_{2^n}$), for certain choices of $d$.  Bent functions of this form can be specified within the $\mathcal{PS}$ class of Dillon, thus specifying $d=2^{n/2}-1$. We summarize these exponents and the sufficient/necessary bent conditions in Table~\ref{tab:bentmonomials}. It is conjectured that the listed exponents are the only (monomial) bent exponents, and for $n \leq 24$, the conjecture is confirmed by computer simulations~\cite{LangevinL_2008_FFAP}; see also the discussion in~\cite[Section 9.1]{Tokareva15_BOOK}. Notice that the last row in Table~\ref{tab:bentmonomials} is a special case of the Dillon exponent.

\begin{table}[htbp!]\setlength{\abovecaptionskip}{0.0cm}
    \centering
	\footnotesize\caption{Known monomial bent functions in trace form 
$f(x)={\rm Tr}^n_1(\alpha x^d)$ over $\mathbb{F}_{2^n}$.}
	\label{tab:bentmonomials}	
	\centering
	\begin{tabular}
    {|m{43pt}|m{44pt}|m{115pt}|m{27pt}|m{31pt}|}
		\hline
		Exponent $d$ & Name/Type  & Conditions &  Source &  In $\mathcal{M}^{\#}$? \\
	\hline
    $2^k + 1$
    & Gold
    & $\gcd(d, 2^{n}-1) \neq 1$, \newline $\alpha \notin \{ y^{\gcd(d, 2^{n}-1)} \; \mid \; y \in \F_{2^n} \}$
    & Folklore & Yes. \\
    \hline
    $k(2^{n/2}-1)$ & Dillon & $\gcd(k, 2^{n/2}+1)=1$, \newline $\alpha \in \F_{2^{n/2}}$, 
     $ K(\alpha) =0$
    & \cite{Dillon74_PhD} & No. \\ 
    \hline
    $2^{2k} - 2^k + 1$ & Kasami & $\gcd(k,n)=1, \newline \alpha \notin \{ \; y^3 \; \mid \; y \in \F_{2^n} \; \}$ & \cite{DILLON2004342,LangevinL_2008_FFAP} &  No.\\ 
    \hline
    $(2^k+1)^2$ & Leander & $n=4k$, $k$ odd \newline 
    $\alpha \in \omega \F_{2^k}$, where $\omega \in \F_4 \setminus \F_2$ & \cite{Leander_2006_IEEE} & Yes. \\ 
    \hline
    $2^{2k} + 2^k + 1$ & Canteaut-Charpin-Kyureghyan & 
    $n=6k$, \newline $\alpha \in \F_{2^{3k}}$, \; ${\rm Tr}^{3k}_k(\alpha)=0$
    & \cite{CanteautCK08_FF} & Yes. \\
    \hline
    $2^{n/2}-1$ & Canteaut-Charpin & $\gamma^{2^m+1}=1$,\; $\gamma=\alpha^{2^m-1}$, \newline  $\alpha$ generates $\F^*_{2^n}$, \newline $\sum_{i=0}^{2^m} {\rm Tr}_1^n(\alpha \gamma^i)=2^{n/2-1}$
    & \cite{CanteautC03_IEEE} & No. \\ \hline
	\end{tabular}
\end{table}

We remark that the proofs of being outside $\cM^\#$ for the Dillon and Kasami exponents are performed using the standard second-order derivative criterion of Dillon, namely Lemma \ref{lem M-M second}, which is relatively easy compared to binomial trace functions. However, this general form of defining $f(x)={\rm Tr}^n_1(\alpha_1 x^{d_1} + \alpha_2 x^{d_2})$ was investigated in \cite{Dobbertin_Leander_etal_2006} by considering the so-called Niho exponents. A power function on $\F_{2^n}$ is called a \emph{Niho power function} \cite{Niho_PhDThesis} if its restriction to $\F_{2^{n/2}}$ is linear, where $n=2m$. This essentially means that the exponents $d_j$ satisfy the congruence $d_j \equiv 2^i \pmod{2^m-1}$, for some $i < n$. Using the normalized form such that $i=0$, these exponents are represented as $d_j=(2^m-1)s_j + 1$, for some $2 \leq s_1 \neq s_2 \leq 2^m$. It was conjectured that $d_1$ necessarily must be chosen as $d_1=\frac{1}{2}(2^m-1) +1$, where $s_1=\frac{1}{2}$ is understood modulo $2^m +1$, that is $s_1=2^{m-1}+1$. Using quite advanced algebraic techniques, the authors introduced three families of binomial bent functions of the above form by specifying suitable $d_2$. One of these families was later generalized by Leander and Kholosha \cite{Leander_2006_IEEE} into a function with $2^r$ Niho exponents. In \cite{CarletM11_Hclass}, Carlet and Mesnager, in particular, considered the bivariate representation of the second function from \cite{Dobbertin_Leander_etal_2006} and calculated its dual. They observed that this dual is not of a Niho type, replying negatively to an open question stated in \cite{Dobbertin_Leander_etal_2006}. A further analysis of binomial bent functions with Niho exponents was performed in \cite{Nihodual}, where it was shown that the bent functions introduced by Leander and Kholosha \cite{Leander_2006_IEEE} belong to $\cM^\#$. On the other hand, two families (out of three, namely Niho I and Niho III in Table~\ref{tab:bentbinomials}) from \cite{Dobbertin_Leander_etal_2006} are provably outside $\cM^\#$, see Table \ref{tab:bentbinomials}. The results on known binomial (trace) bent functions are summarized in Table~\ref{tab:bentbinomials}, where for the Niho type families the (common) bent condition $(\alpha_1 + \alpha_1^{2^{n/2}})^2=\alpha_2^{2^{n/2}+1}$ is omitted. 

In general, binomial trace bent functions that do not use Niho exponents have also been specified. As mentioned in \cite{PottPMB18},  only a few instances of these functions have been explicitly exhibited \cite{CharpinGong08,HelessethKolosha06,MesnagerFlori12}. We notice that the binomial bent functions in \cite{CharpinGong08} and \cite{MesnagerFlori12} are also \emph{hyperbent} (these are bent functions $f$ on $\F_{2^n}$ with the property that $f(x^k)$ is bent for any $k$ coprime with $2^n-1$; see \cite{YoussefG01,CarletG06}), whereas an infinite class of binomial bent functions in \cite{HelessethKolosha06} regards the ternary prime field. More precisely, the binomial trace hyperbent functions of the form $$f(x)=\operatorname{Tr}_1^n(\lambda_1x^{r_1(2^m-1)}+ \lambda_2x^{r_2(2^m-1)}),$$ with Dillon-like exponents and $\lambda_i \in \F_{2^{n/2}}$, in \cite{CharpinGong08} and \cite{MesnagerFlori12} are exactly specified for $n=6$ and $n=54$, respectively. 
Later, Mesnager~\cite{Mesnager11_IEEE,Mesnager11_Hyperbent} proposed two families of hyperbent binomials 
of the form $$f(x)=\operatorname{Tr}_1^n(ax^{r_1(2^m-1)})+ \operatorname{Tr}_1^2(b x^{(2^n-1)/3}),$$ where $b \in \F_4^*$ and $m=n/2$ is odd.
For more details on the design of hyperbent functions the reader is referred to \cite[Section 10]{Mesnager16_BOOK}. Nevertheless, in the above works the question whether the known instances of hyperbent functions belong to $\cM^\#$ was not considered. 

On the other hand, the authors in \cite{PottPMB18}  considered quadratic binomials of the form $$F_\alpha^i(x)=\operatorname{Tr}_1^n(\alpha x^{2^i}(x+x^{2^{n/2}})),$$ for which $F_\alpha^i$ was shown to be bent (belonging to $\cM$), for any $\alpha \in \F_{2^n}\setminus \F_{2^{n/2}}$ and for $i=0,\ldots, n/2-1$. The study of functions $F\colon\F_{2^n}\to\F_{2^n}$ that possess the maximum number of bent components (namely $2^n-2^{n/2}$) has recently attracted considerable attention. In particular, binomials and functions with multiple trace terms have been investigated in~\cite{MesnagerZhangDCC2019,ZhengPengDCC2020,ZhengKanDM2021,AnbaretalIEEE2022}.
\begin{table}[htbp!]\setlength{\abovecaptionskip}{0.0cm}
    \centering
	\footnotesize\caption{Niho binomial trace bent functions 
$f(x)={\rm Tr}^n_1(\alpha_1 x^{d_1} + \alpha_2 x^{d_2})$ over $\mathbb{F}_{2^n}$.}
	\label{tab:bentbinomials}	
	\centering
	\begin{tabular}
    {|m{73pt}|m{44pt}|m{91pt}|m{27pt}|m{36pt}|}
		\hline
		Exponents $d_1,d_2$ & Name/Type  & Conditions &  Source &  In $\mathcal{M}^{\#}$? \\
	\hline
    $d_1=2^{n/2} + 1$,\newline $d_2=3(2^{n/2}-1)+1$
    & Niho I
    & $n \equiv 4 \pmod{8}$ and $\alpha_2 \in \{x^5: x \in {\F}^*_{2^n}\}$ \newline
    no conditions if $n \not \equiv 4 \pmod{8}$  
    & \cite{Dobbertin_Leander_etal_2006,Nihodual} & No. \\
    \hline
    $d_1=2^{n/2}+1$, \newline $d_2=\frac{1}{4}(2^{n/2}-1)+1$ & Niho II & $n \equiv 2 \pmod{4}$ 
    & \cite{Dobbertin_Leander_etal_2006,Nihodual} & Yes. \\[3mm] \hline 
    $d_1=2^{n/2}+1$, \newline $d_2=\frac{1}{6}(2^{n/2}-1)+1$ & Niho III & $n \equiv 0 \pmod{4}$  
    & \cite{Dobbertin_Leander_etal_2006,Nihodual} & No. \\[3mm] 
    \hline
	\end{tabular}
\end{table}

We also mention two families of bent functions with multiple trace terms. Namely, the one in bivariate representation of Canteaut \emph{et al.}  $$f(x,y)={\rm Tr}_1^m(xy^s+\lambda y^{3s})+{\rm Tr}_1^m(x+\beta y){\rm Tr}_1^m(\alpha x+\alpha^{2^k-1}\beta y),$$ with $x,y \in \F_{2^m}$, see Corollary 15 of \cite{CanteautDDL06_DAM}. And, secondly, already mentioned generalization of the Niho type II exponent 
by  Leander and Kholosha \cite{Leander_2006_IEEE} which is given as $$f(x)={\rm Tr}_1^n(ax^{2^{n/2}+1}+ \sum_{i=1}^{2^{r-1}-1}x^{s_i}),$$ where $r >1$ such that $\gcd(r,m)=1$, $a \in  \F_{2^n}$ is such that 
$a + a^{2^{n/2}} = 1$, $s_i=(2^{n/2}-1)\frac{i}{2^r} + 1\pmod{2^{n/2}+1}$, $ i \in \{1,\ldots,  2^{r-1}-1 \}$. We notice that when $r=2$, this function corresponds to the bent function considered by Carlet and Mesnager~\cite{CarletM11_Hclass}, by setting there $a=b=1$. Moreover, several families of multiple trace terms bent functions were derived in \cite{LiHTK2013} from Dillon exponents, but their relationship to the class $\cM^\#$ was not considered (see also the recent related work~\cite{TuMZLJD2027}). 
We also mention a family of bent functions given by $f(x)={\rm Tr}_1^n(aG(x))$, analyzed in
\cite{Bracken_2008}, where for $n$ divisible by 6 the function $G$ on $\F_{2^n}$ is given as $$G(x)=(x + {\rm Tr}_3^n(x^{2(2^i+1)}+x^{4(2^i+1)}) + {\rm Tr}_1^n(x){\rm Tr}_3^n(x^{2^i+1}+x^{2^{2i}(2^i+1)}).$$ It was confirmed in \cite{Bracken_2008}, by computer simulations, that for $n=12$ and $i=1$ the bent function $f$ does not belong to either $\mathcal{PS}_{ap}^\#$ or to $\cM^\#$. 

Further examples of multiple trace term bent functions arise naturally from hyperbent functions. Specifically, hyperbent functions of the form  $$f_{a,b}(x)=\sum_{r \in R} \operatorname{Tr}_1^n(a_r x^{r(2^{n/2}-1)}) + \operatorname{Tr}_1^t(bx^{s(2^{n/2}-1)})$$ were considered in \cite{MesnagerFlori12}, where $R$ is a set of representatives of the cyclotomic classes modulo $2^{n/2}+1$, $a_r \in \F_{2^{n/2}}$, $s$ divides $2^{n/2} + 1$, $t$ is the size of the cyclotomic coset of $s$ modulo $2^{n/2}+1$ and $b \in \F_{2^t}$. For recent contributions to this topic we mention the work of Carlet \emph{et al.} \cite{CarletHan2022}, where  a more general characterization of the bentness of $f_{a_r}(x)=\sum_{r \in R} \operatorname{Tr}_1^n(a_r x^{r(2^{n/2}-1)})$, with $ a_r\in \F_{2^n}$ (in difference to \cite{CharpinGong08} where $a_r \in \F_{2^{n/2}}$), is obtained by utilizing M\"obius transformations. Moreover, the necessary and sufficient conditions for a class of binomial hyper-bent functions with coefficients  $a,b, c \in \F_{2^n}$ of the form $f_{a,b,c}(x)=\operatorname{Tr}_1^n((a + u_0) x^{2^{n/2}-1} + c x^{3(2^{n/2}-1)})$, was specified in 
\cite{CarletHan2022}, where $u_0 \in \F_4 \setminus \F_2$ and $m=n/2$ is odd.
Very recently, Tu \emph{et al.} \cite{Tu2025JoC} studied  the characterization of the bentness of $f_{a,b}^{d}$ with Dillon-like exponents, which is defined on $\F_{2^n}$
by $f_{a,b}^{d}(x)={\rm Tr}^n_1(a(x^{r(2^m-1)}+b)^d)$, where $a,b\in \F_{2^n}^{*}$. 
  
We also notice that vectorial hyperbent functions within the $\mathcal{PS}_{ap}$ class were considered in \cite{RibicIEEE2014}.
Nevertheless, apart from the three Niho infinite classes of binomial trace bent functions and the family of quadratic bent functions in \cite{PottPMB18} (with the above-mentioned extensions), it seems that new explicit (infinite) classes are not easily specified.

\begin{openproblem}\label{op:bentbinomials}
    Specify other infinite families of bent functions represented as binomials or multiple trace terms, especially those that do not belong to $\mathcal{M}^\#$.
\end{openproblem}
Finally, we mention the so-called \emph{permutation based bent functions} which were originally considered by Hou and Langevin \cite{HouLangevin96}. More precisely, by applying a special class of nonlinear permutations (acting non-linearly on a subset of input variables) to the input space of certain families of bent functions, EA-inequivalent functions could be derived. This approach was later extended in \cite{PasalicHZW2019} by considering larger subsets of input variables under this non-linear action. In \cite{LiLHQ23}, certain quadratic permutations $P$ were employed (acting on the whole space) to define bent functions of the form $f \circ P^{-1}$, where certain instances of these functions were shown (computationally) not to belong to $\cM^\#$.

\begin{openproblem}\label{op:permutationbent}
    Develop new infinite families of permutation based bent functions obtained via nonlinear permutations, and determine conditions under which these constructions yield EA-inequivalent bent functions that do not belong to $\mathcal{M}^\#$.
\end{openproblem}

\subsection{Bent functions as indicators of APN and AB functions}\label{sec bent from AB and APN}

Certain power monomials $F(x)=x^d$ give rise to well-analyzed combinatorial objects such as \emph{almost bent (AB)} and \emph{almost perfect nonlinear (APN)} functions; for an overview on these functions, see Carlet \cite[Section 11]{Carlet21_BOOK}. An APN function $F$ on $\F_{2^m}$ satisfies (among other characterizations) the property that $F(x)+ F(x+a)=b$ has either 0 or 2 solutions for any nonzero $a \in \F_{2^m}$ and any $b \in \F_{2^m}$. On the other hand, an AB function over $\F_{2^m}$ is defined for odd $m$ through the property that, for every nonzero $a\in\F_{2^m}$, the Walsh transform of the component function $(a\cdot F)(x)=\operatorname{Tr}_1^m(aF(x))$ satisfies
$W_{a\cdot F}(u)\in\{0,\pm2^{(m+1)/2}\}$, 
for all $u\in\F_{2^m}$. 

When $m$ is odd and $F$ is APN on $\F_{2^m}$, one can associate with $F$ the indicator Boolean function $\gamma_F(a,b)$ on $\F_{2^m}\times\F_{2^m}$, where, for any $a,b\in\F_{2^m}$, $\gamma_F(a,b)=1$ if the equation $F(x)+F(x+a)=b$ has two solutions, and $\gamma_F(a,b)=0$ otherwise. Then, it turns out that this indicator function is actually bent, which was shown in \cite{CarletCZ98_DCC}, see also \cite[Proposition 158]{Carlet21_BOOK}. The bent functions $\gamma_F$ associated with known AB functions have
been investigated in \cite{BudaghyanCH2011}. More precisely, the bent functions $\gamma_F$ associated to Kasami-Welch, Welch, and Niho functions with $m =
7, 9$, are neither in the completed $\cM$ class nor in the completed $\mathcal{PS}_{ap}$ class, see also \cite[p. 229]{Carlet21_BOOK}.

\begin{openproblem}
    Specify new families of AB functions over $\F_{2^m}$, for odd $m >3$, and determine the relationship to $\cM^\#$ of the associated indicator bent function $\gamma_F$. 
\end{openproblem}

\subsection{\texorpdfstring{$\mathcal{PS}$}{PS} class and \texorpdfstring{$\mathcal{H}$}{H} class with its variants}\label{subsec:classH} 
A \emph{partial spread} of $\F_2^n$, with $n=2m$, is a family $P$ of $m$-dimensional subspaces $S_1, S_2,\ldots, S_t \subset \F_2^n$ with pairwise trivial intersection (i.e., for all $i \neq j$ one has $S_i \cap S_j = \{0_n\}$). Further, a partial spread is a \emph{(full) spread} if the union of its subspaces equals the whole space $\F_2^n$. The main result proved by Dillon \cite{Dillon74_PhD} is that one can construct a bent function $f \colon \F_2^n \rightarrow \F_2$, with $n=2m$, from a partial spread $P$ of $\F_2^n$ by defining its support as either the union of the $t=2^{m-1}+1$ subspaces in $P$, i.e.,
$$\operatorname{supp}(f)=\bigcup_{i=1}^{2^{m-1}+1}S_i,$$
or as the union of the $t=2^{m-1}$ subspaces in $P$ excluding the zero vector, i.e.,
$$\operatorname{supp}(f)=\bigcup_{i=1}^{2^{m-1}} \left(S_i \setminus\{0_n\}\right).$$
The class of all bent functions obtained in the first case is denoted by $\mathcal{PS}^{+}$, while the class obtained in the second case is denoted by $\mathcal{PS}^{-}$. Bent functions in $\mathcal{PS}^{-}$ always attain the maximum algebraic degree $m$, while bent functions in $\mathcal{PS}^{+}$ may have lower algebraic degrees~\cite{Dillon74_PhD}. The union of $\mathcal{PS}^-$ and $\mathcal{PS}^+$ gives the whole partial spread class $\mathcal{PS}$.

 In particular, the Desarguesian spread (used by Dillon \cite{Dillon74_PhD} to specify an explicit subclass of $\mathcal{PS}^-$ called $\mathcal{PS}_{ap}$) defines bent functions $f$ as $f(x,y) =g(x/y)$,  where $g$ is a balanced Boolean function on $\F_{2^m}$ such that $g(0)=0$, with the convention that $x/y=0$ if $y=0$. The affine plane $\F_{2^{n/2}} \times \F_{2^{n/2}}$ is equal to the union of its $2^{n/2} + 1$ lines through the origin $$E_0 = \{0\}  \times  \F_{2^{n/2}}\quad \mbox{and}\quad E_a = \{(x, ax) \colon x \in  \F_{2^{n/2}}\}\quad \mbox{with } a \in \F_{2^{n/2}}.$$ Viewed as subsets of $\F_2^n$, these lines are $(n/2)$-dimensional $\F_2$-subspaces and constitute the so-called \emph{Desarguesian spread}. A univariate representation of the elements of Desarguesian spread is $$\{u \F_{2^{n/2}}^*\colon u \in  U\}, \quad \mbox{where} \quad U = \{u \in \F_{2^n}\colon u^{2^{n/2}+1} = 1\}$$  is the cyclic group of $(2^{n/2} + 1)$th roots of unity in $\F_{2^n}$.  This leads to an expression of $\mathcal{PS}_{ap}$ bent functions of the form $h(z^{2^{n/2}-1})$, $z \in  \F_{2^n}$, where $h$ is a Boolean function over $\F_{2^n}$ such that $h(0) = 0$ and whose restriction to $U$ has Hamming weight $2^{n/2-1}$. It is important that, in contrast to Niho bent functions (treated below), bent functions in this class are constant on multiplicative cosets $u {\F}_{2^m}^*$. For more details, the reader is referred to the textbooks of Carlet \cite{Carlet21_BOOK} and/or Mesnager \cite{Mesnager16_BOOK}. In our context, considering the relationship to $\cM^\#$, we mention that all bent functions in $\mathcal{PS}_{ap}$ of the form ${\rm Tr}_1^n(\alpha x^{2^{n/2}-1})$ (for suitably chosen $\alpha \in \F_{2^n}$, where $n\ge 8$) do not belong to $\cM^\#$, as shown by Dillon \cite{Dillon74_PhD}. The fact that the class $\mathcal{PS}_{ap}^\#$ is in general non-intersecting with $\cM^\#$ was also confirmed using the so-called \emph{2-rank} of Boolean functions. In the context of bent functions, this notion was originally proposed by Weng \emph{et al.} in \cite{Weng_2007_ranks} and then further studied in~\cite{Neumann2006,Weng2008,Polujan21_PhD,Kudin23_PhD}.
\begin{definition} \cite{Weng_2007_ranks} For a Boolean function $f\in\mathcal{B}_n$,  its 2-rank is defined as follows
	\begin{equation*}
		\operatorname{2-rank}(f)=\textnormal{rank}_{\F_2}(A_f),\quad\mbox{where } A_f=\left( f(x+y) \right)_{x,y\in\F_2^n}.
	\end{equation*}
\end{definition}
In~\cite{Weng_2007_ranks}, it was shown that for a Boolean function $f\in\mathcal{B}_n$ with $\deg(f)\ge2$,  the $2$-rank is an invariant under EA-equivalence.
In Table~\ref{table: ranks}, we recall the bounds about the 2-ranks of Maiorana-McFarland and Desarguesian partial spread bent functions, as stated in~\cite[Table 4]{Weng_2007_ranks}. 

\begin{table}[ht]
    \caption{Bounds on the 2-rank of some bent functions on $\F_{2}^{2m}$}
    \label{table: ranks}
	\begin{center}
		\begin{tabular}{llll}
			\hline Type & Lower bound & Upper bound & \\
			\hline All bent functions & $2 m+2$ & $\sum_{k=0}^m\binom{2m}{\min \{k, m-k\}}$ &(tight when $\left.m \leqslant 4\right)$ \\
			$f\in\mathcal{M}^\#$ & $2 m+2$ & $2^{m+1}-2$ \\
			$f\in\mathcal{PS}_{ap}^\#$ & $2^{m+1}-2$ & $\sum_{r=0}^m\binom{m}{r} 2^{\min \{r, m-r\}}$ & (tight when $\left.m \leqslant 6\right)$ \\
			\hline
		\end{tabular}
	\end{center}
\end{table}
Then, a bent function $f \in \mathcal{PS}_{ap}^\#$ in $n=2m$ variables may belong to $\cM^\#$ only if $\operatorname{2-rank}(f)=2^{m+1}-2$. This approach was recently used in \cite{KudinPPZ25_IEEE} to design bent functions that are outside $\mathcal{PS}_{ap}^\# \cup \cM^\#$. Moreover, it was shown in~\cite{KudinP22_DCC} that the probability that an $ n$-variable bent function in $\mathcal{PS}_{ap}$ class is also in $\mathcal{C}$ approaches zero as $n$ increases. Finally, we note that the recent works~\cite{AlkanAAEKT2026,AnbarKKMPP2026} investigate the exclusion of the partial spread and related classes of bent functions from $\mathcal{M}^\#$ by means of $\mathcal{M}$-subspaces (see also Section~\ref{sec M-subspaces}).

We recall that the so-called \emph{Niho Boolean functions} are functions of the form $f(\mu x)={\rm Tr}^m_1(x \phi(\mu))$, with $\mu \in U$ and $x \in \F_{2^m}^*$, that are linear on multiplicative cosets $\mu \F_{2^m}^*$. 
Here, we use polar representation of elements in $\F_{2^{2m}}$ so that any element  $z \in \F_{2^{2m}}^*$ can be represented as $z=\mu x$, where $\mu \in U$ and $$U=\{\mu \in \F_{2^{2m}}\colon \mu^{2^m+1}=1\}$$ is a multiplicative group of order $2^m+1$ (or more generally a set of coset representatives) and $\phi\colon U \to \F_{2^m}$. 
Then, the bentness of $f$ is equivalent to  the condition that for any $u \in \F_{2^{2m}}$ we have $$|\{\mu \in U\colon \phi(\mu) + {\rm Tr}_m^{2m}(\mu u)=0 \}| \in \{0,2\},$$ see Corollary 14 in \cite{Carlet21_BOOK} and \cite{CarletM11_Hclass}.

However, the use of univariate trace representation of bent functions can be traced to Dillon's thesis \cite{Dillon74_PhD}, where it was shown that all bent functions of the form ${\rm Tr}^n_1(\alpha x^{2^{n/2}-1})$, $x \in  \F_{2^n}$, are affine inequivalent to the Maiorana-McFarland functions. Moreover, Dillon introduced \emph{class $H$} which contains bent functions of the form 
\begin{equation*}
f(x,y)={\rm Tr}_1^m(y+G(yx^{2^m-2})), \quad x,y \in \F_{2^m}  
\end{equation*}
where $G$ is a permutation of $\F_{2^m}$ such that $G(x)+x \neq 0$,
and moreover $x\mapsto G(x)+bx$ is 2-to-1 for any nonzero $b \in \F_{2^m}$. Removing the linear term $y$ above and noticing that $G(x)+x\neq 0$ is not required for bentness of $f$, a  broader class can be specified as:
\begin{equation}\label{eq:classH}
g(x,y)= \left \{ \begin{array}{ll} 
{\rm Tr}_1^m(xG(\frac{y}{x})) & \textnormal{ if } x \neq 0, \\
0 & \textnormal{ if } x=0.
\end{array}
\right.
\end{equation}
Then, the class $\mathcal{H}$ defines functions of the form given by Eq. \eqref{eq:classH}, where additionally, for $x=0$, $g$ is specified as  ${\rm Tr}_1^m(\mu y)$. However, as shown in \cite{CarletM11_Hclass}, this minor modification at $x=0$ does not affect the main properties since such a function is EA-equivalent to the function $g$ given by Eq. \eqref{eq:classH}, see also \cite[p. 17]{CarletM16_DCC}.

The bijection $G$  satisfying that  $G(x)+bx$ is 2-to-1 for any nonzero $b \in \F_{2^m}$, corresponds to the so-called \emph{oval polynomials}~\cite{KoelschK25}. The known nine families of inequivalent oval polynomials are listed in Carlet's book \cite[p. 220]{Carlet21_BOOK}. Oval polynomials and $2$-to-$1$ mappings were recently considered in~\cite{ZhangCLQ2026}, where they were used to obtain further constructions of bent functions.

As noticed by Dillon \cite{Dillon74_PhD}, the class $H$ (alternatively $\mathcal{H}$) intersects with $\cM^\#$, and the question whether this class might be contained in $\cM^\#$ was left open. In \cite{Nihodual}, Budaghyan \emph{et al.} proved that Niho bent functions (being particular members of $\mathcal{H}$) of type 1) (using the notation in \cite{Nihodual}) are not included in $\cM^\#$. 

We conclude this section by mentioning a recent approach taken in \cite{XieLWZ2025}, where both additive and multiplicative cosets were employed in the design of bent functions. Moreover, it could be shown that infinite families of bent functions generated using these methods might belong to 
 the $\mathcal{PS}_{ap}$, class $\mathcal{H}$ and $\cM^\#$. On the other hand, it was also demonstrated that certain families of bent functions do not belong to any of the classes mentioned above, see Example 5 in \cite{XieLWZ2025}. Since the bent conditions in~\cite{XieLWZ2025} are not given in explicit form, their solution requires the use of computational tools such as the Magma software~\cite{BosmaCP97}.

 To summarize this section, we remark that the above-mentioned classes and their variations are not efficient in capturing the enumeration problem of bent functions since they commonly provide a few equivalence classes outside $\cM^\#$ for a fixed ambient space. In the next section, we will see that the so-called $P_\tau$ property provides more varieties in this context.  
\subsection{Bent functions and the \texorpdfstring{$\mathcal{P}_\tau$}{Ptau} property}\label{sec Pi tau property}

 The so-called \emph{$P_\tau$ property} originates in the work of  Mesnager \cite{Mesnager2014}, where both necessary and sufficient conditions for the function $$h(x)=f(x) + {\rm Tr}^n_1(ax){\rm Tr}^n_1(bx)$$ on $\F_{2^n}$ to be bent were deduced, where $f$ is a bent function. Since $x\mapsto {\rm Tr}^n_1(ax)$ describes a linear function, the above product is in general a quadratic Boolean function (unless $a=b$, $a=0$ or $b=0$). A more general form  
\begin{equation*}\label{eq:hfunctionPtau} 
h(x)=f(x) + F({\rm Tr}^n_1(u_1x), {\rm Tr}^n_1(u_2x), \ldots, {\rm Tr}^n_1(u_\tau x))
\end{equation*}
was considered in \cite[Construction 1]{Tang2017}, where $f$ is bent,  $F(X_1,X_2, \ldots,X_\tau)$ is an arbitrary polynomial in $\F_2[X_1,X_2,\ldots,X_\tau]$, $\tau$ is a positive integer such
that $1 \leq  \tau \leq m$, and $u_1, u_2,\ldots, u_\tau$ are $\tau$ nonzero elements
in $\F_{2^n}$. Then, the bentness of $h$ comes from a certain linearity condition on duals, more precisely 
$$f^*\big (x + \sum_{i=1}^{\tau}u_i w_i \big)=f^*(x)+ \sum_{i=1}^\tau w_i f_i(x), \mbox{ for all } x \in \F_{2^n}\mbox{ and all } (w_1, \ldots,w_\tau) \in \F_2^\tau,$$ where $f_i\in\mathcal{B}_n$ are Boolean functions. For instance, as observed in \cite{Tang2017}, the case $F(X_1,X_2)=X_1X_2$ corresponds to ${\rm Tr}^n_1(ax){\rm Tr}^n_1(bx)$ used in \cite{Mesnager2014}. Moreover, the case $F(X_1,X_2,X_3)=X_1X_2X_3$ covers the case considered by Xu \emph{et al.} in \cite{Xu2017}. Due to the arbitrary choice of $F$, large families of bent functions (not necessarily in the same class as the initial bent function $f$) can be obtained provided that the dual condition is satisfied.

An equivalent characterization of the bentness of $h$ on $\F_2^n$ (called $P_\tau$ property) was later stated in terms of the second-order derivatives of the dual bent function of $f$ \cite{Zheng2020}, in brief requiring that 
\begin{equation*} \label{eq:Ptaucondition}
D_{\mu_i}D_{\mu_j}f^*=0,\quad \mbox{for } 1 \leq i < j \leq \tau,
\end{equation*} 
where the elements $\mu_i \in \F_2^n$  build a vector subspace $\mathcal{U}_\tau=\langle \mu_1, \ldots,\mu_\tau\rangle $. However, the question about the class membership of the deduced families of bent functions was mostly left open. In a recent work \cite{LiKMPTZ22_IEEE}, this approach was further elaborated and the authors provided an example of a bent function outside $\cM^\#$ specified using this method.

In \cite{PasalicBZW24_DCC}, Pasalic \emph{et al.} provided explicit families of bent functions outside $\cM^\#$ within this framework, thereby partially answering the open problem from \cite{LiKMPTZ22_IEEE} (namely Open Problem 2 in \cite{LiKMPTZ22_IEEE}). Once again, the second-order derivatives played an important role, though using bivariate representation in \cite{PasalicBZW24_DCC}.

 More precisely, using  polar coordinates, every element $x\in\F_{2^n}$, $n=2m$, can be represented as $x=y+z\omega$, where $y,z\in\F_{2^m}$ and $\omega\in\F_{2^n}$ satisfies $\omega+\omega^{2^m}=1$.
Noting that any linear function from $\F_{2^ m}\times\F_{2^m}$ to $\F_2$ can be written as $(y,z)\mapsto {\rm Tr}_1^m(\alpha y+\beta z)$, where $\mu=(\alpha,\beta)\in\F_{2^m}\times\F_{2^m}$. Hence, the function $$h(x)=f(x) + \sum_{I \subseteq \{1, \ldots,\tau\}} a_I  \big (\prod_{i \in I }{\rm Tr}^n_1(\mu_i x) \big ),$$ can be written in bivariate form as 
\begin{align}\label{eq:biv_kan}
h(y,z)&=f(y,z)+F({\rm Tr}_1^m(\alpha_1y+\beta_1z),\ldots,{\rm Tr}_1^m(\alpha_{\tau}y+\beta_{\tau}z)),
\end{align}
where $f$ is bent and $F(X_1,X_2, \ldots, X_\tau)=\sum_{I \subseteq [\tau]}a_I (\prod_{i \in I}X_i)$, with $\mu_i=(\alpha_i,\beta_i)\in\F_{2^m}\times\F_{2^m}$ and $D_{\mu_i}D_{\mu_j}f^*=0$, for any $1\leq i<j\leq \tau \leq m$.
 \begin{remark}\label{rem:satisfyingPtau}
 Notice that taking a bent function $f(y,z)={\rm Tr}_1^m(y\pi(z)) \in \cM$, its dual $f^*$ is also in $\cM$ and then there always exists (at least) one $m$-dimensional subspace $V$ such that $D_aD_b f^*=0$, for all $a,b \in V$. Consequently, any subset of the basis $\{\mu_1, \ldots, \mu_\tau\}$ of $V$ can be used to specify $h$. Moreover, 
 for any bent function $f \not \in \cM^\#$ its dual is also outside $\cM^\#$ so that $\dim(V) \leq m-1$, where $V$ is such that $D_aD_b f^* = 0$ for all $a,b \in V$.
 \end{remark}
 \begin{openproblem}\label{op:Ptau_f_outsideM}
  Specify bent functions $h$ in Eq. \eqref{eq:biv_kan} using $f \not \in \cM^\#$ such that $h$ substantially differs from $f$, for instance in terms of its $\cM$-subspaces.   
 \end{openproblem}
For a Boolean function $F\in\mathcal{B}_m$ and  nonzero distinct elements $\alpha_i, \alpha_j \in  \F_{2^m }$, the authors in \cite{PasalicBZW24_DCC}
defined  $$\mathbf{\Xi}(F)=\{i:D_{\alpha_i}D_{\alpha_j}F \neq 0,\ 1\leq i<j\leq  2^m-1 \}.$$
Also, using bivariate form, the \emph{defining set} $\mathcal{U}_\tau=\{\mu_1, \ldots, \mu_\tau \}$ with $\mu_i \in \F_{2^n}$, is specified by replacing $\mu_i=(\alpha_i, \beta_i)$, where $(\alpha_i, \beta_i) \in \F_{2^m} \times \F_{2^m}$. Notice that employing $f(y,z)=\operatorname{Tr}_1^m(y\pi(z)) \in \cM$, the defining set $\mathcal{U}_\tau$ is called \emph{trivial} \cite{PasalicBZW24_DCC} if all $\alpha_i=0$ since $h$ given by Eq. \eqref{eq:biv_kan} belongs to $\cM$, see Remark 3.1 in \cite{PasalicBZW24_DCC}.
 
Using this notation, the following result was deduced in \cite{PasalicBZW24_DCC}.
\begin{theorem}\label{th:partialsol}
Let $f$ be a bent function on $\fF{n}$, $n=2m$, defined with $f(y,z)={\rm Tr}_1^m(y\pi(z))$, with $y,z \in \F_{2^m}$, and let $\mu_i=(\alpha_i,0)\in\fF{m}\times\fF{m}$ for $1\leq i \leq \tau \leq m$ be such that $D_{\mu_i}D_{\mu_j}f^*=0$ for any $i\neq j$, thus constituting a non-trivial defining set for $f^*$. Let $h$ be defined by Eq. \eqref{eq:biv_kan}. If 
\begin{enumerate}
\item ${\rm Tr}_1^m(\lambda\pi)$ has no nonzero linear structures for any $\lambda\in\fF{m}^*$,
\item  $|\mathbf{\Xi}(F)|\geq 2$, 
\end{enumerate}
 then $h$ is a bent function on $\F_{2^n}$ outside $\cM^{\#}$.
\end{theorem}
Most notably, Example 3.2 in \cite{PasalicBZW24_DCC} demonstrates the fact that $\mathcal{P}_{\tau}\not\subset (\mathcal{C}\cup \mathcal{SC}\cup \mathcal{D}_0\cup\mathcal{D}\cup\mathcal{CD})$ and also $\mathcal{P}_{\tau}\not\subset  \cM^\#$, where $P_\tau$ denotes the family of bent functions specified by Eq. \eqref{eq:biv_kan}. On the other hand, it is important to emphasize that the only completed class considered in these relationships is $\cM^\#$. 

 \begin{openproblem}
    Analyze the inclusion of bent functions obtained using the $P_\tau$ property in the completed (super) classes.
 \end{openproblem}

\begin{remark}
    The authors in \cite{PasalicBZW24_DCC} also extended the notion of ``non-trivial defining sets'' and considered an alternative form for $F$ given by 
    \begin{align*}\label{eqF}
    &F(X_1,X_2,\ldots,X_{\tau})\\
=&G_1({\rm Tr}_1^m(\alpha_1y),\ldots,{\rm Tr}_1^m(\alpha_ky))G_2({\rm Tr}_1^m(\beta_{k+1}z),\ldots,{\rm Tr}_1^m(\beta_{\tau}z)) \\
+&G_3({\rm Tr}_1^m(\alpha_1y),\ldots,{\rm Tr}_1^m(\alpha_ky)) ,
\end{align*}
and with respect to $f(y,z)={\rm Tr}_1^m(y\pi(z))$ and its dual  $f^*(y,z)={\rm Tr}_1^m(z\pi^{-1}(y))$,  one could associate two sets 
\begin{equation*}
    \begin{split}
        A&=\{a_i=(\alpha_i,0):1\leq i \leq k,\alpha_i\in\F_{2^m}\}\quad\mbox{and} \\
        B&=\{b_j=(0,\beta_j):k+1\leq j \leq \tau,\beta_j\in\F_{2^m}\},
    \end{split}
\end{equation*}
such that $D_{a_i}D_{a_{i'}}f^*=D_{b_{j}}D_{b_{j'}}f^*=0$ for any $a_i,a_{i'}\in A,b_{j},b_{j'}\in B$, where $i\neq i'$ and $j \neq j'$. Then, the bentness of $h$ defined by Eq. \eqref{eq:biv_kan}, using $F$ above, is given in Lemma 4.1 in \cite{PasalicBZW24_DCC}. 
\end{remark}

Whereas the methods discussed in the previous section are in general not so efficient (again using the condition that the components of a permutation $\pi$ do not admit linear structures), within the same framework another approach was recently proposed in \cite{LiKMPZ25_IEEE}. Essentially, the main difference compared to the function $h$ given by Eq. \eqref{eq:biv_kan}, which uses linear functions (given by ${\rm Tr}_1^m(\alpha_iy+\beta_iz)$ using bivariate representation) as input variables of $F$, is that the function $h$ in Construction 1 in \cite{LiKMPZ25_IEEE} employs the derivatives of $f$ instead of ${\rm Tr}_1^m(\alpha_iy+\beta_iz)$. More precisely, the following construction describes this approach and applies to arbitrary Boolean functions, thus not only bent. 
\begin{cons}\cite{LiKMPZ25_IEEE}\label{cons:CostructionPtau}
 Let $f$ be a Boolean function on $\F_2^n$. Let $\mu_1, \mu_2,\ldots, \mu_r \in  \F_2^n$ be $r$ linearly
independent vectors over $\F_2$ such that $D_{\mu_i}D_{\mu_j} f = 0$ for any $1 \leq  i < j \leq  r$. Then, for any
Boolean function $F$ on $\F_2^r$, define a Boolean function as:
\begin{equation*}\label{eq:functionh2025}
h(x) = f(x) + F(f_1(x), f_2(x),\ldots, f_r(x)),
\end{equation*}
where $f_i(x) = f(x) + f(x + \mu_i)$ for each 
$i \in \{1, 2,\ldots, r\}$. 
\end{cons}
The main connection between $h$ and $f$ is that,  using the univariate representation in Theorem 1 in \cite{LiKMPZ25_IEEE},  $$W_h(\beta)=(-1)^{F({\rm Tr}_1^n(\mu_1\beta), \ldots, {\rm Tr}_1^n(\mu_r \beta))}W_f(\beta),$$ for any $\beta \in \F_{2^n}$ and linearly independent $\mu_1, \ldots, \mu_r \in \F_{2^n}$,  
 which proves that $f$ and $h$ have the same extended Walsh spectrum, i.e., the same distribution of the absolute values of their Walsh transforms. This result was frequently employed in \cite{LiKMPZ25_IEEE} for deriving a few infinite classes of bent functions outside $\cM^\#$.

In Table \ref{tab:summaryPtau}, we summarize the known families (instances) of bent functions outside $\cM^\#$ that were obtained using this $P_\tau$ property. 
\begin{table}[ht!]\setlength{\abovecaptionskip}{0.0cm}
    \centering
    \begin{threeparttable}
	\footnotesize\caption{Known bent functions outside the class $\mathcal{M}^{\#}$ using $P_\tau$ property}
    \label{tab:summaryPtau}	    
	\begin{tabular}{|m{28pt}|m{125pt}|m{90pt}|m{45pt}|}
		\hline
		Family & Function  & Conditions &  Source \\
	\hline

 Pasalic \emph{et al.} & {${\rm Tr}_1^m(y\pi(z))+F({\rm Tr}_1^m(\alpha_1y+\beta_1z),\ldots,{\rm Tr}_1^m(\alpha_{\tau}y+\beta_{\tau}z))$} & {Theorems 3.2 and 4.1 of \cite{PasalicBZW24_DCC}}  & \cite{PasalicBZW24_DCC}\\

\hline
  Li~\emph{et al.} & {\small${\rm Tr}_1^n\big(\alpha(x+{\rm Tr}_k^n(x^{2^k+1}))^3\big)$}  & $n=4k\geq8$ and ${\rm Tr}_k^n(\alpha)\neq0$ & \cite{LiLHQ23}$^{1}$
\\
\hline
  & $f(x,y)+F\circ\phi(y)$ & Theorem 3 of \cite{LiKMPZ25_IEEE} & \cite{LiKMPZ25_IEEE}
\\
\cline{2-4}
 & $x\cdot z^{(r)}+p(z^{(r)})+q(z^{(\bar{r})})$ & $q$ is bent outside $\mathcal{M}^{\#}$ & \cite[Th. 8]{LiKMPZ25_IEEE}
\\
\cline{2-4}
 Li~\emph{et al.}  & {\small$x\cdot z^{(r)}\!+\!p_1(z^{(r)})\!+\!q_1(z^{(\bar{r})})\!+\!(p_1\!+\!p_2)(z^{(r)})(q_1\!+\!q_2)(z^{(\bar{r})})$} & $\deg(p_1+p_2)=r\geq3$, $q_1$, $q_2$ are bent with $q_1$ outside $\mathcal{M}^{\#}$ &\cite[Th. 9]{LiKMPZ25_IEEE}
\\
\cline{2-4}
  &{\small $x\cdot z^{(r)}+p_1(z^{(r)})+q_1(z^{(\bar{r})})+[(p_1+p_2)(z^{(r)})+q_2(z^{(\bar{r})})]\big(\alpha\cdot z^{(\bar{r})}\big)$} & $\deg(p_1+p_2)=r\geq3$, $q_2(z^{(\bar{r})})=\beta\cdot z^{(\bar{r})}$, $q_1$ and $q_1(z^{(\bar{r})})+\big(\alpha\cdot z^{(\bar{r})}\big)\big(\beta\cdot z^{(\bar{r})}\big)$ are bent with the latter one outside $\mathcal{M}^{\#}$ & \cite[Th. 10]{LiKMPZ25_IEEE}
\\
\cline{2-4}
  & $x\cdot \phi(y)+G(y)$ & Theorem 12 of \cite{LiKMPZ25_IEEE} &\cite{LiKMPZ25_IEEE}
\\
\hline

	\end{tabular}

\begin{tablenotes}
    \footnotesize
    \item[$^{1}$] In \cite{LiLHQ23}, bent functions of the form ${\rm Tr}_1^n(F\circ P(x))$ outside $\mathcal{M}^{\#}$ were found using Magma~\cite{BosmaCP97}.
  \end{tablenotes}
\end{threeparttable}
\end{table}

 \begin{openproblem}
    Analyze the inclusion of bent functions constructed using $P_\tau$ property in the classes $\mathcal{C}^\#$ and $\mathcal{D}^\#$. Show that bent functions constructed using the $P_\tau$ property can be $\ell$-optimal.
 \end{openproblem}

We remark that this approach was extended to vectorial Boolean functions in \cite{BapicP_DCC2021} and later used in \cite{BapicPPP24} in the design of vectorial bent functions whose certain number of components are outside $\cM^\#$, as originally considered in \cite{PasalicZKW21_DAM}. In this context, the following research task is quite interesting. 
\begin{openproblem}\label{op:stronglyoutside}
    Use the $P_\tau$ property to define vectorial bent functions that are strongly outside $\cM^\#$.
\end{openproblem}

Although the contributions summarized in this section also imply that there is a large variety of approaches to specifying bent functions outside $\mathcal{M}^\#$, in most cases, the presented families are quite sparse in each fixed dimension. In the following section, we focus on a broader class of Boolean functions, known as the generalized Maiorana-McFarland class, which in many cases gives rise to large families of bent functions provably outside $\mathcal{M}^\#$.

\section{The completed generalized Maiorana-McFarland class}\label{sec GMM class}

The \emph{generalized Maiorana-McFarland class} was first introduced in~\cite{Camion1991}, motivated by the study of resilient Boolean functions, and has since played a central role in the construction of cryptographically significant functions~\cite{Camion1991,Carlet_Indirect,Carlet02_CRYPTO,WGZhang2014} (see also \cite[p. 354]{Carlet10_CHAPTER_BOOLEAN}).

Members of the generalized Maiorana-McFarland class, denoted by $\mathcal{M}_{r,s}$ are exactly the Boolean functions $f \in \mathcal{B}_n$ that are concatenations of $2^s$ affine functions in $r$ variables.  

\subsection{Definitions, basic properties, and the notion of an \texorpdfstring{$\mathcal{M}$-subspace}{M-subspace}}\label{sec GMM class intro}\label{sec GMM basics}
Formally, the generalized Maiorana-McFarland class is defined as follows.

\begin{definition}\label{def:GMMclass}
	Let $r$ and $s$  be  positive integers and $n=r+s$. The set of all Boolean functions $f_{\phi,h} \colon \F_2^{r} \times \F_2^{s} \rightarrow \F_2$
    of the form
	\begin{eqnarray}\label{eq:Mextend}
		f_{\phi,h}(x,y)=x\cdot\phi(y) +
		h(y), \;\; x \in \F_2^{r}, y \in \F_2^{s},
	\end{eqnarray}
	is called the \emph{generalized Maiorana-McFarland class} $\mathcal{M}_{r,s}$,
    where $\phi
	\colon \F_2^{s} \rightarrow \F_2^{r}$  and  $h\in
	\mathcal{B}_{s}$.
    When $n$ is even and $f$ is supposed to be bent, we use the notation $\mathcal{GMM}_{\frac{n}{2}+k}$, where $r=\frac{n}{2}-k$ and $s=\frac{n}{2}+k$, and $k=0,1,\ldots,\frac{n}{2}-1$.
\end{definition}
First, we note that a Boolean function $f_{\phi,h} \colon \F_2^{r} \times \F_2^{s} \rightarrow \F_2$, defined by $f_{\phi,h}(x,y) = x \cdot \phi(y) + h(y)$, where $x \in \F_2^{r}$ and $y \in \F_2^{s}$, is not bent for $r > s$. This fact justifies restricting $k$ to the range $0 \le k \le n/2 - 1$ in the definition above. Using simple arguments (see, e.g., \cite{Carlet21_BOOK,Polujan21_PhD,ZhaoWZPC22_DCC}) one can show that if a Boolean function $f \in \mathcal{B}_n$ belongs to $\mathcal{M}_{r,s}$ (resp. $\mathcal{GMM}_{\frac{n}{2}+k}$), then $f$ belongs to $\mathcal{M}_{r-1,s+1}$ (resp.  $\mathcal{GMM}_{\frac{n}{2}+k+1}$) too. Consequently, these classes form nested families

\[
\begin{array}{ccccccccc}
\mathcal{M}_{n/2,n/2} & \subset & \mathcal{M}_{n/2-1,n/2+1} & \subset & \cdots & \subset & \mathcal{M}_{1,n-1} \\
(\mbox{resp. } \mathcal{GMM}_{n/2} & \subset & \mathcal{GMM}_{n/2+1} & \subset & \cdots & \subset & \mathcal{GMM}_{n-1})
\end{array}.
\]

In the extreme ``left'' case $k=0$, the class defined by Eq.~\eqref{eq:Mextend} corresponds to the $\mathcal{M}$
class of bent functions if and only if $\phi$ permutes $\F_2^{\frac{n}{2}}$. In the extreme ``right'' case  $k = n/2 - 1$, the class defined by Eq.~\eqref{eq:Mextend} includes all Boolean functions on $\F_2^{n}$. Similarly to the class $\mathcal{M}$, the completed version of $\mathcal{M}_{r,s}$ (resp. $\mathcal{GMM}_{n/2+k}$) is denoted by $\mathcal{M}_{r,s}^\#$ (resp. $\mathcal{GMM}_{n/2+k}^{\#}$) and called the \emph{completed generalized Maiorana-McFarland class}. It is obtained by applying all invertible affine transformations that preserve extended-affine equivalence to the functions in the class $\mathcal{M}_{r,s}$ (resp. $\mathcal{GMM}_{n/2+k}$).  

The completed Maiorana-McFarland class $\mathcal{M}^\#$ of bent functions has been characterized in~\cite[p. 10]{Dillon74_PhD} and~\cite[Lemma 33]{CanteautDDL06_DAM}. A similar proof holds for the $\mathcal{M}_{r,s}^\#$ class as well (it can be found, e.g., in Proposition 54 of~\cite{Carlet21_BOOK} and in~\cite{Polujan21_PhD}). Also note that this characterization remains valid when $n$ is odd, although our main focus will continue to be on the even case.

\begin{theorem}\label{theorem: MM} Let $f\in\mathcal{B}_n$ with $n=r+s$. The following statements are equivalent:
	\begin{enumerate}
		\item The function $f$ belongs to the $\MrsClass$ class.
		\item There exists a vector subspace $U$ of dimension $r$ such that the second-order derivatives $D_a D_b f$ vanish for all $a ,b  \in U$, that means $D_a D_b f=0$.
		\item\label{proposition: MM part (iii)} There exists a vector subspace $U$ of dimension $r$ such that the function $f$ is affine on every coset of $U$.
	\end{enumerate}
\end{theorem}

\begin{remark}\label{rem: GMM iff remarks}
    The function $f_{\phi,h}\in\mathcal{B}_n$ defined by Eq.~\eqref{eq:Mextend} is referred to as a \emph{Maiorana-McFarland representation} of a Boolean function $f \in \mathcal{B}_n$ (also known in some literature, e.g.,~\cite{LogachevSY12_BOOK}, as a \emph{linear branching}). In fact, Theorem~\ref{theorem: MM} guarantees there always exists an invertible linear transformation $A$ over $\F_2^n$ (where $A\in\operatorname{GL}(n,2)$) such that $f(zA) = f_{\phi,h}(x, y)$ for suitable $\phi\colon\F_2^s\to\F_2^r$ and $h\in\mathcal{B}_s$.
\end{remark}

In order to introduce an algorithm that can be practically applied to construct a Maiorana-McFarland representation of a Boolean function, we recall the definition of an $\mathcal{M}$-subspace, which is based on the second condition of Proposition~\ref{eq:Mextend}.

\begin{definition}\cite{PolujanP20_DCC}\label{definition: MSubspace}
	A vector subspace $U$ of $\F_2^n$ is called an \emph{$\mathcal{M}$-subspace} of a Boolean function $f\in\mathcal{B}_n$, if for all $a ,b  \in U$ the second-order derivatives $D_a D_b f$ are constant zero functions, i.e., $D_a D_b f=0$. We denote by $\mathcal{MS}_r(f)$ the collection of all $r$-dimensional $\mathcal{M}$-subspaces of $f$ and by $\mathcal{MS}(f)$ the collection 
	$$\mathcal{MS}(f)=\bigcup\limits_{r=1}^{n} \mathcal{MS}_r(f).$$
\end{definition}

We also recall the following EA-invariant, known as the \emph{linearity index}~\cite[p.~82]{Yashchenko97_PPI}, which represents the maximum number of variables in a Maiorana-McFarland representation (given by Eq.~\eqref{eq:Mextend}) of a Boolean function.

\begin{definition}\label{definition: Linearity Index}
	The \emph{linearity index} $\ind(f)$ of a Boolean function $f\in\mathcal{B}_n$ is the maximal possible $r$, such that $f\in\MrsClass$. In terms of $\mathcal{M}$-subspaces, the linearity index of $f$ is given by $\ind(f)=\max\limits_{U\in \mathcal{MS}(f)}\dim(U)$.
\end{definition}

\begin{remark}
    For an arbitrary Boolean function $f \in \mathcal{B}_n$, its linearity index satisfies $1 \le \ind(f) \le n$. Functions attaining the lower bound, $\ind(f) = 1$, are precisely those for which $D_a D_b f \ne 0$ for all linearly independent $a, b \in \mathbb{F}_2^n$. Conversely, functions attaining the upper bound, $\ind(f) = n$, are exactly the affine functions. It is worth noting that the linearity index of a bent function in $n$ variables cannot exceed $n/2$.
\end{remark}

The following simple algorithm (considered in~\cite[Algorithm 1]{PolujanP20_DCC} as a variation of~\cite[Algorithm 2]{CanteautDDL06_DAM}) constructs all $\mathcal{M}$-subspaces of a given dimension $r$ for a Boolean function $f$. It starts with all two-dimensional $\mathcal{M}$-subspaces (provided they exist) and iteratively extends them to higher dimensional ones until reaching dimension $r$ (if possible). The maximum value of $r$ for which $\mathcal{MS}_r(f)\ne\varnothing$ holds is exactly the linearity index $\ind(f)$ of $f\in\mathcal{B}_n$. For recent applications of this algorithm (and its variations) to the analysis of bent functions outside $\mathcal{M}^\#$, we refer to~\cite{QuZL25,LiLHQ23,QianSC26,XieLWZ2025}.

\begin{algorithm}[H]
	\caption{Construct the collection $\mathcal{MS}_r(f)$.}
	\label{algorithm: f in M}
	\begin{algorithmic}[1]
		\REQUIRE A Boolean function $f\colon\F_2^{n}\rightarrow\F_2$ and $2 \le r \le n$.
		\ENSURE The collection $\mathcal{MS}_r(f)$.	
		\STATE \textbf{Construct} $\mathcal{MS}_2(f)=\{ \langle a ,b   \rangle: \dim(U)=2 \mbox{ and }D_a D_b f=0 \}$.
		\FORALL{subspaces $U\in \mathcal{MS}_2(f)$}
		\REPEAT  
		\STATE \textbf{Determine} subspaces $\tilde{U}=\langle U,\tilde{u} \rangle$, for all $\tilde{u}\notin U$, such that for any two-dimensional vector subspace $ \langle a ,b  \rangle\subseteq \tilde U$ the second-order derivative satisfies $D_a D_b f=0$. \ \\[1mm]
		\STATE \textbf{Put} $U\gets \tilde{U}$ for the obtained subspaces $\tilde{U}$. 		
		\UNTIL{$\dim(U)=r$.}		
		\STATE \textbf{Output} subspaces $U$ of dimension $r$.
		\ENDFOR
	\end{algorithmic}
\end{algorithm}

\begin{remark}\label{remark: How to construct a linear mapping}
	In many cases, a bent function $f \in \mathcal{B}_n$ is not directly expressed in a Maiorana-McFarland representation as defined in Eq.~\eqref{eq:Mextend}. To obtain such a representation, one can proceed as follows. Given an $\mathcal{M}$-subspace $U \in \mathcal{MS}_r(f)$ of $f\in\mathcal{B}_n$, the goal is to construct an invertible matrix $A_U$ that transforms $f$ into its Maiorana-McFarland representation. In other words, we seek a linear change of variables such that $f(z A_U) = f_{\phi,h}(x,y) = x \cdot \phi(y) + h(y)$, where $z \in \mathbb{F}_2^n$, $x \in \mathbb{F}_2^r$, and $y \in \mathbb{F}_2^s$.

    The main idea is that the values of the function $f_{\phi,h}(x,y)=x \cdot \phi(y) + h(y)$ on each coset $\mathbb{F}_2^r + y$ for $y \in \mathbb{F}_2^s$ coincide with the values of $f$ on the corresponding coset $U + \bar{u}$ for $\bar{u} \in \bar{U}$. Consequently, this correspondence enables the construction of $A_U$ using the change-of-basis matrix
    \begin{equation}\label{equation: Linear Transform}
    		A_U=
    		\left(
    		\begin{array}{c}
    			\gjb(U) \\ \hline
    			\gjb(\bar{U}) \rule{0pt}{1.02\normalbaselineskip}
    		\end{array}
    		\right),
    	\end{equation}
    where $\gjb(U)$ denotes the \emph{Gauss-Jordan basis} (row echelon form) of the subspace $U$. Here, $\bar{U}$ represents a \emph{complement} of $U$, satisfying $\dim(U) + \dim(\bar{U}) = n$ and $U \cap \bar{U} = \{0_n\}$. It can be uniquely determined using the procedure described in~\cite[Subsection 4]{CanteautDDL06_DAM}. In this way, $A_U$ provides the necessary linear transformation that brings $f$ into a Maiorana-McFarland representation. 
\end{remark}

\begin{example}\label{ex Almost MM cov EBS}
In this example, we consider a Boolean bent function in $n = 8$ variables, constructed in~\cite[Theorem 7]{Clikeman_2019} using a technique called \emph{covering extended building sets (EBS)}. This approach adapts a powerful construction method from difference sets~\cite{DavisJedwab97} to the context of Boolean bent functions and, notably, generalizes the Maiorana-McFarland construction, as discussed in~\cite{Clikeman_2019}. The algebraic normal form of this function $f \in \mathcal{B}_8$ is given as follows:
    \begin{equation*}
        \begin{split}
            f(z)&=z_1 + z_2 + z_1 z_2 + z_3 + z_2 z_3 + z_4 + z_1 z_4 + z_2 z_4 + z_1 z_2 z_4 + z_3 z_4 \\
            &+z_1 z_2 z_3 z_4 + z_1 z_5 + z_2 z_6 + z_1 z_3 z_7 + z_2 z_3 z_7 + z_1 z_2 z_3 z_7 + z_4 z_7 \\
            & + z_1 z_4 z_7 + z_3 z_8 + z_2 z_3 z_8 + z_1 z_4 z_8 + z_2 z_7 z_8 + z_1 z_2 z_7 z_8.
        \end{split}
    \end{equation*}

In~\cite[Theorem 7]{Clikeman_2019}, using an analysis of $\mathcal{M}$-subspaces, the authors showed that the considered function $f \in \mathcal{B}_8$ lies outside the completed Maiorana-McFarland class $\mathcal{M}^\#$. Using Algorithm~\ref{algorithm: f in M}, one can verify that $|\mathcal{MS}_2(f)| = 187$, $|\mathcal{MS}_3(f)| = 15$, and $|\mathcal{MS}_4(f)| = 0$. Hence, $f \notin \mathcal{M}^\#$, but $f \in \mathcal{M}_{3,5}^\#$ since $\ind(f) = 3$. We now construct a Maiorana-McFarland representation $f_{\phi,h}$ of $f$ using the following 3-dimensional $\mathcal{M}$-subspace $U$ of $f$, given together with its complement $\bar{U}$, both expressed in the Gauss-Jordan basis form as follows:
\begin{equation*}
    \gjb(U) =\left(
\begin{array}{cccccccc}
 0 & 0 & 0 & 0 & 1 & 0 & 0 & 0 \\
 0 & 0 & 0 & 0 & 0 & 1 & 0 & 0 \\
 0 & 0 & 0 & 0 & 0 & 0 & 1 & 0 \\
\end{array}
\right)\quad \mbox{and} \quad \gjb(\bar{U}) =\left(
\begin{array}{cccccccc}
 1 & 0 & 0 & 0 & 0 & 0 & 0 & 0 \\
 0 & 1 & 0 & 0 & 0 & 0 & 0 & 0 \\
 0 & 0 & 1 & 0 & 0 & 0 & 0 & 0 \\
 0 & 0 & 0 & 1 & 0 & 0 & 0 & 0 \\
 0 & 0 & 0 & 0 & 0 & 0 & 0 & 1 \\
\end{array}
\right).
\end{equation*}
Constructing the linear mapping $A_U$ as in Eq.~\eqref{equation: Linear Transform}, one can see that Eq.~\eqref{eq:Mextend} $$f(z A_U) = f_{\phi,h}(x,y) = x \cdot \phi(y) + h(y)$$ holds for all $z\in\F_2^8,x\in\F_2^3,y\in\F_2^5$ satisfying $zA_U=(x,y)$, where the mapping $\phi\colon\F_2^5\to\F_2^3$ and the function $h\in\mathcal{B}_6$ are defined by
\begin{equation*}\label{eq MM representation example}
    \begin{split}
        \phi(y)&=(y_1, y_2, y_1 y_3 + y_2 y_3 + y_1 y_2 y_3 + y_4 + y_1 y_4 + y_2 y_5 + y_1 y_2 y_5), \\
        h(y)&= y_1 + y_2 + y_1 y_2 + y_3 + y_2 y_3 + y_4 + y_1 y_4 + y_2 y_4 + y_1 y_2 y_4 \\
        &+ y_3 y_4 + y_1 y_2 y_3 y_4 + y_3 y_5 + y_2 y_3 y_5 + y_1 y_4 y_5.\\
    \end{split}
\end{equation*}
\end{example}

\begin{remark}
    The mapping $\phi$ and the Boolean function $h$ constructed in the previous example possess additional algebraic properties, as they emerge from the Maiorana-McFarland decomposition of a bent function in $n=8$ variables with an almost maximal linearity index, $n/2 - 1 = 3$. Bent functions of this type are known as almost Maiorana-McFarland functions. Their characterization, along with the properties of the associated mappings $\phi$ and $h$, will be considered in detail in Section~\ref{sec Almost MM}.
\end{remark}

\begin{openproblem}
    Provide explicit constructions of bent functions that are provably outside $\mathcal{M}^\#$ on $\F_2^n$ (with a prescribed value of the linearity index) using the EBS construction.
\end{openproblem}

In the following sections, we focus on bent functions with a fixed linearity index. In many situations, the notion of a relaxed linearity index (another EA-invariant introduced in~\cite{PolujanP20_DCC}) proves useful for establishing upper bounds on the linearity index of bent functions, and consequently, for determining their classification with respect to the generalized Maiorana-McFarland class.

\begin{definition}\cite{PolujanP20_DCC} A vector subspace $U$ of $\F_2^n$ is called a \emph{relaxed} $\mathcal{M}$\emph{-subspace} of a Boolean function $f\in\mathcal{B}_n$, if for all $a ,b  \in U$ the second-order derivatives $D_a D_b f$ are either constant zero or constant one functions, i.e., $D_a D_b f=0$ or $D_a D_b f=1$.  We denote by $\mathcal{RMS}_r(f)$ the collection of all $r$-dimensional relaxed $\mathcal{M}$-subspaces of $f$ and by $\mathcal{RMS}(f)$ the collection $$\mathcal{RMS}(f)=\bigcup\limits_{r=1}^{n} \mathcal{RMS}_r(f).$$
Similarly, a \emph{relaxed linearity index} $\rind(f)$ of $f\in\mathcal{B}_n$ is defined by  $\rind(f)=\max\limits_{U\in \mathcal{RMS}(f)}\dim(U)$.
\end{definition}

\begin{remark}
Another important concept closely related to the notion of an \emph{$\mathcal{M}$-subspace} is that of a \emph{normality subspace}~\cite[p.~30]{Leander2005}, which arises from the broader concept of \emph{normality}. Formally, a Boolean function $f \in \mathcal{B}_n$ is called \emph{$k$-normal} if it is constant on some $k$-dimensional affine subspace of $\mathbb{F}_2^n$, and \emph{weakly $k$-normal} if it is affine on such a subspace (sometimes referred to as a \emph{normality subspace}). A function is said to be \emph{normal} (respectively, \emph{weakly normal}) if it is $\lceil n/2 \rceil$-normal (respectively, weakly $\lceil n/2 \rceil$-normal).

The notion of normality was first introduced by Dobbertin~\cite{Dobbertin94_FSE} for even numbers of variables and later extended by Charpin~\cite{Charpin04}. Efficient algorithms for testing the normality of Boolean functions were presented in~\cite{CanteautDDL06_DAM}. A variation of such algorithms was used in~\cite{GillotLP_2025} to demonstrate that all bent functions in eight variables are either normal or weakly normal.

While normality and weak normality describe Boolean functions whose restrictions are constant or affine on affine subspaces, several recent works~\cite{STABILITY,HAUGLAND,GillotLP_2025} extend these notions further by investigating the algebraic degree of such restrictions.
\end{remark}

\subsection{Characterizing bentness in the generalized Maiorana-McFarland class} \label{sec Characterizing bentness in GMM}

In contrast to the standard Maiorana-McFarland class $\mathcal{M}$ of bent functions $f_{\phi,h}$ defined in Eq.~\eqref{eq:Mextend} for $r=s=n/2$, where the function $h\in \mathcal{B}_{n/2}$ is arbitrary, the bentness condition imposes additional restrictions on $h \in \mathcal{B}_{n/2 + k}$, for $k>0$. In this subsection, we summarize the known characterizations of bent functions within this generalized framework, which are referred to as bent functions in the \emph{generalized Maiorana-McFarland class}. To do so, we first recall the notation originally introduced in~\cite[Section~6.3]{LogachevSY2004}; here we refer to the corresponding parts of the book~\cite{LogachevSY12_BOOK}, which is an English translation of~\cite{LogachevSY2004}.
    
    \begin{definition}
        A set $C \subseteq \F_2^n$ is called a \emph{bent set} if the following two conditions are satisfied:  
    \begin{enumerate}
        \item $|C| = 2^l$ for some integer $l$.  
        \item For each $v \in \F_2^n \setminus \{0_n\}$, the set $C_v = \{y \in C \colon y + v \in C\}$ is either empty or has even cardinality.
    \end{enumerate}
    \end{definition}

\begin{definition}\cite[Definition 6.39]{LogachevSY12_BOOK}
	Let $h\in\mathcal{B}_n$ be a Boolean function and let $C\subseteq\F_2^n$ be a bent set. The pair $(h, C)$ is called a \emph{partially defined bent function} if, for every $v \in \mathbb{F}_2^n \setminus \{0_n\}$, the Hamming weight of the restriction of the derivative $D_v h$ to the set $C_v$ is equal to $\tfrac{1}{2}|C_v|$.
\end{definition}

\begin{remark}\cite[p. 173]{LogachevSY12_BOOK}
    1. Any one-element set $C$ is a bent set since $C_v =\varnothing$ for $v \in \F_2^n \setminus\{0_n\}$. 
    
    \noindent 2. If $C$ is a coset of a subspace $L$ of even dimension, then $C$ is bent set since $C_v = \varnothing$ if $v \notin L$, and $C_v = C$ if $v \in L$. Moreover, the space $\F_2^n$ is a bent set for even $n$.
    
    \noindent 3. A partially defined bent function $(h, \F_2^n)$ corresponds to a bent function $h$ on $\F_2^n$. Any Boolean function is a partially defined bent function over a one-element set.

    \noindent 4. In~\cite{ZhangPKP_2026}, it is shown that every subset of $\F_2^n$ of cardinality $2^l$ is a bent set.
\end{remark}

With these definitions, we give the main characterization of bent functions in the $\mathcal{GMM}_{\frac{n}{2}+k}$ class, where $k=0, \ldots, n/2 -1$.

\begin{theorem}\cite[Theorem 6.40]{LogachevSY12_BOOK}\label{th GMM bent Logachev}
	Let $n=2m$. A Boolean function $f_{\phi,h}\in\mathcal{B}_n$ defined by Eq.~\eqref{eq:Mextend} is bent if and only if for all $\alpha \in \F_2^r$ the following conditions are satisfied:
	\begin{enumerate}
		\item $r \leq m$;
		\item $|\phi^{-1}(\alpha)| = 2^{n-2r}$;
		\item $\phi^{-1}(\alpha)$ is a bent set;
		\item $(h, \phi^{-1}(\alpha))$ is a partially defined bent function.
	\end{enumerate}
\end{theorem}

The original proof of this statement can be found in~\cite{Yashchenko97_PPI}, \cite[pp. 245-247]{LogachevSY2004} or~\cite[pp. 173-175]{LogachevSY12_BOOK}. An alternative proof of this statement can be found in~\cite{ZhangPKP_2026}. We also notice that a sufficient condition of bentness in the generalized Maiorana-McFarland class was provided in~\cite{ZhaoWZPC22_DCC}, which is in fact necessary as well (see also Theorem~\ref{th GMM bent Logachev} and~\cite{CARLET_2004_JCompl}).
\begin{theorem}\label{prop:suffic_g_k}
		Let $n$ be even and  $f_{\phi,h}\in \mathcal{M}_
		{r,s}$ with $r=n/2-k,s=n/2+k$ (resp. $f\in \mathcal{GMM}_
		{\frac{n}{2}+k}$) be defined by Eq.~\eqref{eq:Mextend}. If $\phi  \colon \F_2^{\frac{n}{2}+k} \rightarrow
		\F_2^{\frac{n}{2}-k}$ is a $2^{2k}$-to-1 mapping, then $f_{\phi,h}$ is a bent function if and only if for all $\alpha \in \mathbb{F}_2^{\frac{n}{2}-k}$, $\beta\in \mathbb{F}_2^{\frac{n}{2}+k}$,
		\begin{equation*}\label{eq:bentcond}
			\Bigg| \sum_{y\in \mathbb{F}_2^{\frac{n}{2}+k}: \ {\phi}(y)=\alpha} (-1)^{h(y) + \beta \cdot y}\Bigg|=2^k.
		\end{equation*}
	\end{theorem}

For a given bent function $f\in\mathcal{M}^\#$ on $\F_2^n$ it holds that its dual $f^*\in\mathcal{M}^\#$ as well, so $\operatorname{ind}(f)=\operatorname{ind}(f^*)=n/2$; for a detailed discussion, we refer to~\cite{PolujanKP2024WCC,PolujanKP2025_IEEE}. However, in general it is not necessarily the case that $\ind(f)=\ind(f^*)$; for counterexamples, see~\cite[Section~6.2]{KudinPPZ25_JOC}. In view of these observations, we believe that it is interesting to investigate the relation between the linearity index of a bent function $f\in\mathcal{M}^\#_{r,s}$ and the linearity index of its dual $f^*$.

\begin{openproblem}
    Construct bent functions $f\in\mathcal{B}_n$ satisfying $\ind(f)=\ind(f^*)<n/2$.
\end{openproblem}
    
In the following section, we address in detail a special case of the generalized Maiorana-McFarland class corresponding to $k=1$.

\subsection{Almost Maiorana-McFarland bent functions}\label{sec Almost MM}

Maiorana-McFarland bent functions are exactly the bent functions that can be represented as concatenations of affine functions in $n/2$ variables. Now, we consider the next closest case, i.e., bent functions that can be written as concatenations of affine functions in $n/2 - 1$ variables. Such bent functions are called \emph{almost Maiorana-McFarland} in~\cite{KudinPPZZ25_IEEE}, and these are exactly the bent functions $f \in \mathcal{GMM}_{n/2+1}$ of the form $f(x,y) = x \cdot \phi(y) + h(y)$, for $x\in\F_2^{m-1},y\in\F_2^{m+1}$. First, we recall the characterization of almost Maiorana-McFarland bent functions, noting that in this particular framework, this gives a refinement of Theorem~\ref{th GMM bent Logachev}.

\begin{theorem}\cite[Theorem 3.1, Corollary 3.3]{KudinPPZZ25_IEEE} \label{theo: GMM+1condition}
	Let $n=2m$ and let $f\in\mathcal{B}_n$ be a function defined by 
	\begin{eqnarray}\label{eq:GMM+1}
		f(x,y)=x\cdot\phi(y) + 
		h(y), \;\; x \in \F_2^{m-1}, y \in \F_2^{m+1},
	\end{eqnarray}
	where $\phi
	\colon \F_2^{m+1} \rightarrow \F_2^{m-1}$  and  $h \colon \F_2^{m+1} \rightarrow \F_2$. Then, $f$ is a bent function if and only if
	\begin{itemize}
		\item the collection $\left\lbrace \phi^{-1}(a) \colon a \in \F_2^{m-1} \right\rbrace$ is a partition of $\F_2^{m+1}$ into $2$-dimensional affine subspaces (where by $\phi^{-1}(a)$ we denote the set $\{ y \in \F_2^{m+1}  \colon \phi(y)=a \}$), and
		\item for every $a\in\F_2^{m-1}$, the restriction of $h$ on the set $\phi^{-1}(a)$ has odd weight.
	\end{itemize}
	Moreover, if $f$ is a bent function defined by Eq. \eqref{eq:GMM+1}, then the Hamming weight of $h$ satisfies $2^{m-1} \leq \operatorname{wt}(h) \leq 3 \cdot 2^{m-1}$.
\end{theorem}

It is well-known that for a Maiorana-McFarland bent function $f(x,y)=x\cdot\pi(y)+h(y)$, where $x,y\in\F_2^m$, its dual function is given by $f^*(x, y)=y \cdot\pi^{-1}(x) + g(\pi^{-1}(x))$, where $\pi^{-1}$ is the
inverse permutation of $\pi$; for the proof, see~\cite[Proposition 77]{Carlet21_BOOK}. The following result gives an explicit form of the dual for an almost Maiorana-McFarland bent function.

\begin{theorem}\cite[Theorem 3.4]{KudinPPZZ25_IEEE}\label{theo: dualform}
	Let $f$ be a bent function on $\F_2^{2m}$, thus satisfying the conditions in Theorem \ref{theo: GMM+1condition}, defined  by
	\begin{eqnarray*}
		f(x,y)=x\cdot\phi(y) + 
		h(y), \;\; x \in \F_2^{m-1}, y \in \F_2^{m+1},
	\end{eqnarray*}
	where $\phi
	\colon\F_2^{m+1} \rightarrow \F_2^{m-1}$  and  $h \colon \F_2^{m+1} \rightarrow \F_2$. Let $\prec$ denote any total ordering on $\F_2^{m+1}$ (in particular, we can take the lexicographic order on $\F_2^{m+1}$).  Then, the dual function $f^*$ of $f$ is given by
	\begin{equation*} \label{eq:dualform+1}
		f^*(a,b)= \sum_{\substack{y,y' \in \phi^{-1}(a) \\ y \precneqq y'}}(h(y)+y \cdot b)(h(y')+y' \cdot b),
	\end{equation*}
	for all  $a\in \F_2^{m-1}, b \in \F_2^{m+1}$.
\end{theorem}

Clearly, for a fixed permutation $\pi$ of $\mathbb{F}_2^m$, there are $2^{2^m}$ Boolean functions $h \in \mathcal{B}_m$ such that the function $f(x, y) = x \cdot \pi(y) + h(y)$, where $x, y \in \mathbb{F}_2^m$, is bent. Now, we recall a similar counting result for almost Maiorana-McFarland bent functions. Specifically, the following statement specifies, for a fixed 4-to-1 mapping $\phi$ that partitions $\mathbb{F}_2^{n/2+1}$ into 2-dimensional affine subspaces, the exact number of different functions $h$ such that $f(x, y) = x \cdot \phi(y) + h(y)$, where $x \in \mathbb{F}_2^{m-1}$, $y \in \mathbb{F}_2^{m+1}$, is bent.

\begin{theorem}\cite[Corollary 3.6]{KudinPPZZ25_IEEE}\label{th: thenumberofh}
	Let $\phi \colon \F_2^{m+1} \to \F_2^{m-1}$ be a 4-to-1 mapping such that $\left\lbrace \phi^{-1}(a) \colon a \in \F_2^{m-1} \right\rbrace$ is a partition of $\F_2^{m+1}$ into $2$-dimensional affine subspaces. Then, there are exactly $2^{3 \cdot 2^{m-1}}$ functions $h \colon \F_2^{m+1} \to \F_2$ such that the function $f$ defined by
	\begin{eqnarray*} \label{eq: counting}
		f(x,y)=x\cdot\phi(y) +
		h(y), \;\; x \in \F_2^{m-1}, y \in \F_2^{m+1}
	\end{eqnarray*}
	is bent.
\end{theorem}

In~\cite{KudinPPZZ25_IEEE}, it was shown that the decomposition of $\mathbb{F}_2^{m+1}$ into 2-dimensional affine subspaces must be chosen carefully, as an arbitrary selection may lead to bent functions that belong to the $\mathcal{M}^\#$ class. Below, we present an example of such a case.

\begin{theorem}\cite[Proposition IV.1]{KudinPPZZ25_IEEE}\label{th: trivialpartition}
	Let $n=2m$ and let $f \in \mathcal{B}_n$ be a bent function in the $\mathcal{GMM}_{m+1}$ class, thus satisfying the conditions in Theorem \ref{theo: GMM+1condition},  defined by 
    \begin{eqnarray*}\label{eq:GMM+1second}
		f(x,y)=x\cdot\phi(y) +
		h(y), \;\; x \in \F_2^{m-1}, y \in \F_2^{m+1},
	\end{eqnarray*}
	where $\phi
	\colon \F_2^{m+1} \rightarrow \F_2^{m-1}$  and  $h \colon \F_2^{m+1} \rightarrow \F_2$.
	Assume that there exists a  nonzero element $v \in \F_2^{m+1}$, such that, for all $z \in \F_2^{m-1}$, we have that $v \in w_z + \phi^{-1}(z)$, for some $w_z \in \phi^{-1}(z)$. Then, the function $f$ is in the (standard) completed Maiorana-McFarland class.
\end{theorem}

\begin{remark}
    We note that, in general, the properties of the mapping $\phi\colon \F_2^{m+1} \rightarrow \F_2^{m-1}$, and hence of the corresponding partition of $\F_2^{m+1}$ into $2$-dimensional affine subspaces $\lbrace \phi^{-1}(a)\colon a\in\F_2^{m-1} \rbrace$, are decisive for the $\mathcal{M}^\#$-exclusion problem. 
    Indeed, certain choices of the partition inevitably imply the existence of at least one $m$-dimensional $\mathcal{M}$-subspace of the function $f$ defined by Eq.~\eqref{eq:GMM+1}, and hence imply that $f$ belongs to the $\mathcal{M}^\#$ class. 
    Two important classes of such partitions are \emph{trivial} and \emph{non-proper} partitions. 
    
    \noindent 1. Following the terminology of~\cite{KudinPPZZ25_IEEE}, we call a partition of $\F_2^{m+1}$ into 2-dimensional affine subspaces $\{ \phi^{-1}(a) \colon a \in \F_2^{m-1} \}$ \emph{trivial}, if there exists a $2$-dimensional vector subspace $W \subset \F_2^{m+1}$ such that  $\{ \phi^{-1}(a) \colon a \in \F_2^{m-1} \} = \{ c + W \colon c \in \F_2^{m+1} \}$, where $\F_2^{m+1} =\bigcup_{c \in Q}(c+W)$ for the set of coset representative $Q$ with $|Q|=2^{m-1}$. As indicated in~\cite[Corollary IV.4]{KudinPPZZ25_IEEE}, a trivial partition of $ \mathbb{F}_2^{m+1}$ into 2-dimensional affine subspaces induced by preimages of $\phi \colon \F_2^{m+1} \rightarrow \F_2^{m-1}$ defines a bent function $f(x,y)=x\cdot\phi(y) + h(y)$, $x \in \F_2^{m-1}, y \in \F_2^{m+1}$ that has at least $3$ different $\mathcal{M}$-subspaces of dimension $m$ (independently on the choice of a function $h\in\mathcal{B}_{m+1}$).  
    
    \noindent 2. A partition of $\F_2^{m+1}$ into 2-dimensional affine subspaces $\{ \phi^{-1}(a) \colon a \in \F_2^{m-1} \}$ is called in~\cite{KudinPPZZ25_IEEE} \emph{non-proper} if it satisfies the conditions given in Theorem~\ref{th: trivialpartition}; otherwise it is called \emph{proper}. While every non-proper partition induced by the preimages of $\phi \colon \F_2^{m+1} \rightarrow \F_2^{m-1}$ defines a bent function $f\in \mathcal{M}^{\#}$, there also exist some proper partitions that define bent functions in $\mathcal{M}^{\#}$; for examples of such functions, see~\cite[Section V]{KudinPPZZ25_IEEE}. On the other hand, in~\cite{KudinPPZZ25_IEEE}, it was shown that certain proper partitions of $\F_2^5$ into $2$-dimensional affine subspaces yield $2^{78}$ different bent functions on $\F_2^8$ that are outside $\mathcal{M}^\#$.
\end{remark}

\begin{example}
    The mapping $\phi\colon\F_2^5\to\F_2^3$ defined in Example~\ref{ex Almost MM cov EBS} provides an example of a proper partition that yields an almost Maiorana-McFarland bent function outside $\mathcal{M}^\#$.
\end{example}

To conclude this section, we believe that solving the following problems is crucial to understanding almost Maiorana-McFarland bent functions.

\begin{openproblem}
    1. Characterize the partitions of $\F_2^{m+1}$ into 2-dimensional affine subspaces $\{ \phi^{-1}(a) \colon a \in \F_2^{m-1} \}$ that define almost Maiorana-McFarland bent functions $f(x,y)=x\cdot \phi (y)+h(y)$ on $\F_2^{2m}$ outside the $\mathcal{M}^{\#}$ class.
    
    \noindent 2. Determine the mappings $\phi$ such that $f(x,y)=x\cdot \phi (y)+h(y)$ is bent and outside $\mathcal{M}^{\#}$ for all $2^{3 \cdot 2^{m-1}}$ functions $h\in\mathcal{B}_{m+1}$.
\end{openproblem}

\subsection{Bent functions with minimal linearity index}\label{sec Minimal linearity index}

Recall that for a Boolean bent function $f \in \mathcal{B}_n$, its \emph{linearity index} satisfies 
$1 \le \operatorname{ind}(f) \le n/2$. Bent functions $f \in \mathcal{B}_n$ attaining the upper bound are exactly those belonging to the 
completed Maiorana-McFarland class $\mathcal{M}^\#$. Bent functions achieving the lower bound were called \emph{$\ell$-optimal} in~\cite{KudinPPZ25_JOC}.

\begin{definition}\cite{KudinPPZ25_JOC}
	Let $f\in\mathcal{B}_n$ be bent. If $\operatorname{ind}(f)=1$, we say that $f$ is \emph{$\ell$-optimal}, i.e., $f$ has the optimal (i.e., minimal) linearity index.
\end{definition}

One can consider $\ell$-optimal bent functions as the opposite counterpart of Maiorana-McFarland bent functions, since they admit only representations as concatenations of affine functions defined on trivial vector spaces (of dimension one). Notably, examples and constructions of such functions are very limited, the first being provided in Dillon’s thesis.

\begin{example}\cite[p. 102]{Dillon74_PhD}\label{ex: Dillon ell-optimal}
    Consider the $\mathcal{PS}_{ap}$ bent function $f$ on $\F_{2^8}$ given in the trace form as follows
$$f(x) = {\rm Tr}^8_1(x^{15}).$$
One can easily verify that for distinct nonzero $a, b \in \F_{2^8}$, the second-order derivative of $f$ can be expressed as
\begin{equation*}
		D_a D_b f(x) ={\rm Tr}^8_1(c_8 x^8 + c_{12} x^{12} + c_{10} x^{10} + c_9 x^9),
\end{equation*}
where the coefficients $c_t\in\mathbb{F}_{2^8}$ depend on $a$ and $b$, and the coefficient $c_9$ is explicitly given by $c_9 = (ab(a+b))^2$. Clearly, the second-order derivative $D_a D_b f(x)$ vanishes if only if $c_{8} = c_{12} = c_{10} = c_{9} = 0$. However, the condition $c_9 = 0$ implies that either $a$ or $b$ is zero, or that $a = b$, thus contradicting the assumption that $a, b \in \mathbb{F}_{2^8}$ are two distinct nonzero elements.
\end{example}

The study of 4-decompositions of Boolean bent functions naturally led to the investigation of bent functions whose second-order derivatives do not vanish identically, as the following result connecting these two notions demonstrates.

\begin{theorem}\cite{CanteautC03_IEEE}
	\label{Duv=0/1}
	Let $f\in\mathcal{B}_n$ be a bent function with dual $f^*$. Let $u,v$ be linearly 
	independent vectors of $\F_2^n$ and let $S=\langle u,v\rangle^\perp$ be the orthogonal
	complement of $\langle u,v\rangle$. Then,  $f$ has:
	\begin{itemize}
		\item a bent decomposition on $S$ if and only if $D_uD_vf^*= 1$,
		\item a semi-bent decomposition on $S$ if and only if $D_uD_vf^*= 0$, i.e., the second-order derivative vanishes identically,
		\item a $5$-valued decomposition on $S$ if and only if $D_uD_vf^*$ is not a constant function.
	\end{itemize}
\end{theorem}

The first infinite families of Boolean (bent) functions with non-constant second-order derivatives were presented in~\cite{CanteautC03_IEEE} using monomial functions, as indicated by the following result.

\begin{lemma}\label{lemma: Canteauet Charpin Lemma 3}\cite[Lemma 3]{CanteautC03_IEEE}
	Let $ n \geq 8 $, $ n $ even, and consider the Boolean functions in $ \mathcal{B}_n $ of the form
	\begin{equation*}
		f(x)={\rm Tr}^n_1(\alpha x^d), \quad 15 \leq d \leq 2^{n/2} - 1, \alpha \in \F_{2^n}^*,
	\end{equation*}
where the binary weight of $ d $ satisfies $ \operatorname{wt}(d) \geq 4 $. If the 2-ary expansion of $ d $, say $(d_0, \ldots, d_\ell)$, is such that $ d_r = d_s = 1 $ for some $ r $ and $ s $, $ 0 < r < s < \ell $, with $\gcd(s-r, n) = 1$ then the function $ f $ is such that all its second-order derivatives are not constant.
\end{lemma}

The following particular instance of the monomial functions considered above not only provides an explicit family of $\ell$-optimal bent functions, but also generalizes Example~\ref{ex: Dillon ell-optimal} of Dillon in dimension eight.

\begin{theorem}\label{lemma: Canteauet Charpin Theorem 10}\cite[Theorem 10]{CanteautC03_IEEE}
	When a function $f$, defined as in Lemma~\ref{lemma: Canteauet Charpin Lemma 3} is bent, then its dual $f^*$ admits neither a bent 4-decomposition nor a semi-bent 4-decomposition. In particular, this property holds for the infinite class of bent functions of $ \mathcal{B}_n $, $ n = 2m $
	\begin{equation*}
		g_n (x)= {\rm Tr}^n_1(\alpha x^{2^m-1}),
	\end{equation*}
where $\alpha \in \F_{2^m}$ satisfies $\sum_{i=0}^{2^m} {\rm Tr}^n_1(\alpha \gamma^i) = 2^{m-1}$ (where $\gamma$ is a $(2^m + 1)$-th root of unity). Moreover, each $ g_n $ is equal to its dual. So $ g_n $ (as well as $g^*_n$) admits neither a bent 4-decomposition nor a semi-bent 4-decomposition.
\end{theorem}

As soon as an explicit univariate family of bent functions is specified, it is natural to identify such functions within the well-known ``large'' classes of bent functions. The following class of permutations, whose second-order derivatives do not vanish, was introduced in~\cite{PasalicPKZ24_IEEE} and later used in~\cite{KudinPPZ25_JOC} to identify $\ell$-optimal bent functions within Carlet’s $\mathcal{D}_0$ class.

\begin{definition}~\cite{PasalicPKZ24_IEEE}
    A permutation $\pi$ of $\F_2^m$ is said to have~\eqref{eq: P1} property, if
	\begin{equation}\label{eq: P1} \tag{$P_1$}
		D_vD_w\pi\neq0_m \quad \mbox{ holds for all linearly independent } v,w\in\F_2^m.
	\end{equation}
\end{definition}

The following result provides a characterization of bent functions with the minimal relaxed linearity index from the class $\mathcal{D}_0$, that arise from permutations satisfying the property~\eqref{eq: P1}. Consequently, it gives a method for constructing bent functions satisfying $\operatorname{r\text{-}ind}(f) < n/2$, thus providing a solution to~\cite[Open Problem~1]{ZhangPBW24_INFCOMP} concerning construction methods for such functions.

\begin{theorem}\cite{KudinPPZ25_JOC}\label{theo: MofDzero} Let $\pi$ be a permutation of $\F_2^m$, $m \geq 4$ which has the property \eqref{eq: P1}. Define $f\colon\F_2^{m}\times\F_2^{m} \to \F_2$ by $f(x,y)=x \cdot \pi(y) + \delta_{0}(x)$, for all $x,y \in \F_2^{m}$. Then, $\operatorname{r-ind}(f)\leq 2$. Furthermore, $\operatorname{r-ind}(f)=1$, implying that $\operatorname{ind}(f)=1$, if and only if $\pi$ has no components with linear structures.
\end{theorem}

 The following result indicates that non-quadratic monomial permutations $\pi$ on $\F_2^m$ with the~\eqref{eq: P1} property can generate $\ell$-optimal bent functions in the $\mathcal{D}_0$ class.

\begin{corollary}\cite{KudinPPZ25_JOC}\label{cor:loptimal} Let $\pi$ be a permutation of $\F_2^m$, $m \geq 4$ which has the property \eqref{eq: P1} and additionally satisfies the condition that $ a_1 \cdot  D_{b_2}\pi(y) \neq 0$, for any nonzero $a_1,b_2 \in \F_2^m$.  Then, the bent function $f(x,y)=x \cdot \pi(y) + \delta_{0}(x)$, where $x,y \in \F_2^{m}$, is outside 
	$\cM^\#$ and also $\ell$-optimal since its $\operatorname{r-ind}(f)=1$. In particular, the same is true if  $\deg(\pi)>2$ and $\pi$ is a monomial permutation satisfying \eqref{eq: P1}.
\end{corollary}

Despite significant progress towards the understanding of the generalized completed Maiorana-McFarland class (e.g., \cite{ZhaoWZPC22_DCC,KudinPPZ25_IEEE,KudinPPZZ25_IEEE}), the following problems remain largely open.

\begin{openproblem}
1. Provide explicit constructions of bent functions $f_{\phi,h}\in\mathcal{M}_{r,s}$ on $\F_2^n$ (with $r<n/2$) having linearity index $r$. \ \\ 
\noindent 2. Provide necessary and sufficient conditions for bent functions $f_{\phi,h}\in\mathcal{B}_n$ in $\mathcal{M}^\#_{r-1,s+1}$ to be outside $\mathcal{M}^\#_{r,s}$.
\end{openproblem}

\section{\texorpdfstring{$\mathcal{M}$-Subspaces}{M-Subspaces} of bent functions} \label{sec M-subspaces}
In this section, we consider in detail the $\mathcal{M}$-subspaces of bent functions, their algebraic and combinatorial properties, as well as construction methods for bent functions with a unique or with many $\mathcal{M}$-subspaces.
\subsection{Algebraic and combinatorial properties}\label{sec M-subspaces basics}

In the previous sections, we mentioned that for a Boolean function $f \in \mathcal{B}_n$, its linearity index, defined as the maximum dimension of an $\mathcal{M}$-subspace, satisfies $1 \le \ind(f) \le n$. The next result gives a connection between the size of a subspace $V \subseteq \mathbb{F}_2^n$ on whose cosets $f$ is affine and the maximal absolute value of its Walsh transform (the variations of this statement can be found in the literature, but to the best of our knowledge the considered formulation is more general). Based on this characterization, one immediately deduces that the linearity index of a bent function in $n$ variables is at most $n/2$.

\begin{proposition}\cite{Dilawar}
Let $f\in\mathcal{B}_n$ be a Boolean function, and let $$M = \max_{u \in \mathbb{F}_2^n} |W_f(u)|.$$ Let $V$ be a subspace of $\mathbb{F}_2^n$. If $|V| > M$, then $f$ cannot be affine on any coset of $V$. Consequently, for a bent function $f\in\mathcal{B}_n$, it holds that $\ind(f)\le n/2$.
\end{proposition}

\begin{remark}
    Another way to upper bound the size and cardinality of $\mathcal{M}$-subspaces is through the relaxed linearity index. For these EA-invariants, it is straightforward to verify that every $f \in \mathcal{B}_n$ satisfies $\ind(f) \le \rind(f)$. Moreover, the total number of $r$-dimensional $\mathcal{M}$-subspaces of a function $f \in \mathcal{B}_n$, denoted $\left| \mathcal{MS}_r(f) \right|$ (as well as the total number of $r$-dimensional relaxed $\mathcal{M}$-subspaces $\left| \mathcal{RMS}_r(f) \right|$), is invariant under EA-equivalence; see~\cite{PolujanP20_DCC} for details.
\end{remark}

For bent functions $f \in \mathcal{B}_n$, an upper bound on the number of $\mathcal{M}$-subspaces of maximal dimension $n/2$ is known; moreover, the bent functions achieving this bound have been characterized. For details, we refer to~\cite[Theorem 2]{Kolomeec17_DCC} and~\cite[Proposition 1.38]{Polujan21_PhD}.

\begin{theorem}\label{th: max nr of msubspaces}
    Let $f\in\mathcal{B}_n$ be a bent function. The number of $n/2$-dimensional $\mathcal{M}$-subspaces of $f$ is at most
    $$|\mathcal{MS}_{n/2}(f)| \leq \prod_{i=1}^{n/2} \left(2^i + 1\right),$$
    with equality if and only if the function $f$ is quadratic.
\end{theorem}

\begin{remark}
    Let $f$ be a quadratic Boolean function on $\mathbb{F}_2^{n}$, i.e., it can be written in the form $f(x)=x Qx^\top+c x^\top+\varepsilon$, where $Q$ is an upper-triangular $n\times n$ matrix over $\mathbb{F}_2$ with zero diagonal, $c\in\F_2^n,\epsilon\in\F_2$. Define the bilinear form
$$
B_f(a,b)=f(a+b)+f(a)+f(b)+f(0_n)=a Qb^\top.
$$
The function $f$ is \emph{bent} if and only if $B_f$ is nondegenerate, i.e., for each nonzero $a$ there exists $b$ with $B_f(a,b)=1$. An $\mathcal{M}$-subspace $V$ of maximal dimension $n/2$ is precisely a subspace satisfying $B_f(a,b)=0$ for all $a,b\in V$,
which is called a \emph{totally isotropic subspace} with respect to the alternating form $B_f$; see also~\cite[Section 2]{Cameron1999}. From this geometric perspective, non-quadratic Maiorana-McFarland bent functions can be considered as generalizations of non-degenerate quadratic forms.
\end{remark}

While Theorem~\ref{th: max nr of msubspaces} establishes the ``global'' maximum number of $\mathcal{M}$-subspaces within the class of all bent functions, no such bound is known for bent functions of a fixed algebraic degree. This motivates the following open problem.

\begin{openproblem}
    Determine the maximum number of $\mathcal{M}$-subspaces for bent functions $f\in\mathcal{B}_n$ of algebraic degree $3\le d\le n/2$.
\end{openproblem}

Note that quadratic Maiorana-McFarland bent functions are not the only ones that may possess multiple $\mathcal{M}$-subspaces of the maximum dimension. In the following result, we summarize some known construction methods of bent functions with this property.

\begin{proposition}\cite{PasalicPKZ24_IEEE}
    Let $\pi$ be a permutation of $\F_2^m$, $h\in\mathcal{B}_m$ be a Boolean function, and $f_{\pi,h} (x,y)=x \cdot \pi(y)+ h(y)$ be a Maiorana-McFarland bent function on $\F_2^m\times\F_2^m$.
    \begin{itemize}
        \item[1.] If $\pi$ has a nonzero linear structure, then $f_{\pi,h}$ has at least two $m$-dimensional $\mathcal{M}$-subspaces.
        \item[2.] Let $h=0$ and assume there exists an $(m-1)$-dimensional subspace $S \subset \F_2^m$ for which $D_aD_b \pi=0_m$, for all $a,b \in S$. Then, $f_{\pi,0}$ has at least two $m$-dimensional $\mathcal{M}$-subspaces.
    \end{itemize}
\end{proposition}

Similarly, almost Maiorana-McFarland bent functions of the form $f(x,y) = x \cdot \phi(y) + h(y)$, defined for $x\in\F_2^{m-1}$ and $y\in\F_2^{m+1}$, may possess many $\mathcal{M}$-subspaces of maximal dimension $m$, provided that the mapping $\phi$ and the Boolean function $h$ are appropriately chosen. For a detailed discussion of such sufficient conditions, we refer to~\cite[Section~IV.A]{KudinPPZZ25_IEEE}.

\begin{remark}
We also note that $\mathcal{M}$-subspaces play an important role in the analysis of vectorial bent functions. In~\cite{KPP_2026_WCC}, the authors proved that any two quadratic Boolean functions on $\F_2^n$ (including bent functions) share at least $2^{n/2}+1$ common $\mathcal{M}$-subspaces of dimension $n/2$. This result implies that every quadratic vectorial bent function $F\colon\F_2^n\to\F_2^2$ has, up to extended-affine equivalence, the canonical form
$$(x,y)\in\F_2^{n/2}\times\F_2^{n/2}
\mapsto
(x\cdot y+g(y),\,x\cdot L(y)+h(y)),$$
where $L$ is a linear complete permutation of $\F_2^{n/2}$ and $g,h\in\mathcal{B}_{n/2}$ satisfy $\deg(g)$, $\deg(h)\leq 2$.
Here, a permutation $L$ is called \emph{complete} if both mappings $y\mapsto L(y)$ and $y\mapsto L(y)+y$ are permutations of $\F_2^{n/2}$.
\end{remark}

\subsection{Bent functions with a unique \texorpdfstring{$\mathcal{M}$-subspace}{M-subspace}, permutations with the \texorpdfstring{$(P_1)$}{(P1)} and \texorpdfstring{$(P_2)$}{(P2)} properties, and their applications} \label{se M-subspace unique}

Dillon, in his thesis~\cite{Dillon74_PhD}, observed that every Maiorana-McFarland bent function
$f(x,y)=x\cdot \pi(y)+h(y)$ on $\F_2^m\times \F_2^m$ has an $m$-dimensional $\mathcal{M}$-subspace
$U=\F_2^m\times\{0_m\}$, which is called in~\cite{PasalicPKZ24_IEEE} the \emph{canonical $\mathcal{M}$-subspace} of dimension~$m$. As indicated in Section~\ref{sec M-subspaces basics}, a (Maiorana-McFarland) bent function may possess many ($m$-dimensional) $\mathcal{M}$-subspaces. In this section, we summarize construction methods for bent functions in $\mathcal{B}_{2m}$ having a unique $m$-dimensional $\mathcal{M}$-subspace and discuss applications of such functions, namely flexible constructions of bent functions outside $\mathcal{M}^\#$, as well as counting the number of bent functions in $\mathcal{M}^\#$.

The first construction method of Maiorana-McFarland bent functions with the unique canonical $\mathcal{M}$-subspace of the maximal dimension is based on the use of permutations satisfying the property~\eqref{eq: P1}.

\begin{theorem}\cite[Theorem III.1]{PasalicPKZ24_IEEE}\label{theo: unique P1} Let $\pi$ be a permutation of $\F_2^m$ with the~\eqref{eq: P1} property, i.e., $D_vD_w\pi\neq0_m$,  for all linearly independent $v,w\in\F_2^m$. Define $f\colon\F_2^{m}\times\F_2^{m} \to \F_2$ by $$f(x,y)=x \cdot \pi(y) + h(y),$$ for all $x,y \in \F_2^{m}$, where $h\in\mathcal{B}_m$ is an arbitrary Boolean function. Then, the following hold:
	\begin{itemize}
		\item[1)] The permutation $\pi$ has no linear structures.
		\item[2)] The vector space $V=\F_2^m \times \{0_m \}$ is the only $m$-dimensional $\mathcal{M}$-subspace of $f$.
	\end{itemize}
\end{theorem}

Theorem~\ref{theo: unique P1} provides a sufficient condition for $\pi$ under which the function $f(x,y)=x \cdot \pi(y) + h(y)$ admits only the canonical $m$-dimensional $\mathcal{M}$-subspace, independently of the choice of $h$. If the function $h$ is specified, or if additional assumptions are imposed on $h$, then it is possible to derive less restrictive conditions for the permutation $\pi$, so that $f$ admits only one $m$-dimensional $\mathcal{M}$-subspace.
For example, the following result gives a necessary and sufficient condition for a permutation
$\pi$ of $\F_2^m$, with $\deg(\pi) < m-1$, under which the function
$
f(x,y) = x \cdot \pi(y) + h(y),$ where 
$h(y) = \delta_0(y) = \prod_{i=1}^m (y_i + 1),$
admits only the canonical $\mathcal{M}$-subspace
$\F_2^m \times \{0_m\}$.
  
\begin{proposition}\cite[Proposition 12]{PasalicPKZ24_IEEE}
	Let $\pi$ be a permutation of $\F_2^m$ with $\deg(\pi)<m-1$, and let $f\colon \F_2^{m}\times\F_2^{m} \to \F_2$ be the function defined by 
	$$f(x,y)= x \cdot \pi(y) + \delta_0(y), \text{ for all } x,y \in \F_2^{m}.$$
	Then, $f$ has only one $m$-dimensional $\mathcal{M}$-subspace if and only if $\pi$ has no linear structures.
\end{proposition}

\begin{remark}
Determining possible $\mathcal{M}$-subspaces for $\mathcal{PS}_{ap}$ bent functions in $n$ variables seems to be a challenging problem. As the constructions in~\cite{CanteautC03_IEEE,AnbarKKMPP2026,AlkanAAEKT2026} indicate, the number and dimensions of $\mathcal{M}$-subspaces of bent functions in $\mathcal{PS}_{ap}$ class can vary significantly. For example, they can have a unique $\mathcal{M}$-subspace of dimension $n/2$, or have only the trivial ones of dimension at most $1$.
\end{remark}

The property~\eqref{eq: P1} ensures that Maiorana-McFarland bent functions of the form $f(x,y)= x \cdot \pi(y)$ have a unique canonical $\mathcal{M}$-subspace, due to the non-vanishing property of their second-order derivatives. However, numerous permutations $\pi$ that yield Maiorana-McFarland bent functions $f(x,y)= x \cdot \pi(y)$ with only the canonical $\mathcal{M}$-subspace of the maximal dimension do not satisfy this property. The following concept aims to capture a broader class of such permutations.

\begin{definition}\cite{PasalicPKZ24_IEEE}
	Let $\pi$ be a permutation of $\F_2^m$. Let $S$ be a subspace of $\F_2^m$ of dimension $m-k$, with $1 \leq k  \leq m-1$, such that $D_aD_b\pi=0_m$ for all $a,b \in S$. Then, $\pi$ satisfies the property \eqref{eq: P2} with respect to $S$ if there does not exist a vector subspace $V$ of $\F_2^m$  with $\dim(V) = k$ such that
	\begin{equation}\label{eq: P2} \tag{$P_2$}
			v \cdot D_a\pi(y) = 0,  \textnormal{ for all } a \in S, \; \textnormal{ all } y \in \F_2^m,  \textnormal{ and  for all } v \in V.
	\end{equation}
	If $\pi$ satisfies this property with respect to all vector subspaces $S$ of $\F_2^m$ of arbitrary dimension $1 \leq \dim(S)  \leq m-1$, then we simply say that $\pi$  satisfies \eqref{eq: P2}.  
\end{definition}

The following statement provides an alternative characterization of the~\eqref{eq: P2} property, which simplifies a specification of such permutations.

\begin{proposition}\label{prop:Vs(f)}\cite{KudinPPZ25_JOC}
Let $\pi$ be a permutation of $\F_2^m$ and let $S$ be a $k$-dimensional subspace of $\F_2^m$, $k\in \lbrace 1, 2, \ldots , m-1 \rbrace$, such that $D_aD_b \pi = 0_m$, for all $a,b \in S$. Denote by $V_S(f)$ the subspace of $\F_2^m$ generated by the set 
$ \lbrace D_a \pi(y) : a \in S \text{ and } y \in \F_2^m \rbrace.$ 
 Then, the permutation $\pi$ satisfies the property~\eqref{eq: P2} with respect to the subspace $S$ if and only if $\dim(V_S(f)) > \dim(S)=k$. 
\end{proposition} 

\begin{remark}
    As shown in~\cite[Theorem 3.8]{KudinPPZ25_JOC}, \label{th:refineddegreeP2} for a permutation $\pi$ over $\F_2^m$, to satisfy the property  ~\eqref{eq: P2}, it is enough to verify that it satisfies \eqref{eq: P2} for all subspaces $S$ such that $\dim(S) \leq m-\deg(\pi)+1$, where $D_aD_b \pi =0_m$, for all $a,b \in S$. In particular, when $\deg(\pi)=m-1$, to verify the property~\eqref{eq: P2} it is enough to check that it satisfies \eqref{eq: P2} for all 2-dimensional subspaces and that it has no linear structures.
\end{remark}

With this property, it becomes possible to completely characterize non-quad\-ra\-tic ``scalar product-like'' Maiorana-McFarland bent functions that possess a unique canonical $\mathcal{M}$-subspace of maximal dimension.

\begin{proposition}\cite[Proposition 14]{PasalicPKZ24_IEEE}\label{prop:suffcondunique}	
	Let $\pi$  be a non-affine permutation of $\F_2^m$ and $f(x,y)= x \cdot \pi(y)$ be a bent function on $\F_2^m \times \F_2^m$ in $\cM$.
	Then, the permutation $\pi$ has the property \eqref{eq: P2} if and only if the only $m$-dimensional $\mathcal{M}$-subspace of $f$ is $\F_2^m \times \{0_m\}$.
\end{proposition}

We now discuss the relationship between the properties~\eqref{eq: P1} and~\eqref{eq: P2} and the known construction methods for permutations satisfying them. In~\cite[Proposition 3.9]{KudinPPZ25_JOC}, it was shown that if a permutation $\pi$ of $\F_2^m$ has the property \eqref{eq: P1}, then it also has the property \eqref{eq: P2}. However, as indicated in~\cite[Remark 16]{PasalicPKZ24_IEEE}, there exist many examples of (quadratic) permutations (on $\F_2^5$), which possess the~\eqref{eq: P2} property, but not~\eqref{eq: P1}. 

One straightforward class of permutations satisfying the property~\eqref{eq: P1} was identified in~\cite[Corollary~18]{PasalicPKZ24_IEEE}, where it was shown that a quadratic permutation $\pi$ on $\F_2^m$ satisfies~\eqref{eq: P1} if and only if $\pi$ is a quadratic APN permutation of $\F_2^m$; for examples and constructions of APN permutations, we refer to~\cite{BeierleCLP22,LiK24}. Besides that, permutations with the property~\eqref{eq: P1} can be obtained either by concatenating suitably chosen permutations~\cite[Proposition~22]{PasalicPKZ24_IEEE} or by deriving new ones from the known ones satisfying the~\eqref{eq: P1} property~\cite[Corollary~23]{PasalicPKZ24_IEEE}. An explicit construction of permutations satisfying the property~\eqref{eq: P2} was given in~\cite[Proposition~21]{PasalicPKZ24_IEEE} using monomials. For example, it was observed that the permutation $\pi(y)=y^{2^t+1}$ of $\F_{2^m}$, where $s=\gcd(t,m)=2$, $m=2r$, and $r\ge 3$ is odd, satisfies~\eqref{eq: P2}. Secondary construction techniques for permutations satisfying the property~\eqref{eq: P2} are analogous to those developed for the property~\eqref{eq: P1}; see~\cite{PasalicPKZ24_IEEE,KudinPPZ25_JOC}. Finally, additional examples of permutations satisfying the properties~\eqref{eq: P1} and~\eqref{eq: P2} in small dimensions can be found among the affine-inequivalent permutations of $\F_2^4$ listed in~\cite{DeCanniere2007}, as well as in the classification of affine-inequivalent quadratic permutations of $\F_2^5$ given in~\cite{BozilovBS17}.

\begin{remark}
    In \cite{AlkanAAEKT2026}, the authors introduced the property~$(P_3)$ to characterize the uniqueness of the $m$-dimensional $\mathcal{M}$-subspace of a Maiorana-McFarland bent function $f(x,y)=\operatorname{Tr}_1^m(x\pi(y))$ on $\F_{2^m}\times\F_{2^m}$, which, moreover, contains all non-trivial $\mathcal{M}$-subspaces of $f$. For the (technical) definition of this property and its connection to properties~\eqref{eq: P1} and~\eqref{eq: P2}, we refer the reader to~\cite{AlkanAAEKT2026}.
\end{remark}

\begin{remark}
    There exist two important applications of bent functions with a unique $\mathcal{M}$-subspace.
    
    \noindent\emph{1.} In~\cite{KolomeecB25}, using an estimate on the number of Maiorana-McFarland bent functions with a unique $\mathcal{M}$-subspace, the authors obtained asymptotically tight bounds for the cardinality of the $\mathcal{M}^\#$ class in $n$ variables for every even $n$.
    
    \noindent\emph{2.} In~\cite{PasalicPKZ24_IEEE}, the authors presented generic methods for constructing bent functions $f \in \mathcal{B}_{n+2}$ outside $\mathcal{M}^\#$, for even $n \ge 6$, using the 4-concatenation $f = f_1 || f_2 || f_3 || f_4$ of functions $f_i \in \mathcal{B}_n$. One direct approach relies on concatenating four suitably chosen Maiorana-McFarland bent functions $f_1, f_2, f_3, f_4 \in \mathcal{B}_n$ that share a unique $\mathcal{M}$-subspace of maximal dimension $n/2$ (which can be constructed using permutations with the~\eqref{eq: P1} and~\eqref{eq: P2} properties). Since all bent functions in dimension $n = 6$ belong to $\mathcal{M}^\#$, this provides a better explanation of the origin of bent functions outside $\mathcal{M}^\#$ in dimension $n = 8$. 
\end{remark}

Instead of concatenations $f = f_1||f_2||f_3||f_4$ of suitably chosen bent functions $f_i \in \mathcal{B}_n$ in $\mathcal{M}^\#$ with a common maximal $\mathcal{M}$-subspace, it is natural to consider concatenations of functions $f_i \in \mathcal{B}_n$ (not necessarily Maiorana-McFarland or even bent) which either have no common $\mathcal{M}$-subspaces, or only $\mathcal{M}$-subspaces of sufficiently small dimensions. Such construction approaches were considered in detail in~\cite{PasalicPKZ24_IEEE,PolujanPKZ24_CCDS}.  To get a better understanding of such building blocks with small $\mathcal{M}$-subspaces, a generalization of the~\eqref{eq: P1} property was introduced in~\cite{KudinPPZZ25_IEEE}, which was proven to be useful in characterizing Boolean functions with the unique canonical $\mathcal{M}$-subspace in the generalized Maiorana-McFarland class.

\begin{theorem}\label{th:uniqueMsubspace}\cite[Theorem V.3]{KudinPPZZ25_IEEE}
	Let $\phi\colon\F_2^{m+k} \rightarrow \F_2^{m-k} $, where $ 0 < k < m$, be a mapping which satisfies the extended~\eqref{eq: P1ext} property, i.e.,
    \begin{equation}\label{eq: P1ext} \tag{$P_1^*$}
	D_vD_w \phi  \neq 0_{m-k} \mbox{ for all linearly independent } v,w\in\F_2^{m+k}.
	\end{equation}
    Define the function $f(x,y)=x \cdot \phi(y) + h(y)$  on $\F_2^{2m}$, where $x \in \F_2^{m-k} $, $y \in \F_2^{m+k}$.  Then, for any $\cM$-subspace $V$ of $f$ it holds that $\dim(V) \leq m-k$. Furthermore, $f$ has the unique $\cM$-subspace $V=\F_2^{m-k} \times \{0_{m+k} \}$ of the maximal dimension $m-k$ if and only if there does not exist a nonzero $b_2 \in \F_2^{m+k}$ and an $(m-k-1)$-dimensional subspace $U$ of $\F_2^{m-k}$ such that $u\cdot D_{b_2}\phi(y) = 0$, for all $u \in U$ and $y \in \F_2^{m-k}$. 
    In particular, when $f$ is bent then $f$ is outside $\cM^\#$. 
\end{theorem}

Note that, in general, bent functions in the generalized Maiorana-McFarland class do not necessarily arise from mappings with the extended~\eqref{eq: P1ext} property~\cite[Remark V.4]{KudinPPZZ25_IEEE}. We believe that providing additional constructions of such bent functions remains an interesting open problem.

\begin{openproblem}
    Provide further construction methods of bent functions on $\F_2^n$ in~$\mathcal{GMM}_{n/2+k}$ with the unique canonical $\mathcal{M}$-subspace.
\end{openproblem}

The set of all bent functions on $\F_2^n$ is a union of bent functions from the completed generalized Maiorana-McFarland classes $\mathcal{GMM}^\#_{n/2+k}$, for $k=0,\ldots,n/2-1$. Thus, to count all bent functions, which is a well-known open problem~\cite{PotapovTT24}, it suffices to determine the cardinalities of the completed generalized Maiorana-McFarland classes. For $k=0$, such results were obtained in~\cite{KolomeecB25} due to a better understanding of Maiorana-McFarland bent functions with a unique canonical $\mathcal{M}$-subspace. A natural next step is to consider the remaining cases.

\begin{openproblem}
Establish the cardinalities (or their ``good'' estimates) of the completed generalized Maiorana-McFarland classes $\mathcal{GMM}_{n/2+k}$ on $\F_2^n$, for $k>0$. In particular, solve this problem for almost Maiorana-McFarland bent functions $(k=1)$.
\end{openproblem}

\subsection{\texorpdfstring{$\mathcal{M}$-Subspaces}{M-Subspaces} of secondary constructions of bent functions}\label{sec M-subspaces primary secondary}

In the previous sections, we considered whether the classes of bent functions derived from primary constructions belong to $\mathcal{M}^\#$. We now turn to secondary constructions of bent functions, which, thanks to the flexibility in choosing their building blocks, allow not only to determine whether the resulting functions lie inside or outside $\mathcal{M}^\#$, but also, in many cases, to determine their linearity index and even to describe their $\mathcal{M}$-subspaces completely. In general, there are many ways to construct new bent functions from known ones; for a detailed discussion, we refer the reader to~\cite[Section~6.1.16]{Carlet21_BOOK}. In what follows, we focus on the secondary constructions of bent functions that are best understood with respect to the membership of the induced functions in the $\mathcal{M}^\#$ class: direct and indirect sum, as well as 2- and 4-concatenations.

The \emph{direct sum} $h(x,y)=f(x)+g(y)$, where $x\in\F_2^{n}$ and $y\in\F_2^{m}$, is bent iff $f\in\mathcal{B}_n$ and $g\in\mathcal{B}_m$ are both bent, as originally observed by Dillon~\cite{Dillon74_PhD} and Rothaus~\cite{Rothaus76_JCTA}. In his thesis, Dillon used the direct sum construction to prove that for all $n\ge 8$, there exist bent functions on $\F_2^n$ that are outside the completed Maiorana-McFarland class $\mathcal{M}^\#$. His construction was based on the concatenation of arbitrary bent functions with the function from Example~\ref{ex: Dillon ell-optimal} that has no vanishing second-order derivatives. However, since the degree of bent functions constructed with this approach is at least $4$, it was not clear whether cubic bent functions on $\F_2^n$ can be outside $\mathcal{M}^\#$ for infinitely many $n$. To prove that cubic bent functions on $\F_2^n$ outside $\mathcal{M}^\#$ exist for all $n\ge 10$ (opposite to dimensions $n=6$ and $n=8$), further connections between $\mathcal{M}$-subspaces and relaxed $\mathcal{M}$-subspaces were investigated in~\cite{PolujanP20_DCC,Polujan21_PhD}. For example, Corollary~4.6 of~\cite{PolujanP20_DCC} provides sufficient conditions in terms of relaxed linearity indexes of bent functions $f\in\mathcal{B}_n$ and $g\in\mathcal{B}_m$ such that the direct sum is outside $\mathcal{M}^\#$. In particular, if $\rind(f)<n/2$ and $\rind(g)\le m/2$, then the direct sum $f+g\in\mathcal{B}_{n+m}$ is outside $\mathcal{M}^\#$. In~\cite{PasalicBZW23_DCC,ZhangPBW24_INFCOMP}, these conditions were further relaxed, and in~\cite[Theorem 7]{LiKMPZ25_IEEE} it was shown that the direct sum $f+g$ is outside $\mathcal{M}^\#$ if and only if at least one of the functions $f$ or $g$ is outside $\mathcal{M}^\#$.

The \emph{indirect sum} $h(x,y)=f_1(x)+g_1(y)+(f_1+f_2)(x)\,(g_1+g_2)(y)$, where $x\in\F_2^{n}$ and $y\in\F_2^{m}$, is bent if $f_1,f_2\in\B_n$ and $g_1,g_2\in\B_m$ are bent, as observed by Carlet~\cite{Carlet_Indirect}. Similarly to the direct sum, the indirect sum requires no additional structure on the bent functions used in its definition. There exist many variations and generalisations of this construction (see~\cite[Section~6.1.16]{Carlet21_BOOK} and the works~\cite{CarletGenIndirect2012,CarletGL_JCTA_2014,Zhang2016AAECC}). We now summarize the known results concerning the exclusion of bent functions obtained through this construction and its variants from~$\mathcal{M}^\#$.  In~\cite[Theorem~4]{ZhFR2015IT}, in the context of constructing bent-negabent functions outside the $\mathcal{M}^\#$ class, (technical) sufficient conditions on $f_1,f_2,g_1,g_2$ were provided so that the indirect sum does not belong to $\mathcal{M}^\#$. Further design methods of bent functions outside~$\mathcal{M}^\#$ were proposed in~\cite{ZhangPBW24_INFCOMP}, where sufficient conditions on the bent functions $f_i$ and $g_i$ so that $h$ defined by the indirect sum is provably outside~$\mathcal{M}^\#$ were given. Notably, the same work proposes methods that either impose assumptions on the class inclusion of the building blocks in terms of the linearity index and the relaxed linearity index (see Subsection~3.2.1 of~\cite{ZhangPBW24_INFCOMP}), or do not rely on such assumptions. For example, it was shown that if $f_1+f_2$ and $g_1+g_2$ are bent, then the resulting indirect sum $h$ does not belong to~$\mathcal{M}^\#$ (Theorem~3.4 of~\cite{ZhangPBW24_INFCOMP}). The latter approach was further developed in~\cite{PolujanKP2025_IEEE}, where the authors provided sufficient conditions on the building blocks of the variation of the indirect sum introduced in~\cite{CarletGL_JCTA_2014} to give the first constructions of rotation-symmetric bent functions outside $\mathcal{M}^\#$. Further variations of the indirect sum construction were also studied in~\cite[Theorem~3]{Zhang2016AAECC}.

The \emph{2-concatenation} $h=f_1||f_2$, defined by  
$h(x,y)=(1+y)f_1(x)+y f_2(x)$, where $x\in\F_2^n$ and $y\in\F_2$, is bent iff $f_1,f_2\in\mathcal{B}_n$ are two disjoint-spectra semi-bent functions~\cite[Theorem~6]{ZhengZhang_IEEE_2001}, where $n$ is odd. In~\cite{KudinPPZ25_IEEE}, the authors provided a complete description of $\mathcal{M}$-subspaces of the 2-concatenation $h=f_1||f_2$ (where $f_1$ and $f_2$ are arbitrary). As a consequence of this result, the authors obtained a necessary and sufficient condition for inclusion of bent functions $h=f_1||f_2$ in the $\mathcal{M}^\#$ class.

\begin{corollary} \cite[Corollary 1]{KudinPPZ25_IEEE} \label{cor:completecharact} Let $f_1,f_2 \in \mathcal{B}_{n}$, $n=2k+1$, be Boolean functions such that $f=f_1 \vert \vert f_2 \in \mathcal{B}_{n+1}$ is a bent function. Then, the function $f$ is outside the $\cM^{\#}$ class if and only if the following conditions hold:
	\begin{enumerate}
		\item The functions $f_1$ and $f_2$ do not share a common $(k+1)$-dimensional $\cM$-subspace;
		\item For every vector $u \in \F_2^n$ and every $k$-dimensional $\cM$-subspace $V \subset \F_2^n$ of both $f_1$ and $f_2$, there is $a \in V$ such that
		$D_af_1(z)+D_af_2(z+u) \neq 0, \text{ for some } z \in \F_2^n$.
	\end{enumerate}
\end{corollary}

Finally, we consider the \emph{4-concatenation} $h=f_1||f_2||f_3||f_4$, whose algebraic expression is given by 
$$h(x,y,z)=(y+1)(z+1)f_1(x)+(y+1)zf_2(x)+y(z+1)f_3(x)+yzf_4(x),$$
where $x\in\F_2^n$ and $y,z\in\F_2$. This function is bent if and only if certain technical conditions on the Walsh spectra of the involved functions $f_i\in\mathcal{B}_n$ are satisfied~\cite[Theorem~III.1]{HodzicPZ_IEEE_2019}. Note that, in the case that $h$ is bent, the four Boolean functions $f_1,f_2,f_3,f_4$ are simultaneously bent, semi-bent, or $5$-valued spectrum functions~\cite{CanteautC03_IEEE}. If $f_i\in\B_n$ are all bent, then the necessary and sufficient condition for the function $h=f_1||f_2||f_3||f_4\in\mathcal{B}_{n+2}$ to be bent is that the \emph{dual bent condition} is satisfied~\cite{HodzicPZ_IEEE_2019}, i.e., $f_1^*+f_2^*+f_3^*+f_4^*=1$. In general, it is not an easy task to find four bent functions satisfying the dual bent condition (for examples of such constructions, we refer to~\cite{PolujanPKZ24_CCDS}). However, some simple choices of the functions $f_i$, e.g., $h=f_1||f_2||f_1||(f_2+1)$ and $h=f_1||f_1||f_1||(f_1+1)$ are often considered, since in this case the dual bent condition is automatically satisfied.

In~\cite{PasalicBZW23_DCC,PasalicPKZ24_IEEE,KudinPPZ25_JOC}, various sufficient conditions on the functions $f_i$ were provided that guarantee the exclusion of the 4-concatenation from the $\mathcal{M}^\#$ class. Similarly to the 2-concatenation, in~\cite{KudinPPZ25_IEEE} the authors provided a complete description of the $\mathcal{M}$-subspaces of the 4-concatenation $h=f_1||f_2||f_3||f_4$ (where the $f_i$ are arbitrary). Using this characterization, the following result gives a necessary and sufficient condition for the membership of bent functions constructed via 4-concatenations in the $\mathcal{M}^\#$ class.

\begin{corollary}\label{cor: 4concatenationoutsideMM}\cite[Corollary 2]{KudinPPZ25_IEEE}
	Let $f=f_1 \vert \vert f_2 \vert \vert f_3 \vert \vert f_4\in\mathcal{B}_{n+2}$ be the concatenation of $f_1, \ldots ,f_4 \in \mathcal{B}_n$ and assume that $f$ is bent; thus $f_i$ are bent, semi-bent or five-valued spectra functions. Then, $f$ is outside of the $\cM^{\#}$ class if and only if the following conditions hold:
	\begin{enumerate}[a)]
		\item The functions $f_1, \dots ,f_4$ do not share a common $(n/2+1)$-dimensional $\cM$-subspace; 
		\item There are no common $(n/2)$-dimensional $\cM$-subspaces $V \subset \F_2^n$ of $f_1, \dots ,f_4$ such that there is an element $a\in \F_2^n$ for which
		\begin{equation*}
			\begin{gathered}
				D_vf_1+D_vf_2^a= D_vf_3+D_vf_4^a=0, \text{ for all } v \in V, \text{ or} \\
				D_vf_1+D_vf_3^a= D_vf_2+D_vf_4^a=0, \text{ for all } v \in V, \text{ or} \\
				D_vf_1+D_vf_4^a= D_vf_2+D_vf_3^a=0, \text{ for all } v \in V.\phantom{\text{ or}}
			\end{gathered}
		\end{equation*}
		\item There are no common $(n/2-1)$-dimensional $\cM$-subspaces $V \subset \F_2^n$ of $f_1, \dots ,f_4$ such that there are elements $a,b\in \F_2^n$ (not necessarily different), for which
		\begin{equation*}
			\begin{gathered}
				D_vf_1+D_vf_3^a= D_vf_2+D_vf_4^a=D_vf_1+D_vf_2^b =\\
				D_vf_3+D_vf_4^b=0, \text{ for all } v \in V, \text{ and} \\
				f_1(x)+f_2(x+b)+f_3(x+a)+f_4(x+a+b)=0,
			\end{gathered}
		\end{equation*}
		for all $x \in \F_2^n$.
	\end{enumerate}
\end{corollary}

An important application of the algebraic characterization of $\mathcal{M}$-subspaces is the following characterization of bent functions outside $\mathcal{M}^\#$ constructed via 4-concatenation using only two functions.

\begin{corollary}\cite[Corollary 4]{KudinPPZ25_IEEE}\label{cor: insideMMgh} Let $g,h \in\mathcal{B}_n$ be two arbitrary bent functions. Then, the function $f=g||h||g||(h+1)\in\mathcal{B}_{n+2}$ is bent and $f$ is in the $\cM^{\#}$ class if and only if the functions $g$ and $h$ have a common $(n/2)$-dimensional $\cM$-subspace, thus $g, h \in \cM^\#$.
\end{corollary}

With this result, one can easily construct bent functions outside $\mathcal{M}^\#$, as indicated in~\cite[Corollary~5]{KudinPPZ25_IEEE}, by taking a bent function $g \in \mathcal{B}_n$ outside $\cM^\#$, with $n \geq 8$, and an arbitrary bent function $h \in \mathcal{B}_n$. Then, the 4-concatenation $f=g||h||g||(h+1)\in \mathcal{B}_{n+2}$ is a bent function outside $\mathcal{M}^{\#}$. As a consequence of this result, one obtains the following bound on the number of bent functions outside $\mathcal{M}^\#$.

\begin{theorem}\cite[Theorem 5]{KudinPPZ25_IEEE}\label{theo:nmbbent}
	For $n \geq 6$, the number of bent functions outside $\mathcal{M}^\#$ in $n+2$ variables is always strictly greater than the number of all bent functions in $n$ variables.
\end{theorem}

Finally, we note that further investigation of secondary constructions of bent functions is essential, particularly regarding their linearity indexes and the generic description of their possible $\mathcal{M}$-subspaces. 

\begin{openproblem}
    Determine the $\mathcal{M}$-subspaces and the linearity index of secondary constructions of bent functions in terms of the properties of their building blocks. Provide necessary and sufficient conditions ensuring that the families of bent functions obtained via secondary constructions are outside the $\mathcal{M}^{\#}$ and $\mathcal{GMM}_{n/2+k}^\#$ classes, for $k>0$.
\end{openproblem}

\section{Conclusion and open problems}\label{sec Conclusion}

In this article, we summarized some of the recent advances related to the understanding of bent functions provably outside the completed Maiorana-McFarland class, specifically focusing on the $\mathcal{C}$ and $\mathcal{D}$ classes, bent functions constructed with the $P_\tau$ property, the generalized completed Maiorana-McFarland class, and the theory of $\mathcal{M}$-subspaces. To conclude the article, we would like to give a few more open problems in addition to those mentioned in the previous sections.

As indicated by numerous works considered in this article, the research on bent functions has been largely driven by two questions: first, to show that a given Boolean function is bent; second, to determine whether it belongs to some of the well-established classes. In most cases, the second question was asked for the Maiorana-McFarland class, thanks to its algebraic characterization, which is not available for most of the other classes. We believe that a natural next question is to describe the algebraic structure of newly constructed bent functions, which can be quantified by the linearity index. This question is particularly interesting for bent functions obtained via finite-field constructions, as well as for bent functions possessing additional cryptographically relevant properties. Since the Maiorana-McFarland construction is widely used for constructing various classes of cryptographically relevant functions, similar questions arise for related classes of functions beyond the bent function setting.

\begin{openproblem}
    For a given infinite family of (bent) functions on $\mathbb{F}_{2}^{n}$, determine the linearity index of each member and its $\mathcal{M}$-subspaces.
\end{openproblem}

In recent years, there has been an extended interest in the investigation of $p$-ary bent functions; for definitions and recent results, we refer to the survey articles \cite{Meidl22_CCDS,AKM_Survey}, as well as some recent works indicating that problems considered for classical Boolean bent functions are also of interest in the more general odd-$p$ setting~\cite{MesnagerOS18,KoelschP24_COMB}. Since the $\mathcal{C}$ and $\mathcal{D}$ classes can be defined when $p$ is odd~\cite{Nyberg93_CD_p_odd}, and since the characterization of the (generalized) Maiorana-McFarland class remains valid in that setting~\cite{Budaghyan14_BOOK}, it is natural to extend the binary theory to the $p$-odd case as well.

\begin{openproblem}
    1. Determine the intersection of the $\mathcal{C}$, $\mathcal{D}$, and $\mathcal{D}_0$ classes with the $\mathcal{M}^\#$ class for $p$-ary bent functions. 2. Extend the theory of $\mathcal{M}$-subspaces to the $p$-odd case.
\end{openproblem}

Boolean bent functions are indicators of Hadamard difference sets in elementary abelian 2-groups~\cite{Dillon74_PhD}. In the theory of difference sets, one of the main questions is to determine which groups contain difference sets~\cite{Applebaum_DS_2023}. We believe that answering the following question would help transfer some of the known constructions of bent functions (outside $\mathcal{M}^\#$) to the context of difference sets.

\begin{openproblem}
    Use methods for constructing bent functions outside $\mathcal{M}^\#$ to construct difference sets that are provably different (up to equivalence) from the McFarland ones~\cite{McFarland73_JCTA}.
\end{openproblem}

Finally, we emphasize that questions concerning bent functions within $\mathcal{M}^\#$ remain of interest, as some of them are still open. Although the cardinality of this class was recently estimated in~\cite{KolomeecB25}, the number of EA-equivalence classes is not known (except for the cases of six and eight variables, determined computationally in~\cite{Rothaus76_JCTA} and~\cite{LP_BFA_2024}).

\begin{openproblem}
1. Provide a good estimate of the number of EA-equivalence clas\-ses inside the (completed) Maiorana-McFarland class.
2. In general, how should permutations $\pi _{1},\pi _{2}$ and Boolean functions $h_{1},h_{2}\in \mathcal{B}_{m}$ be selected to ensure that the resulting Maiorana-McFarland bent functions, $f_{\phi _{1},h_{1}}$ and $f_{\phi _{2},h_{2}}$, are EA-in\-eq\-uiva\-lent on $\mathbb{F}_{2}^{m}\times \mathbb{F}_{2}^{m}$?
\end{openproblem}

\section*{Acknowledgments} 
Enes Pasalic and Sadmir Kudin are supported in part by the Slovenian Research Agency (research program P1-0404 and research project J1-60012). Fengrong Zhang is supported in part by the Natural Science Foundation of China (No. 62372346), and the Youth Innovation Team of Shaanxi Universities. Alexandr Polujan is supported by the Deutsche Forschungsgemeinschaft (DFG, German Research Foundation) --- project number 581467925.

The authors would like to thank the reviewers for their careful reading of the manuscript and for their valuable comments and suggestions, which have helped improve both the presentation and the clarity of the paper.


\end{document}